\documentclass[review,3p,11pt]{elsarticle}

\usepackage{epstopdf}
\usepackage{bookmark}
\usepackage{amsfonts}
\usepackage{lineno,hyperref}
\usepackage{float}
\usepackage[font=small,labelfont=bf,labelsep=none]{caption}
\usepackage{graphicx}
\usepackage{color}
\usepackage{subfigure}

\usepackage{amsmath,amsthm,mathrsfs}
\usepackage{amssymb}

\usepackage{geometry}

\allowdisplaybreaks

\newtheorem{theorem}{Theorem}
\newtheorem{lemma}{Lemma}
\newtheorem{corollary}{Corollary}
\newtheorem{remark}{Remark}

\begin{document}

\captionsetup[figure]{labelfont={bf},labelformat={default},labelsep=period,name={Fig.}}
\captionsetup[table]{labelsep=newline, singlelinecheck=false}

\begin{frontmatter}

\title{{Structure-Preserving Oscillation-Eliminating Discontinuous Galerkin Schemes for Ideal MHD Equations: Locally Divergence-Free and Positivity-Preserving}}
\author[mymainaddress]{Mengqing Liu}
\ead{liumq@sustech.edu.cn}

\author[mymainaddress,mysecondaryaddress,mythirdaddress]{Kailiang Wu\corref{mycorrespondingauthor}}
\cortext[mycorrespondingauthor]{Corresponding author}
\ead{wukl@sustech.edu.cn}

\address[mymainaddress]{Department of Mathematics and Shenzhen International Center for Mathematics, Southern University of Science and Technology, Shenzhen 518055, China}
\address[mysecondaryaddress]{National Center for Applied Mathematics Shenzhen (NCAMS), Shenzhen 518055, China}
\address[mythirdaddress]{Guangdong Provincial Key Laboratory of Computational Science and Material Design, Shenzhen 518055, China}

\date{April 27, 2024}

\begin{abstract}
Numerically simulating magnetohydrodynamics (MHD) poses notable challenges, including the suppression of spurious oscillations near discontinuities (e.g., shocks) and preservation of essential physical structures (e.g., the divergence-free constraint of  magnetic field and the positivity of density and pressure). This paper develops  structure-preserving oscillation-eliminating discontinuous Galerkin (OEDG) schemes for ideal MHD, as a sequel to our recent work [M.~Peng,~Z. Sun \& K.~Wu, arXiv:2310.04807, 2023]. The schemes leverage a locally divergence-free (LDF) oscillation-eliminating (OE) procedure to suppress spurious oscillations while retaining the LDF property of magnetic field and many desirable attributes of original DG schemes, such as conservation, local compactness, and optimal convergence rates. The OE procedure is based on the solution operator of a novel damping equation, a linear system of ordinary differential equations that are exactly solvable without any discretization. The OE procedure is performed after each Runge--Kutta stage and does not impact DG spatial discretization, facilitating its easy integration into existing DG codes as an independent module. Moreover, this paper presents a rigorous positivity-preserving (PP) analysis of the LDF OEDG schemes on Cartesian meshes, utilizing the optimal convex decomposition technique [S.~Cui, S.~Ding \& K.~Wu, {\em SIAM J.~Numer.~Anal.}, 62:775--810, 2024] and the geometric quasi-linearization (GQL) approach [K.~Wu \& C.-W.~Shu, {\em SIAM Review}, 65:1031--1073, 2023]. Efficient PP LDF OEDG schemes are derived by incorporating appropriate discretization of Godunov--Powell source terms into only the discrete equations of cell averages, under a condition achievable through a simple PP limiter. Several one- and two-dimensional MHD tests verify the accuracy, effectiveness, and robustness of the proposed structure-preserving OEDG schemes.

\end{abstract}

\begin{keyword}
magnetohydrodynamics (MHD)\sep oscillation-eliminating discontinuous Galerkin (OEDG) \sep locally divergence-free (LDF)\sep positivity-preserving (PP) \sep multidimensional
\end{keyword}

\end{frontmatter}

\renewcommand\baselinestretch{2}

\section{Introduction}\label{sec:Introd}

The magnetohydrodynamic (MHD) equations are extensively employed to model macro-scale behaviors across various domains, such as laboratory, space, and astrophysical plasmas. These equations find notable applications in the study of phenomena including coronal mass ejections \cite{Yang_2023}, the evolution of solar coronal magnetic fields \cite{Jiang_2021}, and the dynamics within Earth's magnetosphere \cite{GUO2016543}. As a nonlinear hyperbolic system, the ideal compressible MHD equations are inherently complex, making analytical solutions challenging to derive. Thus, numerical methods have become a primary approach in this research area. Among these methods, discontinuous Galerkin (DG) methods are especially preferred due to their ability to achieve high-order accuracy, excellent local compactness, and adaptability to arbitrary mesh configurations (cf.~\cite{COCKBURN1998199,COCKBURN2004588}). The effectiveness of DG methods in MHD numerical simulations is well-documented in various studies (e.g., \cite{Li2005,LI20114828,zhao2017runge,Mocz2013,Guillet2019,Balsara2021}). However, in scenarios with strong discontinuities like shocks, conventional high-order DG methods may generate spurious oscillations near these discontinuities, leading to non-physical results and potential simulation failures.

To suppress spurious oscillations in the DG methods for hyperbolic conservation laws, various techniques have been developed, including but not limited to the total variation bounded limiter \cite{COCKBURN1998199}, the weighted essentially non-oscillatory limiter \cite{Qiu2005}, and artificial viscosity techniques  (e.g.,~\cite{ZINGAN2013479,Huang2020An}). These strategies effectively control oscillations and enhance the stability of DG schemes. However, some of them may compromise certain desirable attributes of the original DG schemes, such as optimal convergence rates, or depend on problem-specific parameters.
Recently, Lu, Liu, and Shu \cite{Jianfang2021} introduced the so-called oscillation-free DG (OFDG) method for scalar conservation laws, by ingeniously incorporating a damping term into semi-discrete DG schemes. The OFDG method can effectively mitigate numerical oscillations near strong discontinuities while retaining many beneficial features of the original DG schemes, including conservation, superconvergence, and optimal error estimates.
Subsequently, the OFDG method was successfully extended to hyperbolic systems \cite{Liu2022} with characteristic decomposition and applied to shallow water equations \cite{Liu2022shallow}, nonlinear degenerate parabolic equations \cite{Tao2023para}, and multi-component chemically reacting flows \cite{Du2023multi}. However, the characteristic decomposition in the OFDG method may be computationally expensive and require careful adjustments to the eigenvectors due to a lack of scale invariance. Furthermore, as Liu, Lu, and Shu mentioned \cite{Liu2022}, the stiffness of the large damping terms near strong discontinuities makes the Courant--Friedrichs--Lewy (CFL) condition more restrictive, necessitating the use of modified exponential Runge–Kutta methods to alleviate this restriction.
More recently, motivated by \cite{Jianfang2021}, a novel oscillation-eliminating DG (OEDG) method was systematically developed in \cite{peng2023oedg}. This method alternately evolves the conventional semi-discrete DG scheme and a novel damping equation. Unlike the OFDG method, the OEDG approach decouples the damping terms from the DG formulations, facilitating non-intrusive integration into existing DG codes, simplifying implementation, and enhancing efficiency. Additionally, the OEDG method does not require characteristic decomposition for hyperbolic systems, maintains stability under the standard CFL condition without the need for exponential Runge--Kutta discretization, and exhibits desirable properties such as scale invariance and evolution invariance across various scales and wave speeds.

In ideal MHD, the magnetic field must be divergence-free, and both density and internal energy (or pressure) must remain positive. Preserving these physical structures---zero divergence of the magnetic field and positivity of density and pressure---presents significant challenges in MHD simulations. Large divergence errors can introduce nonphysical features or cause numerical instabilities \cite{Brackbill1980The,Balsara1999A}. To mitigate these divergence errors, several techniques have been developed, including but not limited to the projection method \cite{Brackbill1980The}, the constrained transport method \cite{Evans1988Simulation,GARDINER2005509}, the eight-wave method  \cite{Powell1999Regular}, the hyperbolic divergence cleaning method \cite{Dedner2002Hyperbolic}, the locally divergence-free (LDF) method \cite{Li2005,YAKOVLEV201380,LIU2021110694,Ding2024A}, and the globally divergence-free (GDF) method \cite{BALSARA2015687,LI20114828,XU2016203,Fu2018globally,Balsara2021}. These methods have been widely applied across various schemes due to their effectiveness in controlling divergence errors. Positivity of density and pressure can be compromised in challenging MHD simulations involving low density or pressure, high Mach numbers, large magnetic energies, or strong discontinuities. To address this, early research primarily focused on developing robust one-dimensional approximate Riemann solvers; see, for example,  \cite{BALSARA1999Maintaining,JANHUNEN2000A,Bouchut2007A,Bouchut2010A}.
Waagan \cite{WAAGAN20098609} introduced a second-order MUSCL-Hancock scheme using relaxation Riemann solvers \cite{Bouchut2007A,Bouchut2010A} that ensures the positivity of density and pressure; the robustness of this scheme was further validated through extensive testing \cite{WAAGAN20113331}. There has been notable progress in the development of positivity-preserving (PP) limiters for higher-order schemes \cite{ZHANG20103091,ZHANG20108918}.
Balsara \cite{BALSARA2012Self} proposed a self-adjusting PP limiter to maintain the positivity of reconstructed density and pressure in ideal MHD. Cheng et al.~\cite{CHENG2013255} extended the PP limiter of Zhang--Shu \cite{ZHANG20108918} to ideal MHD and demonstrated the PP property of one-dimensional (central) DG schemes under certain assumptions. Christlieb et al.~\cite{Christlieb2015PP,CHRISTLIEB2016218} developed high-order PP finite difference methods for ideal MHD, via a parameterized flux limiter \cite{XIONG2013310}. These PP techniques have significantly enhanced the robustness of MHD schemes, as evidenced by numerous numerical tests. However, the inherent complexity of MHD equations and the unclear relationship between the PP property and the divergence-free condition have limited earlier research in proposing numerical schemes that can be rigorously proven to preserve positivity for ideal MHD.

In recent years, important advances \cite{Wu2018Positivity,Wu2018A,Wu2019Provably,Wu2023Geometric,Wu2023provab} have been made in revealing the theoretical connections between the PP property and the divergence-free condition for MHD systems.  For the ideal MHD equations, the work \cite{Wu2018Positivity} first rigorously proved that a discrete divergence-free condition is intrinsically linked to achieving the PP property. As discovered in \cite{Wu2018Positivity}, even slight deviations from the discrete divergence-free condition could compromise the PP property. Furthermore, it has been found in  \cite{Wu2018A} that even the exact smooth solutions of the conservative MHD system with nonzero magnetic divergence may result in negative pressure.
Fortunately, as shown in \cite{Wu2019Provably}, this issue does not arise in the symmetrizable MHD system, which incorporates the Godunov--Powell source term \cite{Godunov1972SymmetricFO,Powell1995AnUS}. A provably PP DG scheme for this system was first developed in \cite{Wu2018A}. Subsequently, \cite{Wu2019Provably} established a framework for provably positive high-order finite volume and DG schemes on general meshes.
These advancements leverage an innovative approach termed geometric quasilinearization (GQL) \cite{Wu2018Positivity,Wu2023Geometric}, which equivalently transforms complex nonlinear constraints---such as ensuring the positivity of internal energy or pressure in the MHD system---into manageable linear constraints through the integration of appropriate auxiliary variables. Capitalizing on the geometric characteristics of convex sets, a comprehensive GQL framework was recently proposed in \cite{Wu2023Geometric} to address bound-preserving schemes under nonlinear constraints.
 For developments in the field of relativistic MHD, refer to \cite{Wu2017RMHD,Wu2021RMHD,ding2024gqlbased}.
More recently, a novel second-order structure-preserving finite element method was introduced in \cite{dao2024structure}, based on convex limiting and a novel splitting technique.

This paper aims to develop the LDF OEDG method for ideal MHD and then present a rigorous PP analysis. Our analysis establishes the probable PP property for the LDF OEDG method through the optimal convex decomposition technique \cite{CUI2023111882,Cui2024On} and the GQL approach \cite{Wu2023Geometric}. The contributions, novelty, and significance of this work are outlined as follows: 
\begin{itemize}
  \item We design an LDF OE procedure, which effectively eliminates spurious oscillations while locally maintaining the divergence-free magnetic field. This procedure is integrated as an independent module within the LDF DG framework, applied after each Runge--Kutta stage, resulting in the LDF OEDG method. The procedure is built on the evolution operator of a novel damping equation, whose exact solution is explicitly formulated without any discretization. It serves as an LDF filter to the DG modal coefficients and is easily implemented and efficient, involving merely the multiplication of modal coefficients by scalars. The LDF OEDG method retains the desirable features of original DG schemes, such as conservation, local compactness, and optimal convergence rates, and is naturally extensible to general meshes.
  \item A rigorous PP analysis of the LDF OEDG schemes with the Harten--Lax--van Leer (HLL) flux is carried out on Cartesian meshes. This analysis employs the recently developed GQL approach \cite{Wu2018Positivity,Wu2023Geometric}, which skillfully transforms the nonlinear constraint of pressure into a family of equivalent linear constraints by introducing several auxiliary variables. To ensure the weak PP property for the cell averages of the OEDG solution, we incorporate suitable upwind discretization of Godunov--Powell source terms into the schemes. This addresses the impact of divergence errors on the PP property. Unlike our previous work \cite{Wu2018A,Wu2019Provably}, we establish a couple of ``two-state" inequalities for estimating the effect of the HLL flux on PP on rectangular cells, which simplifies/improves the estimates in \cite{Wu2019Provably} for general meshes when reduced to Cartesian meshes. Another difference is that we incorporate the discrete Godunov--Powell source terms solely into the evolution equations of the cell averages, as the weak PP property relates only to these discrete equations, thereby enhancing simplicity and efficiency.
  \item A key aspect of the PP analysis is the convex decomposition, where the cell averages of the numerical solution are decomposed into a convex combination of its point values at certain quadrature points. We provide a comprehensive PP analysis via general convex decomposition \cite{Cui2024On} on Cartesian meshes and derive PP CFL conditions through specific convex decompositions, including the Zhang--Shu convex decomposition \cite{ZHANG20103091} and the optimal convex decomposition \cite{CUI2023111882,Cui2024On}. The optimal convex decomposition leads to a milder theoretical PP CFL condition and requires fewer nodes, thus improving efficiency further. Our theoretical analyses establish the weak PP property of the LDF OEDG method, implying that employing a simple scaling PP limiter can enforce the pointwise PP property of the OEDG solution at any point of interest.
\end{itemize}
Some one- and two-dimensional (2D) MHD examples are utilized to test the accuracy, effectiveness, and robustness of the proposed PP LDF OEDG method.
The numerical results demonstrate that our method achieves the expected optimal convergence order, effectively suppresses oscillations, and performs robustly in challenging scenarios.

The paper is organized as follows.
Section \ref{sec:2Dgovequcon} presents the LDF OEDG method for 2D ideal MHD equations.
In section \ref{sec:PPLDFOEDGMHD}, we construct a PP LDF OEDG scheme for 2D ideal MHD equations and prove the PP property.
Section \ref{Sec:numericaltest} provides several numerical examples for MHD in one and two dimensions.
Section \ref{Sec:conclude} gives the concluding remarks of this paper.

\section{Locally divergence-free OEDG method for ideal MHD}\label{sec:2Dgovequcon}

This section presents the LDF OEDG method for the ideal MHD system. Without loss of generality, we focus on the 2D MHD equations, while our method is directly applicable to the one-dimensional MHD system and can be readily extended to the three-dimensional case.

\subsection{Governing equations of 2D MHD system}

The 2D conservative MHD system can be written as
\begin{equation}\label{Equ:govequ2D}
	\frac{\partial \boldsymbol{U}}{\partial t}+\frac{\partial\boldsymbol{F}_{1}\left(\boldsymbol{U}\right)}{\partial{x}}
	+\frac{\partial\boldsymbol{F}_{2}\left(\boldsymbol{U}\right)}{\partial{y}}= \boldsymbol{0},
\end{equation}
where
$\boldsymbol{U}=\left(\rho, \rho\boldsymbol{u}, \boldsymbol{B}, E\right)^{\top}=\left(\rho, \boldsymbol{m}, \boldsymbol{B}, E\right)^{\top}$
is the conservation vector.  In alignment with existing studies  
in the literature, this paper considers the general configurations where $\boldsymbol{m}, \boldsymbol{B} \in \mathbb R^3$ rather than $\mathbb R^2$, to accommodate the general 2D case that is a reduction from the 3D scenario. 
The fluxes in the $x$- and $y$-directions are respectively given by
\begin{equation*}\label{Equ:f12D}
	\boldsymbol{F}_{1}\left(\boldsymbol{U}\right)=
	\begin{pmatrix}
		\rho u_{1} \\
		\rho {u_{1}}^{2}+\left(p+{\|\boldsymbol{B}\|}^{2}/2\right)-
		B_{1}^{2} \\
		\rho u_{1} u_{2}-B_{1}B_{2} \\
		\rho u_{1} u_{3}-B_{1}B_{3} \\
		0 \\
		u_{1}B_{2}-B_{1}u_{2} \\
		u_{1}B_{3}-B_{1}u_{3} \\
		\left(E+p+{\|\boldsymbol{B}\|}^{2}/2\right)u_{1}-
		B_{1}\left(\boldsymbol{B}\cdot\boldsymbol{u}\right)\\
	\end{pmatrix},\quad
	\boldsymbol{F}_{2}\left(\boldsymbol{U}\right)=
	\begin{pmatrix}
		\rho u_{2} \\
		\rho u_{2} u_{1}-B_{2}B_{1} \\
		\rho {u_{2}}^{2}+\left(p+{\|\boldsymbol{B}\|}^{2}/2\right)-
		B_{2}^{2} \\
		\rho u_{2} u_{3}-B_{2}B_{3} \\
		u_{2}B_{1}-B_{2}u_{1} \\
		0 \\
		u_{2}B_{3}-B_{2}u_{3} \\
		\left(E+p+{\|\boldsymbol{B}\|}^{2}/2\right)u_{2}-
		B_{2}\left(\boldsymbol{B}\cdot\boldsymbol{u}\right)\\
	\end{pmatrix}.
\end{equation*}
Here, $\rho$ represents the fluid density, $\boldsymbol{m}$ is the momentum vector, and $\boldsymbol{B} = (B_{1}, B_{2}, B_{3})$ is the magnetic field vector. 
The variable $p$ stands for thermal pressure, and $E$ is the total energy density, expressed by the ideal state equation as
\begin{equation*}
	E = \frac{p}{\gamma-1} + \frac{1}{2\rho}\|\boldsymbol{m}\|^{2} + \frac{1}{2}
	\|\boldsymbol{B}\|^{2},
\end{equation*}
where $\gamma$ is the adiabatic index and $\| \cdot \|$ denotes the vector-2 norm.  Additionally, the physical solutions of the MHD equations require a divergence-free constraint on the magnetic field:
\begin{equation}\label{eq:DF}
	\nabla \cdot \boldsymbol{B} : = \frac{\partial B_1}{\partial x} + \frac{\partial B_2} {\partial y} = 0,
\end{equation}
which implies the absence of magnetic monopoles.

Divide the conservation variables $\boldsymbol{U}$ into two parts:
\begin{equation*}\label{Equ:2DRandQ}
	\boldsymbol{R}=\left(\rho, \rho u_{1}, \rho u_{2}, \rho u_{3}, B_{3}, E\right)^{\top}
	\quad \text{and} \quad
	\boldsymbol{Q}=\left(B_{1}, B_{2}\right)^{\top},
\end{equation*}
and correspondingly divide the flux $\boldsymbol{F}_{i}(\boldsymbol{U})$ into $\boldsymbol{F}_{i}^{R}(\boldsymbol{U})$ and $\boldsymbol{F}_{i}^{Q}(\boldsymbol{U})$. The MHD system \eqref{Equ:govequ2D} is then split into
\begin{equation}\label{Equ:2DgovequR}
	\frac{\partial \boldsymbol{R}}{\partial t}+\frac{\partial\boldsymbol{F}_{1}^{R}\left(\boldsymbol{U}\right)}{\partial{x}}
	+\frac{\partial\boldsymbol{F}_{2}^{R}\left(\boldsymbol{U}\right)}{\partial{y}}= \boldsymbol{0},
\end{equation}
\begin{equation}\label{Equ:2DgovequQ}
	\frac{\partial \boldsymbol{Q}}{\partial t}+\frac{\partial\boldsymbol{F}_{1}^{Q}\left(\boldsymbol{U}\right)}{\partial{x}}
	+\frac{\partial\boldsymbol{F}_{2}^{Q}\left(\boldsymbol{U}\right)}{\partial{y}}= \boldsymbol{0}.
\end{equation}
The objective of the LDF OEDG method is to find an approximate solution $\boldsymbol{U}_h$, which comprises $\boldsymbol{R}_h$ and $\boldsymbol{Q}_h$  obtained through different OEDG discretizations.

Next, we will elaborate on our LDF OEDG method, including the LDF finite element spaces, the semi-discrete LDF DG formulation, and the fully discrete LDF OEDG discretization. It incorporates an LDF OE procedure after each Runge--Kutta stage to eliminate spurious oscillations. The PP property of the LDF OEDG method will be systematically addressed in section \ref{sec:PPLDFOEDGMHD}.


\subsection{LDF DG finite element spaces}

Let $\Omega=\bigcup_{ij}I_{ij}$ with mesh element $I_{ij}=\left[x_{i-\frac{1}{2}}, x_{i+\frac{1}{2}}\right]\times\left[y_{j-\frac{1}{2}}, y_{j+\frac{1}{2}}\right]$   be a partition of the 2D spatial domain $\Omega$.
Denote by $\Delta x$ and $\Delta y$ the spatial step-sizes in the $x$- and $y$-directions, respectively.
We now introduce two finite element DG spaces in which we seek approximations to $\boldsymbol{R}$ and $\boldsymbol{Q}$, respectively.

The approximation space for $\boldsymbol{R}$
is defined as follows:
\begin{equation*}\label{Equ:2DspaceV}
	\mathbb{V}_{R}^k:=\left\{\boldsymbol{v}\left(\boldsymbol{x}\right) \in L^2\left(\Omega\right):\left.\boldsymbol{v}\left(\boldsymbol{x}\right)\right|_{I_{ij}} \in \left[\mathbb{P}^k\left(I_{ij}\right)\right]^{6}~~ \forall i,j\right\},
\end{equation*}
where $\boldsymbol{x}=\left(x, y\right)$ denotes the spatial coordinate vector, and $\mathbb{P}^k(I_{ij})$ is the space of polynomials of total degree less than or equal to $k$ on the element $I_{ij}$.
An orthogonal basis of $\mathbb{P}^k(I_{ij})$ is denoted as $\{\phi_{ij}^{(\boldsymbol{\alpha})}(\boldsymbol{x}): |\boldsymbol{\alpha}| \leq k \}$, where $\boldsymbol{\alpha} = (\alpha_{1}, \alpha_{2})$ is the multi-index vector with $|\boldsymbol{\alpha}| = \alpha_1 + \alpha_2$.
For example, we adopt the following local orthogonal Legendre basis over cell $I_{ij}$:
\begin{equation*}\label{Equ:2DLegenbasis}
	1, \quad \xi, \quad \eta, \quad
	\xi^2-\frac{1}{3}, \quad \xi \eta, \quad \eta^2-\frac{1}{3}, \quad \cdots
\end{equation*}
where $\xi=2\left(x-x_{i}\right)/\Delta x$ and $\eta=2\left(y-y_{j}\right)/\Delta y$;
the corresponding multi-index vectors $\boldsymbol{\alpha}$ are
\begin{equation*}\label{Equ:2Dindex}
	\left(0, 0\right), \quad \left(1, 0\right), \quad \left(0, 1\right),
	\quad \left(2, 0\right), \quad \left(1, 1\right), \quad \left(0, 2\right), \quad \dots
\end{equation*}
The dimensionality of $\mathbb{P}^k(I_{ij})$ is $D_R^k = (k+1)(k+2)/2$. Using the Legendre basis, the approximate solution $\boldsymbol{R}_h(\boldsymbol{x}, t)$  on element $I_{ij}$ can be expressed as
\begin{equation}\label{Equ:2DsolutionR}
	\boldsymbol{R}_h\left(\boldsymbol{x}, t\right)=\sum_{\mu=0}^{k}\sum_{\left|
		\boldsymbol{\alpha}\right|=\mu}\boldsymbol{R}_{ij}^{\left(\boldsymbol{\alpha}\right)}
	\left(t\right)\phi_{ij}^{\left(\boldsymbol{\alpha}\right)}\left(\boldsymbol{x}\right) \quad \text{for} \quad \boldsymbol{x} \in I_{ij},
\end{equation}
where $\boldsymbol{R}_{ij}^{(\boldsymbol{\alpha})}(t)$ represents a moment of the approximate solution $\boldsymbol{R}_h(\boldsymbol{x}, t)$.

We seek a divergence-free approximation to $\boldsymbol{Q}$
in the following LDF space \cite{Li2005}:
\begin{equation*}\label{Equ:2DdivspaceW}
	\mathbb{V}_{Q}^k:=\left\{\boldsymbol{v}\left(\boldsymbol{x}\right)=
	\left(v_{1}\left(\boldsymbol{x}\right), v_{2}\left(\boldsymbol{x}\right)\right)^{\top}:
	v_{1}\left(\boldsymbol{x}\right), v_{2}\left(\boldsymbol{x}\right) \in \mathbb{P}^k\left(I_{ij}\right),
	\left.\nabla \cdot \boldsymbol{v}\left(\boldsymbol{x}\right)\right|_{I_{ij}}=0 ~~ \forall i,j
	\right\}.
\end{equation*}
The local orthogonal basis of $\mathbb{V}_{Q}^k$ on element $I_{ij}$ is denoted as $\{\boldsymbol{\varphi}_{ij}^{\left(\mu\right)}\left(\boldsymbol{x}\right),
\mu=0, 1, \cdots, D_{Q}^{k}-1\}$,
where $D_{Q}^{k}=\left(k+1\right)\left(k+4\right)/2$ is the dimensionality of the LDF polynomial space defined on $I_{ij}$.
We employ the following local orthogonal LDF basis
\begin{equation*}\label{Equ:2Ddivbasis}
	\begin{gathered}
		\left(\begin{array}{c}
			0 \\
			1 \\
		\end{array}\right),
		\quad
		\left(\begin{array}{c}
			1 \\
			0 \\
		\end{array}\right),
		\quad
		\left(\begin{array}{c}
			0 \\
			\xi \\
		\end{array}\right),
		\quad
		\left(\begin{array}{c}
			\eta \\
			0 \\
		\end{array}\right),
		\quad
		\left(\begin{array}{c}
			\Delta x\xi \\
			-\Delta y\eta \\
		\end{array}\right),
		\\
		\left(\begin{array}{c}
			\eta^{2}-1/3 \\
			0 \\
		\end{array}\right),
		\quad
		\left(\begin{array}{c}
			0 \\
			\xi^{2}-1/3 \\
		\end{array}\right),
		\quad
		\left(\begin{array}{c}
			\Delta x\left(\xi^{2}-1/3\right) \\
			-2\Delta y\xi\eta \\
		\end{array}\right),
		\quad
		\left(\begin{array}{c}
			-2\Delta x\xi\eta \\
			\Delta y\left(\eta^{2}-1/3\right) \\
		\end{array}\right),
		\quad
		\cdots,
	\end{gathered}
\end{equation*}
where $\xi=2\left(x-x_{i}\right)/\Delta x$ and $\eta=2\left(y-y_{j}\right)/\Delta y$.
With this basis, we express the approximate solution $\boldsymbol{Q}_h\left(\boldsymbol{x}, t\right)$ over element $I_{ij}$ as
\begin{equation}\label{Equ:2DapprosoluQ}
	\boldsymbol{Q}_h\left(\boldsymbol{x}, t\right)=\sum_{\mu=0}^{k}\sum_{\eta=D_{Q}^{\mu-1}}^{D_{Q}^{\mu}-1}
	Q_{ij}^{\left(\eta\right)}\left(t\right)\boldsymbol{\varphi}_{ij}^{\left(
		\eta\right)}\left(\boldsymbol{x}\right) \quad \text{for} \quad \boldsymbol{x} \in I_{ij}.
\end{equation}
where ${Q}_{ij}^{\left(\eta\right)}
\left(t\right)$ represents a moment of approximate polynomial solution $\boldsymbol{Q}_h\left(\boldsymbol{x}, t\right)$.
Note that the approximate magnetic field $\boldsymbol{Q}_{h}$ given in the above form naturally satisfies the divergence-free constraint \eqref{eq:DF} within each element $I_{ij}$.

\subsection{Semi-discrete LDF DG formulation}

The LDF DG method uses different spatial discretizations for $\boldsymbol{R}$ and $\boldsymbol{Q}$.

For any test function $v\left(\boldsymbol{x}\right) \in \mathbb{P}^k\left(I_{ij}\right)$, the approximate solution $\boldsymbol{R}_h\left(\boldsymbol{x}, t\right)$ is expected to satisfy
\begin{equation}\label{2DRapprosolu}
	\begin{aligned}
		&\frac{\mathrm{d}}{\mathrm{d} t} \int_{I_{ij}} \boldsymbol{R}_{h}\left(\boldsymbol{x}, t\right) v\left(\boldsymbol{x}\right)\mathrm{~d}\boldsymbol{x}
		=\int_{I_{ij}} \boldsymbol{F}_{1}^{R}\left(\boldsymbol{U}_{h}\left(\boldsymbol{x}, t\right)\right) \frac{\partial v\left(\boldsymbol{x}\right)}{\partial x} \mathrm{~d} \boldsymbol{x}+\int_{I_{ij}} \boldsymbol{F}_{2}^{R}\left(\boldsymbol{U}_{h}\left(\boldsymbol{x}, t\right)\right) \frac{\partial v\left(\boldsymbol{x}\right)}{\partial y} \mathrm{~d} \boldsymbol{x}\\
		&-\int_{y_{j-\frac{1}{2}}}^{y_{j+\frac{1}{2}}}\left(
		v\left(x_{i+\frac{1}{2}}^{-}, y\right)\widehat{\boldsymbol{F}_{1}^{R}}_{, i+\frac{1}{2}}
		\left(y, t\right)
		-v\left(x_{i-\frac{1}{2}}^{+}, y\right)\widehat{\boldsymbol{F}_{1}^{R}}_{, i-\frac{1}{2}}
		\left(y, t\right)
		\right)\mathrm{~d}y \\
		&-\int_{x_{i-\frac{1}{2}}}^{x_{i+\frac{1}{2}}}\left(
		v\left(x, y_{j+\frac{1}{2}}^{-}\right)\widehat{\boldsymbol{F}_{2}^{R}}_{, j+\frac{1}{2}}
		\left(x, t\right)
		-v\left(x, y_{j-\frac{1}{2}}^{+}\right)\widehat{\boldsymbol{F}_{2}^{R}}_{, j-\frac{1}{2}}
		\left(x, t\right)
		\right)\mathrm{~d}x,
	\end{aligned}
\end{equation}
where the superscripts ``$-$'' and ``$+$'' indicate the left- and right-hand-side limits at the cell interfaces, respectively, and
\begin{equation*}\label{Equ:2DF1equ}
	\widehat{\boldsymbol{F}_{1}^{R}}_{, i+\frac{1}{2}}\left(y, t\right):=
	\widehat{\boldsymbol{F}_{1}^{R}}\left(
	\boldsymbol{U}_{h}\left(x_{i+\frac{1}{2}}^{-}, y, t\right),
	\boldsymbol{U}_{h}\left(x_{i+\frac{1}{2}}^{+}, y, t\right)\right), \\
\end{equation*}
\begin{equation*}\label{Equ:2DF2equ}
	\widehat{\boldsymbol{F}_{2}^{R}}_{, j+\frac{1}{2}}\left(x, t\right):=
	\widehat{\boldsymbol{F}_{2}^{R}}\left(
	\boldsymbol{U}_{h}\left(x, y_{j+\frac{1}{2}}^{-}, t\right),
	\boldsymbol{U}_{h}\left(x, y_{j+\frac{1}{2}}^{+}, t\right)\right). \\
\end{equation*}
Here, $\widehat{\boldsymbol{F}_{1}^{R}}$ and $\widehat{\boldsymbol{F}_{2}^{R}}$ denote the numerical fluxes in the $x$- and $y$-directions, respectively.
In this work, we use the PP HLL fluxes
\begin{equation}\label{Equ:2DHLLfluxR}
	\widehat{\boldsymbol{F}_{\ell}^{R}}\left(\boldsymbol{U}^{-},\boldsymbol{U}^{+}\right)=
	\frac{\mathscr{V}_{\ell}^{+}\boldsymbol{F}_{\ell}^{R}\left(\boldsymbol{U}^{-}\right)-
		\mathscr{V}_{\ell}^{-}\boldsymbol{F}_{\ell}^{R}\left(\boldsymbol{U}^{+}\right)+
		\mathscr{V}_{\ell}^{-}\mathscr{V}_{\ell}^{+}\left(\boldsymbol{R}^{+}-\boldsymbol{R}^{-}\right)}
	{\mathscr{V}_{\ell}^{+}-\mathscr{V}_{\ell}^{-}}, \quad \ell=1, 2,
\end{equation}
where $\mathscr{V}_{\ell}^{-}$ and $\mathscr{V}_{\ell}^{+}$ denote the (properly estimated) non-positive minimum and non-negative maximum wave speeds, respectively,
with $\ell=1$ for the $x$-direction and $\ell=2$ for the $y$-direction.
The estimates of $\mathscr{V}_{\ell}^{-}$ and $\mathscr{V}_{\ell}^{+}$ are crucial for the PP property; see section \ref{Sec:PPHLLflux} for details.
Taking the test function $v\left(\boldsymbol{x}\right)$ as the (orthogonal) Legendre basis function $\phi_{ij}^{\left(\boldsymbol{\alpha}\right)}\left(\boldsymbol{x}\right)$ and using the Gauss quadrature of sufficiently high accuracy to approximate the right-hand side of (\ref{2DRapprosolu}), we obtain the following semi-discrete scheme
\begin{equation}\label{Equ:2DsemidiscreR}
	\begin{aligned}
		&\left(\int_{I_{ij}} \left(\phi_{ij}^{\left(\boldsymbol{\alpha}\right)}\left(\boldsymbol{x}\right)
		\right)^{2} \mathrm{d} \boldsymbol{x}\right)
		\frac{\mathrm{d}\boldsymbol{R}_{ij}^{\left(\boldsymbol{\alpha}\right)}
			\left(t\right)}{\mathrm{d} t}
		=\Delta x \Delta y\sum_{\mu=1}^{q}\sum_{\eta=1}^{q}\omega_{\mu}^{G}\omega_{\eta}^{G} \boldsymbol{F}_{1}^{R}\left(\boldsymbol{U}_{h}\left(x_{i}^{(\mu)}, y_{j}^{(\eta)}, t\right)\right) \partial_{x}\phi_{ij}^{\left(\boldsymbol{\alpha}\right)}
		\left(x_{i}^{(\mu)}, y_{j}^{(\eta)}\right)\\
		&+\Delta x \Delta y\sum_{\mu=1}^{q}\sum_{\eta=1}^{q}\omega_{\mu}^{G}\omega_{\eta}^{G} \boldsymbol{F}_{2}^{R}\left(\boldsymbol{U}_{h}\left(x_{i}^{(\mu)}, y_{j}^{(\eta)}, t\right)\right) \partial_{y}\phi_{ij}^{\left(\boldsymbol{\alpha}\right)}
		\left(x_{i}^{(\mu)}, y_{j}^{(\eta)}\right)
		\\
		&-\Delta y \sum_{\mu=1}^{q}\omega_{\mu}^{G}\left(
		\phi_{ij}^{\left(\boldsymbol{\alpha}\right)}\left(x_{i+\frac{1}{2}}^{-}, y_{j}^{(\mu)}\right)\widehat{\boldsymbol{F}_{1}^{R}}_{, i+\frac{1}{2}}
		\left(y_{j}^{(\mu)}, t\right)
		-\phi_{ij}^{\left(\boldsymbol{\alpha}\right)}\left(x_{i-\frac{1}{2}}^{+}, y_{j}^{(\mu)}\right)\widehat{\boldsymbol{F}_{1}^{R}}_{, i-\frac{1}{2}}
		\left(y_{j}^{(\mu)}, t\right)
		\right)\\
		&-\Delta x\sum_{\mu=1}^{q}\omega_{\mu}^{G}\left(
		\phi_{ij}^{\left(\boldsymbol{\alpha}\right)}\left(x_{i}^{(\mu)}, y_{j+\frac{1}{2}}^{-}\right)\widehat{\boldsymbol{F}_{2}^{R}}_{, j+\frac{1}{2}}
		\left(x_{i}^{(\mu)}, t\right)
		-\phi_{ij}^{\left(\boldsymbol{\alpha}\right)}\left(x_{i}^{(\mu)}, y_{j-\frac{1}{2}}^{+}\right)\widehat{\boldsymbol{F}_{2}^{R}}_{, j-\frac{1}{2}}
		\left(x_{i}^{(\mu)}, t\right)
		\right), \quad \forall
		|\boldsymbol{\alpha}| \le k,
	\end{aligned}
\end{equation}
where $\left\{x_{i}^{(\mu)}\right\}_{\mu=1}^{q}$ and $\left\{y_{j}^{(\mu)}\right\}_{\mu=1}^{q}$
are the $q$-point Gauss quadrature nodes with $2q-1\geq2k$ in $\left[x_{i-\frac{1}{2}}, x_{i+\frac{1}{2}}\right]$ and $\left[y_{j-\frac{1}{2}}, y_{j+\frac{1}{2}}\right]$, respectively, and $\left\{\omega_{\mu}^{G}\right\}_{\mu=1}^{q}$ are the associated weights for the interval $[-\frac12,\frac12]$.

For any test function $\boldsymbol{v}\left(\boldsymbol{x}\right) \in \mathbb{V}_{Q}^k$, the approximate solution $\boldsymbol{Q}_h\left(\boldsymbol{x}, t\right)$ is expected to satisfy
\begin{equation}\label{2DQapprosolu}
	\begin{aligned}
		&\frac{\mathrm{d}}{\mathrm{d} t} \int_{I_{ij}} \boldsymbol{Q}_{h}\left(\boldsymbol{x}, t\right) \cdot \boldsymbol{v}\left(\boldsymbol{x}\right)\mathrm{~d}\boldsymbol{x}
		=\int_{I_{ij}} \boldsymbol{F}_{1}^{Q}\left(\boldsymbol{U}_{h}\left(\boldsymbol{x}, t\right)\right) \cdot \frac{\partial \boldsymbol{v}\left(\boldsymbol{x}\right)}{\partial x} \mathrm{~d} \boldsymbol{x}+\int_{I_{ij}} \boldsymbol{F}_{2}^{Q}\left(\boldsymbol{U}_{h}\left(\boldsymbol{x}, t\right)\right) \cdot \frac{\partial \boldsymbol{v}\left(\boldsymbol{x}\right)}{\partial y} \mathrm{~d} \boldsymbol{x}\\
		&-\int_{y_{j-\frac{1}{2}}}^{y_{j+\frac{1}{2}}}\left(
		\boldsymbol{v}\left(x_{i+\frac{1}{2}}^{-}, y\right)\cdot\widehat{\boldsymbol{F}_{1}^{Q}}_{, i+\frac{1}{2}}
		\left(y, t\right)
		-\boldsymbol{v}\left(x_{i-\frac{1}{2}}^{+}, y\right)\cdot\widehat{\boldsymbol{F}_{1}^{Q}}_{, i-\frac{1}{2}}
		\left(y, t\right)
		\right)\mathrm{~d}y \\
		&-\int_{x_{i-\frac{1}{2}}}^{x_{i+\frac{1}{2}}}\left(
		\boldsymbol{v}\left(x, y_{j+\frac{1}{2}}^{-}\right)\cdot\widehat{\boldsymbol{F}_{2}^{Q}}_{, j+\frac{1}{2}}
		\left(x, t\right)
		-\boldsymbol{v}\left(x, y_{j-\frac{1}{2}}^{+}\right)\cdot\widehat{\boldsymbol{F}_{2}^{Q}}_{, j-\frac{1}{2}}
		\left(x, t\right)
		\right)\mathrm{~d}x,
	\end{aligned}
\end{equation}
where
\begin{equation*}\label{Equ:2DF1Qequ}
	\widehat{\boldsymbol{F}_{1}^{Q}}_{, i+\frac{1}{2}}\left(y, t\right):=
	\widehat{\boldsymbol{F}_{1}^{Q}}\left(
	\boldsymbol{U}_{h}\left(x_{i+\frac{1}{2}}^{-}, y, t\right),
	\boldsymbol{U}_{h}\left(x_{i+\frac{1}{2}}^{+}, y, t\right)\right), \\
\end{equation*}
\begin{equation*}\label{Equ:2DF2Qequ}
	\widehat{\boldsymbol{F}_{2}^{Q}}_{, j+\frac{1}{2}}\left(x, t\right):=
	\widehat{\boldsymbol{F}_{2}^{Q}}\left(
	\boldsymbol{U}_{h}\left(x, y_{j+\frac{1}{2}}^{-}, t\right),
	\boldsymbol{U}_{h}\left(x, y_{j+\frac{1}{2}}^{+}, t\right)\right), 
\end{equation*}
with $\widehat{\boldsymbol{F}_{1}^{Q}}$ and $\widehat{\boldsymbol{F}_{2}^{Q}}$ denoting the numerical fluxes in the $x$- and $y$-directions, respectively.
Here, we use the PP HLL fluxes
\begin{equation}\label{Equ:2DHLLfluxQ}
	\widehat{\boldsymbol{F}_{\ell}^{Q}}\left(\boldsymbol{U}^{-},\boldsymbol{U}^{+}\right)=
	\frac{\mathscr{V}_{\ell}^{+}\boldsymbol{F}_{\ell}^{Q}\left(\boldsymbol{U}^{-}\right)-
		\mathscr{V}_{\ell}^{-}\boldsymbol{F}_{\ell}^{Q}\left(\boldsymbol{U}^{+}\right)+
		\mathscr{V}_{\ell}^{-}\mathscr{V}_{\ell}^{+}\left(\boldsymbol{Q}^{+}-\boldsymbol{Q}^{-}\right)}
	{\mathscr{V}_{\ell}^{+}-\mathscr{V}_{\ell}^{-}}, \quad \ell=1, 2,
\end{equation}
where $\mathscr{V}_{\ell}^{\pm}$ are the same as those in equation (\ref{Equ:2DHLLfluxR}).
Taking the test function $\boldsymbol{v}\left(\boldsymbol{x}\right)$ as the orthogonal LDF basis function $\boldsymbol{\varphi}_{ij}^{\left(m\right)}\left(\boldsymbol{x}\right)$ and using the same Gauss quadrature to approximate the right-hand side of  (\ref{2DQapprosolu}), we can obtain the following semi-discrete scheme
\begin{equation}\label{Equ:2DsemidiscreQ}
	\begin{aligned}
		&\left(\int_{I_{ij}} \left(\boldsymbol{\varphi}_{ij}^{\left(m\right)}\left(\boldsymbol{x}\right)
		\right)^{2} \mathrm{d} \boldsymbol{x}\right)
		\frac{\mathrm{d}Q_{ij}^{\left(m\right)}
			\left(t\right)}{\mathrm{d} t}
		=\Delta x \Delta y\sum_{\mu=1}^{q}\sum_{\eta=1}^{q}\omega_{\mu}^{G}\omega_{\eta}^{G} \boldsymbol{F}_{1}^{Q}\left(\boldsymbol{U}_{h}\left(x_{i}^{(\mu)}, y_{j}^{(\eta)}, t\right)\right)\cdot \partial_{x}\boldsymbol{\varphi}_{ij}^{\left(m\right)}
		\left(x_{i}^{(\mu)}, y_{j}^{(\eta)}\right)\\
		&+\Delta x \Delta y\sum_{\mu=1}^{q}\sum_{\eta=1}^{q}\omega_{\mu}^{G}\omega_{\eta}^{G} \boldsymbol{F}_{2}^{Q}\left(\boldsymbol{U}_{h}\left(x_{i}^{(\mu)}, y_{j}^{(\eta)}, t\right)\right)\cdot \partial_{y}\boldsymbol{\varphi}_{ij}^{\left(m\right)}
		\left(x_{i}^{(\mu)}, y_{j}^{(\eta)}\right)
		\\
		&-\Delta y \sum_{\mu=1}^{q}\omega_{\mu}^{G}\left(
		\boldsymbol{\varphi}_{ij}^{\left(m\right)}\left(x_{i+\frac{1}{2}}^{-}, y_{j}^{(\mu)}\right)\cdot\widehat{\boldsymbol{F}_{1}^{Q}}_{, i+\frac{1}{2}}
		\left(y_{j}^{(\mu)}, t\right)
		-\boldsymbol{\varphi}_{ij}^{\left(m\right)}\left(x_{i-\frac{1}{2}}^{+}, y_{j}^{(\mu)}\right)\cdot\widehat{\boldsymbol{F}_{1}^{Q}}_{, i-\frac{1}{2}}
		\left(y_{j}^{(\mu)}, t\right)
		\right)\\
		&-\Delta x\sum_{\mu=1}^{q}\omega_{\mu}^{G}\left(
		\boldsymbol{\varphi}_{ij}^{\left(m\right)}\left(x_{i}^{(\mu)}, y_{j+\frac{1}{2}}^{-}\right)\cdot\widehat{\boldsymbol{F}_{2}^{Q}}_{, j+\frac{1}{2}}
		\left(x_{i}^{(\mu)}, t\right)
		-\boldsymbol{\varphi}_{ij}^{\left(m\right)}\left(x_{i}^{(\mu)}, y_{j-\frac{1}{2}}^{+}\right)\cdot\widehat{\boldsymbol{F}_{2}^{Q}}_{, j-\frac{1}{2}}
		\left(x_{i}^{(\mu)}, t\right)
		\right), \\
		&\quad\quad\quad\quad\quad\quad\quad\quad\quad\quad\quad\quad\quad\quad
		\quad\quad\quad\quad\quad\quad\quad\quad\quad\quad\quad\quad
		m=0, 1, \cdots, D_{Q}^{k}-1,
	\end{aligned}
\end{equation}
where the Gauss quadrature nodes and weights are the same as those in \eqref{Equ:2DsemidiscreR}.

\subsection{Fully discrete OEDG schemes}

The combination of semi-discrete schemes (\ref{Equ:2DsemidiscreR}) and (\ref{Equ:2DsemidiscreQ}) can be regarded as
a system of ordinary differential equations (ODEs) 
with respect to time $t$.
For brevity, we rewrite equations (\ref{Equ:2DsemidiscreR}) and (\ref{Equ:2DsemidiscreQ}) in the following ODE form:
\begin{equation}\label{Equ:2DRQODEform}
	\frac{\mathrm{d}}{\mathrm{d} t}\boldsymbol{U}_{h}=\boldsymbol{L}\left(\boldsymbol{U}_{h}\right),
\end{equation}
where $\boldsymbol{U}_{h}$ denotes the combination of $\boldsymbol{R}_{h}$ and $\boldsymbol{Q}_{h}$.
This ODE system \eqref{Equ:2DRQODEform} is typically discretized in time by a Runge--Kutta method, for example, the explicit third-order strong-stability-preserving (SSP) Runge--Kutta method:
\begin{equation*}
	\left\{
	\begin{aligned}
		\boldsymbol{U}_{h}^{n, 0} &= \boldsymbol{U}_{h}^{n}, \\
		\boldsymbol{U}_{h}^{n, 1} &= \boldsymbol{U}_{h}^{n, 0}+\Delta t_{n} \boldsymbol{L}\left(\boldsymbol{U}_{h}^{n, 0}\right), \\
		\boldsymbol{U}_{h}^{n, 2} &= \frac{3}{4}\boldsymbol{U}_{h}^{n, 0}+\frac{1}{4}\left(\boldsymbol{U}_{h}^{n, 1}+\Delta t_{n} \boldsymbol{L}\left(\boldsymbol{U}_{h}^{n, 1}\right)\right), \\
		\boldsymbol{U}_{h}^{n, 3} &= \frac{1}{3}\boldsymbol{U}_{h}^{n, 0}+\frac{2}{3}\left(\boldsymbol{U}_{h}^{n, 2}+\Delta t_{n} \boldsymbol{L}\left(\boldsymbol{U}_{h}^{n, 2}\right)\right), \\
		\boldsymbol{U}_{h}^{n+1} &= \boldsymbol{U}_{h}^{n, 3}, \\
	\end{aligned}
	\right.
\end{equation*}
where $\boldsymbol{U}_{h}^{n}$ denotes the numerical solution at the $n$th time level, and $\Delta t_{n}$ is the time step-size.

The LDF Runge--Kutta DG method described above works well in smooth regions, but can produce serious numerical nonphysical oscillations around strong discontinuities.
To address this issue, we propose an LDF OE procedure after each Runge--Kutta stage
\begin{equation*}\label{Equ:1DOEprocedure}
	\boldsymbol{U}_\sigma^{n, \ell}= \mathcal{F}_{\text{OE}} \boldsymbol{U}_h^{n, \ell}, \quad \ell=1, 2, 3,
\end{equation*}
resulting in the fully discrete LDF OEDG method 
\begin{equation}\label{Equ:2DRKmethodRQ}
	\left\{
	\begin{aligned}
		\boldsymbol{U}_{\sigma}^{n, 0} &= \boldsymbol{U}_{\sigma}^{n}, \\
		\boldsymbol{U}_{h}^{n, 1} &=
		\boldsymbol{U}_{\sigma}^{n, 0}+\Delta t_{n} \boldsymbol{L}\left(\boldsymbol{U}_{\sigma}^{n, 0}\right), \qquad \boldsymbol{U}_{\sigma}^{n, 1} = \mathcal{F}_{\text{OE}} \boldsymbol{U}_{h}^{n, 1}, \\
		\boldsymbol{U}_{h}^{n, 2} &=\frac{3}{4}\boldsymbol{U}_{\sigma}^{n, 0}+\frac{1}{4}\left(\boldsymbol{U}_{\sigma}^{n, 1}+\Delta t_{n} \boldsymbol{L}\left(\boldsymbol{U}_{\sigma}^{n, 1}\right)\right), \qquad \boldsymbol{U}_{\sigma}^{n, 2} = \mathcal{F}_{\text{OE}} \boldsymbol{U}_{h}^{n, 2},\\
		\boldsymbol{U}_{h}^{n, 3} & =\frac{1}{3}\boldsymbol{U}_{\sigma}^{n, 0}+\frac{2}{3}\left(\boldsymbol{U}_{\sigma}^{n, 2}+\Delta t_{n} \boldsymbol{L}\left(\boldsymbol{U}_{\sigma}^{n, 2}\right)\right), \qquad \boldsymbol{U}_{\sigma}^{n, 3} = \mathcal{F}_{\text{OE}}  \boldsymbol{U}_{h}^{n, 3},\\
		\boldsymbol{U}_{\sigma}^{n+1} &= \boldsymbol{U}_{\sigma}^{n, 3}, \\
	\end{aligned}
	\right.
\end{equation}
where $\boldsymbol{U}_{\sigma}^{n}$ is the OE numerical solution at the $n$-th time level, and we define $\boldsymbol{U}_{\sigma}^{0} = \boldsymbol{U}_{h}^{0}$ or $\boldsymbol{U}_{\sigma}^{0} = \mathcal{F}_{\text{OE}} \boldsymbol{U}_{h}^{0}$. 

In the following, we will elaborate on the OE procedures $\mathcal{F}_{\text{OE}}$ for $\boldsymbol{R}_{h}$ and $\boldsymbol{Q}_{h}$, respectively.  It is important to note that the LDF OE procedure we will introduce is specifically designed to maintain the LDF constraint on the approximate magnetic field.

\subsubsection{Component-wise OE procedure for $\boldsymbol{R}_{h}$}

The OE procedure applied for $\boldsymbol{R}_h$ is denoted as
\begin{equation}\label{Equ:2DOEprocedureR}
	\boldsymbol{R}_\sigma^{n, \ell}= \mathcal{F}_{\text{OE}}^R \boldsymbol{R}_h^{n, \ell}, \quad \ell=1, 2, 3.
\end{equation}
The OE operator $\mathcal{F}_{\text{OE}}^R$ is defined by $\left(\mathcal{F}_{\text{OE}}^R \boldsymbol{R}_h^{n, \ell}\right)\left(\boldsymbol{x}\right)=\boldsymbol{R}_\sigma\left(
\boldsymbol{x},  \tau \right)\big|_{\tau = \Delta t_{n}}$ with
$\boldsymbol{R}_\sigma\left(\boldsymbol{x}, \tau\right) \in \mathbb{V}_{R}^k$ being the solution to the following damping equations:
\begin{equation}\label{Equ:2DOERdampequ}
	\left\{
	\begin{aligned}
		&\frac{\mathrm{d}}{\mathrm{d} \tau} \int_{I_{ij}} \boldsymbol{R}_\sigma v \mathrm{~d} \boldsymbol{x}+\sum_{m=0}^k \boldsymbol{\delta}_{ij}^m\left(\boldsymbol{R}_h^{n, \ell}\right)\circ \int_{I_{ij}}\left(\boldsymbol{R}_\sigma-P^{m-1} \boldsymbol{R}_\sigma\right) v \mathrm{~d} \boldsymbol{x}=\boldsymbol{0} \quad \forall v \in \mathbb{P}^k\left(I_{ij}\right), \\
		&\boldsymbol{R}_\sigma\left(\boldsymbol{x}, 0\right)=\boldsymbol{R}_h^{n, \ell}\left(\boldsymbol{x}\right),\\
	\end{aligned}
	\right.
\end{equation}
where
$\tau$ is a pseudo-time different from $t$, $``\circ"$ denotes the Hadamard product of two vectors,
$P^{m}$ is the standard $L^{2}$-projection operator into $\mathbb{V}_{R}^m$ for $m\ge 0$, and we define $P^{-1}=P^{0}$.
The vector $\boldsymbol{\delta}_{ij}^m\left(\boldsymbol{R}_h^{n, \ell}\right)$ is defined by
\begin{equation}\label{Equ:2DRdampterm}
	\boldsymbol{\delta}_{ij}^m\left(\boldsymbol{R}_h\right):=
	\frac{\beta_{ij}^{x}\left(\boldsymbol{\sigma}_{i+\frac{1}{2},j}^m\left(\boldsymbol{R}_h\right)
		+\boldsymbol{\sigma}_{i-\frac{1}{2},j}^m\left(\boldsymbol{R}_h\right)\right)}
	{\Delta x}+\frac{
		\beta_{ij}^{y}\left(\boldsymbol{\sigma}_{i,j+\frac{1}{2}}^m\left(\boldsymbol{R}_h\right)+
		\boldsymbol{\sigma}_{i,j-\frac{1}{2}}^m\left(\boldsymbol{R}_h\right)\right)}
	{\Delta y},
\end{equation}
where $\beta_{ij}^{x}$ and $\beta_{ij}^{y}$ denote suitable estimates of the local maximum wave speed in the $x$- and $y$-directions, respectively,
and $\boldsymbol{\sigma}_{i+\frac{1}{2},j}^m\left(\boldsymbol{R}_h\right)$ and $\boldsymbol{\sigma}_{i,j+\frac{1}{2}}^m\left(\boldsymbol{R}_h\right)$ correspond to the damping coefficients on the $x=x_{i+\frac{1}{2}}$ and $y=y_{j+\frac{1}{2}}$ edges, respectively.
In practice, we take
$\beta_{ij}^{x}$ and $\beta_{ij}^{y}$ as the spectral radii of
the Jacobian matrices $\frac{\partial \boldsymbol{F}_{1}}{\partial {\boldsymbol{U}}}\left( \bar{\boldsymbol{U}}_{ij}^{n, \ell}\right)$ and $\frac{\partial \boldsymbol{F}_{2}}{\partial {\boldsymbol{U}}}\left( \bar{\boldsymbol{U}}_{ij}^{n, \ell}\right)$,
respectively, and $\bar{\boldsymbol{U}}_{ij}^{n, \ell}$ denotes the cell average of
$\boldsymbol{U}_{h}^{n, \ell}\left(\boldsymbol{x}\right)$ over $I_{ij}$.
In equation (\ref{Equ:2DRdampterm}), the damping coefficients are given by
\begin{equation*}\label{Equ:2DRdampcoexx}
	\boldsymbol{\sigma}_{i+\frac{1}{2},j}^m\left(\boldsymbol{R}_h\right)=\left(\sigma_{i+\frac{1}{2},j}^m
	\left(R_h^{\left(1\right)}\right), \sigma_{i+\frac{1}{2},j}^m\left(R_h^{\left(2\right)}\right), \cdots, \sigma_{i+\frac{1}{2},j}^m\left(R_h^{(6)}\right)\right)^{\top}
\end{equation*}
\begin{equation*}\label{Equ:2DRdampcoeyy}
	\boldsymbol{\sigma}_{i,j+\frac{1}{2}}^m\left(\boldsymbol{R}_h\right)=\left(\sigma_{i,j+\frac{1}{2}}^m
	\left(R_h^{\left(1\right)}\right), \sigma_{i,j+\frac{1}{2}}^m\left(R_h^{\left(2\right)}\right), \cdots, \sigma_{i,j+\frac{1}{2}}^m\left(R_h^{\left(6\right)}\right)\right)^{\top},
\end{equation*}
with $R_h^{(\mu)}$ being the $\mu$-th component of $\boldsymbol{R}_h$, and
\begin{align}\label{Equ:2DRdampcoefcompoyy}
	\sigma_{i+\frac{1}{2},j}^m\left(R_h^{(\mu)}\right)=
	\left\{
	\begin{aligned}
		&0, &\text { if } R_h^{\left(\mu\right)} \equiv \operatorname{avg}\left(R_h^{\left(\mu\right)}\right), \\
		&\frac{\left(2 m+1\right) {(\Delta x)}^{m}}{2\left(2 k-1\right) m !} \sum_{\left|\boldsymbol{\alpha}\right|=m}
		\frac
		{\frac{1}{\Delta y}\displaystyle\int_{y_{j-\frac{1}{2}}}^{y_{j+\frac{1}{2}}}\left|\left[\!\left[{\partial}^{\boldsymbol{\alpha}}R_{h}^{\left(\mu\right)}
			\right]\!\right]_{i+\frac{1}{2},j}\right|\mathrm{d}y}
		{\left\|R_h^{\left(\mu\right)}-\operatorname{avg}\left(R_h^{\left(\mu\right)}\right)
			\right\|_{L^{\infty}\left(\Omega\right)}},&\text {otherwise,}
	\end{aligned}
	\right.
	\\
	\label{Equ:2DRdampcoefcompoxx}
	\sigma_{i,j+\frac{1}{2}}^m\left(R_h^{(\mu)}\right)=
	\left\{
	\begin{aligned}
		&0, &\text { if } R_h^{\left(\mu\right)} \equiv \operatorname{avg}\left(R_h^{\left(\mu\right)}\right), \\
		&\frac{\left(2 m+1\right) {(\Delta y)}^{m}}{2\left(2 k-1\right) m !} \sum_{\left|\boldsymbol{\alpha}\right|=m}
		\frac{\frac{1}{\Delta x}\displaystyle\int_{x_{i-\frac{1}{2}}}^{x_{i+\frac{1}{2}}}\left|\left[\!\left[{\partial}^{\boldsymbol{\alpha}}R_{h}^{\left(\mu\right)}
			\right]\!\right]_{i,j+\frac{1}{2}}\right|\mathrm{d}x}
		{\left\|R_h^{\left(\mu\right)}-\operatorname{avg}\left(R_h^{\left(\mu\right)}\right)
			\right\|_{L^{\infty}\left(\Omega\right)}},&\text {otherwise.}
	\end{aligned}
	\right.
\end{align}
In equations (\ref{Equ:2DRdampcoefcompoyy})-(\ref{Equ:2DRdampcoefcompoxx}),
the vector $\boldsymbol{\alpha} = \left(\alpha_{1}, \alpha_{2}\right)$ is the
multi-index,
${\partial}^{\boldsymbol{\alpha}}R_{h}^{\left(\mu\right)}$ is defined as
\begin{equation*}\label{Equ:2DpartialR}
	\partial^{\boldsymbol{\alpha}} R_h^{\left(\mu\right)}\left(\boldsymbol{x}\right):=
	\frac{\left|\boldsymbol{\alpha}\right|!}{\alpha_1 ! \alpha_2 !}
	\frac{\partial^{\left|\boldsymbol{\alpha}\right|}}
	{{\partial x}^{\alpha_1} {\partial y}^{\alpha_2}} R_h^{\left(\mu\right)}\left(\boldsymbol{x}\right),
\end{equation*}
and
${\left[\!\left[{\partial}^{\boldsymbol{\alpha}}R_{h}^{\left(\mu\right)}\right]\!\right]}_
{i+\frac{1}{2},j}$ and ${\left[\!\left[{\partial}^{\boldsymbol{\alpha}}R_{h}^{\left(\mu\right)}
	\right]\!\right]}_{i,j+\frac{1}{2}}$ represent the jump values of ${\partial}^{\boldsymbol{\alpha}}R_{h}^{\left(\mu\right)}$ across the $x=x_{i+\frac{1}{2}}$ and $y=y_{j+\frac{1}{2}}$ edges, respectively.

Note that the coefficient $\boldsymbol{\delta}_{ij}^m\left(\boldsymbol{R}_h^{n, \ell}\right)$ in the damping equations \eqref{Equ:2DOERdampequ} only depends on the ``initial'' value $\boldsymbol{R}_\sigma\left(\boldsymbol{x}, 0\right)=\boldsymbol{R}_h^{n, \ell}\left(\boldsymbol{x}\right)$.
Consequently, the damping equations \eqref{Equ:2DOERdampequ} are actually a linear ODE system, and its exact
solution can be analytically given.
We now derive the exact solver of damping equations \eqref{Equ:2DOERdampequ} for $\boldsymbol{R}_\sigma\left(\boldsymbol{x}, \tau\right)$.
Based on the local orthogonal basis $\left\{\phi_{ij}^{\left(\boldsymbol{\alpha}\right)}\left(\boldsymbol{x}\right):  \left|
\boldsymbol{\alpha}\right| \le k\right\}$ of $\mathbb{P}^k\left(I_{ij}\right)$, the ``initial'' value $\boldsymbol{R}_{h}^{n,\ell}\left(\boldsymbol{x}\right)$ of the  damping equations  (\ref{Equ:2DOERdampequ}) can be expanded as
\begin{equation*}\label{Equ:2DsolutionInitR}
	\boldsymbol{R}_{h}^{n, \ell}\left(\boldsymbol{x}\right)=\sum_{\mu=0}^{k}\sum_{\left|
		\boldsymbol{\alpha}\right|=\mu}\boldsymbol{R}_{ij}^{\left(\boldsymbol{\alpha}\right)}
	\left(0\right)
	\phi_{ij}^{\left(\boldsymbol{\alpha}\right)}\left(\boldsymbol{x}\right) \quad \text{for} \quad \boldsymbol{x} \in I_{ij}.
\end{equation*}
Assume that the solution of (\ref{Equ:2DOERdampequ}) can be expressed as
\begin{equation}\label{Equ:2Dsolution1R}
	\boldsymbol{R}_\sigma\left(\boldsymbol{x}, \tau\right)=\sum_{\mu=0}^{k}\sum_{\left|
		\boldsymbol{\alpha}\right|=\mu}\hat{\boldsymbol{R}}_{ij}^{\left(\boldsymbol{\alpha}\right)}
	\left(\tau\right)
	\phi_{ij}^{\left(\boldsymbol{\alpha}\right)}\left(\boldsymbol{x}\right) \quad \text{for} \quad \boldsymbol{x} \in I_{ij}.
\end{equation}
Note that
\begin{equation}\label{Equ:2DprojectionR}
	\left(\boldsymbol{R}_\sigma-P^{m-1} \boldsymbol{R}_\sigma\right)\left(\boldsymbol{x}, \tau\right)=\sum_{\mu=\max \{m, 1\}}^k \sum_{\left|\boldsymbol{\alpha}\right|=\mu} \hat{\boldsymbol{R}}_{ij}^{\left(\boldsymbol{\alpha}\right)}\left(\tau\right) \phi_{ij}^{\left(\boldsymbol{\alpha}\right)}\left(\boldsymbol{x}\right).
\end{equation}
Substitute (\ref{Equ:2Dsolution1R})-(\ref{Equ:2DprojectionR}) into (\ref{Equ:2DOERdampequ}) and take $v =\phi_{ij}^{\left(\boldsymbol{\alpha}\right)}\left(\boldsymbol{x}\right) $. For $\left|\boldsymbol{\alpha}\right|=\mu \geq 1$, one has
\begin{equation*}\label{Equ:2DdampingR}
	\begin{aligned}
		&\frac{\mathrm{d}}{\mathrm{d} \tau} \hat{\boldsymbol{R}}_{ij}^{\left(\boldsymbol{\alpha}\right)}\left(\tau\right) \int_{I_{ij}}\left(\phi_{ij}^{\left(\boldsymbol{\alpha}\right)}\left(\boldsymbol{x}\right)\right)^2 \mathrm{~d} \boldsymbol{x}+\sum_{m=0}^\mu \boldsymbol{\delta}_{ij}^{m} \left(\boldsymbol{R}_h^{n, \ell}\right)\circ
		\hat{\boldsymbol{R}}_{ij}^{\left(\boldsymbol{\alpha}\right)}\left(\tau\right)  \int_{I_{ij}}\left(\phi_{ij}^{\left(\boldsymbol{\alpha}\right)}\left(\boldsymbol{x}\right)\right)^2
		\mathrm{~d} \boldsymbol{x}=\boldsymbol{0}, 
	\end{aligned}
\end{equation*}
which implies
\begin{equation}\label{Equ:2DsimpdampingR}
	\frac{\mathrm{d}}{\mathrm{d} \tau} \hat{\boldsymbol{R}}_{ij}^{\left(\boldsymbol{\alpha}\right)}\left(\tau\right) +\sum_{m=0}^\mu \boldsymbol{\delta}_{ij}^{m} \left(\boldsymbol{R}_h^{n, \ell}\right)\circ
	\hat{\boldsymbol{R}}_{ij}^{\left(\boldsymbol{\alpha}\right)}\left(\tau\right)  =\boldsymbol{0}, \quad \mu=1,2, \cdots, k,
\end{equation}
Integrating (\ref{Equ:2DsimpdampingR}) from $\tau =0$ to $\Delta t_n$ gives
\begin{equation}\label{Equ:2DsimpdampsoluR} \hat{\boldsymbol{R}}_{ij}^{\left(\boldsymbol{\alpha}\right)}\left(\Delta t_{n}\right)=
	\exp\left( -\Delta t_{n} \sum_{m=0}^\mu \boldsymbol{\delta}_{ij}^{m} \left(\boldsymbol{R}_h^{n, \ell}\right) \right)\circ\boldsymbol{R}_{ij}^{\left(\boldsymbol{\alpha}\right)}\left(0\right), \quad \mu=1,2, \cdots, k.
\end{equation}
For $\left|\boldsymbol{\alpha}\right|=0$, because
$$\int_{I_{ij}} \left(\boldsymbol{R}_\sigma-P^{m-1} \boldsymbol{R}_\sigma\right)\phi_{ij}^{\left(\boldsymbol{\alpha}\right)}\left(\boldsymbol{x}\right)
\mathrm{~d}\boldsymbol{x} =\boldsymbol{0} \quad \forall m \ge 0,$$
and thus
\begin{equation}\label{Equ:2Dsimpdampsolu0R} \hat{\boldsymbol{R}}_{ij}^{\left(\boldsymbol{\alpha}\right)}\left(\Delta t_{n}\right)=
	\boldsymbol{R}_{ij}^{\left(\boldsymbol{\alpha}\right)}\left(0\right) \quad \mbox{for} \quad \boldsymbol{\alpha}={\bf 0} = (0,0).
\end{equation}
Hence, with (\ref{Equ:2DsimpdampsoluR}) and (\ref{Equ:2Dsimpdampsolu0R}), one has
\begin{equation}\label{Equ:2DfinalR}
	\begin{split}
		\boldsymbol{R}_\sigma^{n,\ell}\left(\boldsymbol{x}\right)
		&=\left(\mathcal{F}_{\text{OE}}^R \boldsymbol{R}_h^{n, \ell}\right)\left(\boldsymbol{x}\right)=\boldsymbol{R}_\sigma\left(\boldsymbol{x}, \Delta t_{n}\right)\\
		&=\boldsymbol{R}_{ij}^{\left(\boldsymbol{0}
			\right)}
		\left(0\right)\phi_{ij}^{\left(\boldsymbol{0}\right)}\left(\boldsymbol{x}\right)+
		\sum_{\mu=1}^{k}\exp\left( -\Delta t_{n} \sum_{m=0}^\mu \boldsymbol{\delta}_{ij}^{m} \left(\boldsymbol{R}_h^{n, \ell}\right)\right)\circ\sum_{\left|\boldsymbol{\alpha}\right|=\mu}
		\boldsymbol{R}_{ij}^{\left(\boldsymbol{\alpha}
			\right)}\left(0\right)\phi_{ij}^{\left(\boldsymbol{\alpha}\right)}\left(\boldsymbol{x}\right).
	\end{split}
\end{equation}

\begin{remark}\label{rem1:OE}
In essence, the OE procedure (\ref{Equ:2DfinalR}) suppresses potential spurious oscillations in the DG solution by damping its modal coefficients. Since the OE procedure does not interfere with DG spatial discretization or RK stage updates, it can be easily incorporated into existing DG codes as an independent module and easily extensible onto general meshes. Thanks to the exact solver of the OE procedure, its implementation is straightforward and efficient. Furthermore, the OEDG method remains stable under the normal CFL condition, and many desirable properties of the original DG method are preserved, including conservation and optimal convergence rates (see the theoretical justification in \cite{peng2023oedg} and the numerical evidence in section \ref{Sec:numericaltest}). Note that the OE procedure is applied directly to the conservative variables without the need for characteristic decomposition, which greatly reduces computational costs. As with the OE technique in \cite{peng2023oedg}, the present OE procedure is also scale-invariant and evolution-invariant, thereby being capable of effectively eliminating spurious oscillations for problems across different scales and wave speeds.

\end{remark}

\subsubsection{LDF OE procedure for $\boldsymbol{Q}_{h}$}

After each Runge--Kutta stage, the LDF OE procedure applied for $\boldsymbol{Q}_h$ is denoted as
\begin{equation}\label{Equ:2DOEprocedureQ}
	\boldsymbol{Q}_\sigma^{n, \ell}= \mathcal{F}_{\text{OE}}^Q \boldsymbol{Q}_h^{n, \ell}, \quad \ell=1, 2, 3.
\end{equation}
The LDF OE operator $\mathcal{F}_{\text{OE}}^Q$ is defined by $\left(\mathcal{F}_{\text{OE}}^Q \boldsymbol{Q}_h^{n, \ell}\right)\left(\boldsymbol{x}\right)=\boldsymbol{Q}_\sigma\left(
\boldsymbol{x}, \Delta t_{n}\right)$ with
$\boldsymbol{Q}_\sigma\left(\boldsymbol{x}, \tau\right) \in \mathbb V^k_Q$ ($0\le \tau \le \Delta t_n$) being the solution to the following initial value problem:
\begin{equation}\label{Equ:2DOEQdampequ}
	\left\{
	\begin{aligned}
		&\frac{\mathrm{d}}{\mathrm{d} \tau} \int_{I_{ij}} \boldsymbol{Q}_\sigma \cdot \boldsymbol{v} \mathrm{~d} \boldsymbol{x}+\sum_{m=0}^k \delta_{ij}^m\left(\boldsymbol{Q}_h^{n, \ell}\right) \int_{I_{ij}}\left(\boldsymbol{Q}_\sigma-\mathcal{P}^{m-1} \boldsymbol{Q}_\sigma\right) \cdot \boldsymbol{v} \mathrm{~d} \boldsymbol{x}=0 \quad \forall \boldsymbol{v} \in \mathbb{V}_{Q}^k, \\
		&\boldsymbol{Q}_\sigma\left(\boldsymbol{x}, 0\right)=\boldsymbol{Q}_h^{n, \ell}\left(\boldsymbol{x}\right),\\
	\end{aligned}
	\right.
\end{equation}
where
$\tau$ is a pseudo-time different from $t$,
$\mathcal{P}^{m}$ is the $L^{2}$-projection operator into $\mathbb{V}_{Q}^m$ for $m\ge 0$, and we define $\mathcal{P}^{-1}=\mathcal{P}^{0}$.
The coefficient $\delta_{ij}^m\left(\boldsymbol{Q}_h^{n, \ell}\right)$ is defined as
\begin{equation*}\label{Equ:2DQdampterm}
	\delta_{ij}^m\left(\boldsymbol{Q}_h\right):=
	\frac{\beta_{ij}^{x}\left(\sigma_{i+\frac{1}{2},j}^m\left(\boldsymbol{Q}_h\right)
		+\sigma_{i-\frac{1}{2},j}^m\left(\boldsymbol{Q}_h\right)\right)}
	{\Delta x}+\frac{
		\beta_{ij}^{y}\left(\sigma_{i,j+\frac{1}{2}}^m\left(\boldsymbol{Q}_h\right)+
		\sigma_{i,j-\frac{1}{2}}^m\left(\boldsymbol{Q}_h\right)\right)}
	{\Delta y},
\end{equation*}
where $\beta_{ij}^{x}$ and $\beta_{ij}^{y}$ are the same as those in equation (\ref{Equ:2DRdampterm}),
and $\sigma_{i+\frac{1}{2},j}^m\left(\boldsymbol{Q}_h\right)$ and $\sigma_{i,j+\frac{1}{2}}^m\left(\boldsymbol{Q}_h\right)$ denote the damping coefficients on the $x=x_{i+\frac{1}{2}}$ and $y=y_{j+\frac{1}{2}}$ edges, respectively, given by
\begin{equation*}\label{Equ:2DQdampcoexxyy}
	\sigma_{i+\frac{1}{2},j}^m\left(\boldsymbol{Q}_h\right)=\max_{1\leq\mu\leq 2}\sigma_{i+\frac{1}{2},j}^m\left(Q_h^{(\mu)}\right), \qquad
	\sigma_{i,j+\frac{1}{2}}^m\left(\boldsymbol{Q}_h\right)=\max_{1\leq\mu\leq 2}\sigma_{i,j+\frac{1}{2}}^m\left(Q_h^{(\mu)}\right),
\end{equation*}
with $Q_h^{(\mu)}$ being the $\mu$-th component of $\boldsymbol{Q}_h$, and
\begin{align}\label{Equ:2DQdampcoefcompoyy}
	\sigma_{i+\frac{1}{2},j}^m\left(Q_h^{(\mu)}\right)=
	\left\{
	\begin{aligned}
		&0, &\text { if } Q_h^{\left(\mu\right)} \equiv \operatorname{avg}\left(Q_h^{\left(\mu\right)}\right), \\
		&\frac{\left(2 m+1\right) (\Delta x)^{m}}{2\left(2 k-1\right) m !} \sum_{\left|\boldsymbol{\alpha}\right|=m}
		\frac{\frac{1}{\Delta y}\displaystyle\int_{y_{j-\frac{1}{2}}}^{y_{j+\frac{1}{2}}}\left|
			\left[\!\left[{\partial}^{\boldsymbol{\alpha}}Q_{h}^{\left(\mu\right)}
			\right]\!\right]_{i+\frac{1}{2},j}\right|\mathrm{d}y}
		{\left\|Q_h^{\left(\mu\right)}-\operatorname{avg}\left(Q_h^{\left(\mu\right)}\right)
			\right\|_{L^{\infty}\left(\Omega\right)}},&\text {otherwise,}
	\end{aligned}
	\right.
	\\
	\label{Equ:2DQdampcoefcompoxx}
	\sigma_{i,j+\frac{1}{2}}^m\left(Q_h^{(\mu)}\right)=
	\left\{
	\begin{aligned}
		&0, &\text { if } Q_h^{\left(\mu\right)} \equiv \operatorname{avg}\left(Q_h^{\left(\mu\right)}\right), \\
		&\frac{\left(2 m+1\right) (\Delta y)^{m}}{2\left(2 k-1\right) m !} \sum_{\left|\boldsymbol{\alpha}\right|=m}
		\frac{\frac{1}{\Delta x}\displaystyle\int_{x_{i-\frac{1}{2}}}^{x_{i+\frac{1}{2}}}\left|\left[\!\left[{\partial}^{\boldsymbol{\alpha}}Q_{h}^{\left(\mu\right)}
			\right]\!\right]_{i,j+\frac{1}{2}}\right|\mathrm{d}x}
		{\left\|Q_h^{\left(\mu\right)}-\operatorname{avg}\left(Q_h^{\left(\mu\right)}\right)
			\right\|_{L^{\infty}\left(\Omega\right)}},&\text {otherwise.}
	\end{aligned}
	\right.
\end{align}
In equations (\ref{Equ:2DQdampcoefcompoyy})-(\ref{Equ:2DQdampcoefcompoxx}),
the vector $\boldsymbol{\alpha} = \left(\alpha_{1}, \alpha_{2}\right)$ is the
multi-index,
${\partial}^{\boldsymbol{\alpha}}Q_{h}^{\left(\mu\right)}$ is defined as
\begin{equation*}\label{Equ:2DpartialQ}
	\partial^{\boldsymbol{\alpha}} Q_h^{\left(\mu\right)}\left(\boldsymbol{x}\right):=
	\frac{\left|\boldsymbol{\alpha}\right|!}{\alpha_1 ! \alpha_2 !}
	\frac{\partial^{\left|\boldsymbol{\alpha}\right|}}
	{{\partial x}^{\alpha_1} {\partial y}^{\alpha_2}} Q_h^{\left(\mu\right)}\left(\boldsymbol{x}\right),
\end{equation*}
and
${\left[\!\left[{\partial}^{\boldsymbol{\alpha}}Q_{h}^{\left(\mu\right)}\right]\!\right]}_
{i+\frac{1}{2},j}$ and ${\left[\!\left[{\partial}^{\boldsymbol{\alpha}}Q_{h}^{\left(\mu\right)}
	\right]\!\right]}_{i,j+\frac{1}{2}}$ represent the jump values of ${\partial}^{\boldsymbol{\alpha}}Q_{h}^{\left(\mu\right)}$ across the $x=x_{i+\frac{1}{2}}$ and $y=y_{j+\frac{1}{2}}$ edges, respectively.

We now  formulate explicitly the exact solver of the damping equation (\ref{Equ:2DOEQdampequ}) for $\boldsymbol{Q}_\sigma\left(\boldsymbol{x}, \tau\right)$.
Based on the local orthogonal basis  $\left\{\boldsymbol{\varphi}_{ij}^{\left(\mu\right)}\left(\boldsymbol{x}\right),
\mu=0, 1, \cdots, D_{Q}^k-1\right\}$ of $\mathbb{V}_{Q}^k$ over cell $I_{ij}$,
the ``initial'' value $\boldsymbol{Q}_\sigma\left(\boldsymbol{x}, 0\right)=\boldsymbol{Q}_h^{n, \ell}\left(\boldsymbol{x}\right)$ of the damping equation (\ref{Equ:2DOEQdampequ}) can be expressed as
\begin{equation*}\label{Equ:2DsolutionInitQ}
	\boldsymbol{Q}_{h}^{n, \ell}\left(\boldsymbol{x}\right)=\sum_{\mu=0}^{k}\sum_{\eta=D_{Q}^{\mu-1}}^{D_{Q}^\mu-1}Q_{ij}^{
		\left(\eta\right)}\left(0\right)\boldsymbol{\varphi}_{ij}^{\left(\eta\right)}
	\left(\boldsymbol{x}\right) \quad \text{for} \quad \boldsymbol{x} \in I_{ij}.
\end{equation*}
Assume that the solution of (\ref{Equ:2DOEQdampequ}) can be represented as
\begin{equation}\label{Equ:2DsolutionQ}
	\boldsymbol{Q}_\sigma\left(\boldsymbol{x}, \tau\right)=\sum_{\mu=0}^{k}\sum_{\eta=D_{Q}^{\mu-1}}^{D_{Q}^\mu-1}
	\hat{Q}_{ij}^{
		\left(\eta\right)}\left(\tau\right)\boldsymbol{\varphi}_{ij}^{\left(\eta\right)}
	\left(\boldsymbol{x}\right) \quad \text{for} \quad \boldsymbol{x} \in I_{ij}.
\end{equation}
Note that
\begin{equation}\label{Equ:2DprojectionQ}
	\left(\boldsymbol{Q}_\sigma-\mathcal{P}^{m-1} \boldsymbol{Q}_\sigma\right)\left(\boldsymbol{x}, \tau\right)=\sum_{\mu=\max \{m, 1\}}^k\sum_{\eta=D_{Q}^{\mu-1}}^{D_{Q}^\mu-1}\hat{Q}_{ij}^{
		\left(\eta\right)}\left(\tau\right)\boldsymbol{\varphi}_{ij}^{\left(\eta\right)}
	\left(\boldsymbol{x}\right) \quad \text{for} \quad \boldsymbol{x} \in I_{ij}.
\end{equation}
Substitute (\ref{Equ:2DsolutionQ})-(\ref{Equ:2DprojectionQ}) into (\ref{Equ:2DOEQdampequ}) and take $\boldsymbol{v} =\boldsymbol{\varphi}_{ij}^{\left(\eta \right)}\left(\boldsymbol{x}\right) $. For any $\mu \ge 1$ and $D_{Q}^{\mu-1} \le \eta \le D_Q^{\mu}-1$, one obtains
\begin{equation*}\label{Equ:2DdampingQ}
	\begin{split}
		&\frac{\mathrm{d}}{\mathrm{d} \tau} \hat{Q}_{ij}^{\left(\eta \right)}\left(\tau\right)
		\int_{I_{ij}}\boldsymbol{\varphi}_{ij}^{\left(\eta \right)}\left(\boldsymbol{x}\right)\cdot
		\boldsymbol{\varphi}_{ij}^{\left(\eta \right)}
		\left(\boldsymbol{x}\right) \mathrm{~d} \boldsymbol{x}+\sum_{m=0}^\mu \delta_{ij}^{m} \left(\boldsymbol{Q}_h^{n, \ell}\right)
		\hat{Q}_{ij}^{\left(\eta \right)}\left(\tau\right)  \int_{I_{ij}}\boldsymbol{\varphi}_{ij}^{\left(\eta \right)}\left(\boldsymbol{x}\right)\cdot
		\boldsymbol{\varphi}_{ij}^{\left(\eta\right)}\left(\boldsymbol{x}\right)\mathrm{~d} \boldsymbol{x}=0,
	\end{split}
\end{equation*}
which can be simplified as
\begin{equation}\label{Equ:2DsimpdampingQ}
	\frac{\mathrm{d}}{\mathrm{d} \tau} \hat{Q}_{ij}^{\left(\eta\right)}\left(\tau\right) +\sum_{m=0}^\mu \delta_{ij}^{m} \left(\boldsymbol{Q}_h^{n, \ell}\right)
	\hat{Q}_{ij}^{\left(\eta\right)}\left(\tau\right)  =0, \quad D_{Q}^{\mu-1} \le \eta \le D_Q^{\mu}-1.
\end{equation}
Solving (\ref{Equ:2DsimpdampingQ}) gives
\begin{equation}\label{Equ:2DsimpdampsoluQ}
	\hat{Q}_{ij}^{\left(\eta\right)}\left(\Delta t_{n}\right)=
	\exp\left( -\Delta t_{n} \sum_{m=0}^\mu \delta_{ij}^{m} \left(\boldsymbol{Q}_h^{n, \ell}\right) \right)Q_{ij}^{\left(\eta\right)}\left(0\right), \quad D_{Q}^{\mu-1} \le \eta \le D_Q^{\mu}-1.
\end{equation}
For $0 \le \eta \le D_{Q}^0 -1 =1$,
because
$$\int_{I_{ij}} \left(\boldsymbol{Q}_\sigma-\mathcal{P}^{m-1} \boldsymbol{Q}_\sigma\right)\cdot \boldsymbol{\varphi}_{ij}^{\left(\eta\right)}\left(\boldsymbol{x}\right)
\mathrm{~d}\boldsymbol{x} \equiv0 \quad \forall m \ge 0,$$ we have
$$\frac{\mathrm{d}}{\mathrm{d} \tau} \hat{Q}_{ij}^{\left(\eta\right)}\left(\tau\right) \int_{I_{ij}}\boldsymbol{\varphi}_{ij}^{\left(\eta\right)}\left(\boldsymbol{x}\right)
\cdot \boldsymbol{\varphi}_{ij}^{\left(\eta\right)}\left(\boldsymbol{x}\right)\mathrm{~d} \boldsymbol{x}=0, \quad 0 \le \eta \le D_{Q}^0 -1 =1,$$
and thus
\begin{equation}\label{Equ:2Dsimpdampsolu0Q}
	\hat{Q}_{ij}^{\left(\eta\right)}\left(\Delta t_{n}\right)=
	Q_{ij}^{\left(\eta\right)}\left(0\right), \quad \eta=0, \cdots, D_{Q}^0-1.
\end{equation}
Hence, with (\ref{Equ:2DsimpdampsoluQ}) and (\ref{Equ:2Dsimpdampsolu0Q}), one has
\begin{equation}\label{Equ:2DfinalQ}
	\begin{split}
		\boldsymbol{Q}_\sigma^{n,\ell}\left(\boldsymbol{x}\right)
		&=\left(\mathcal{F}_{\text{OE}}^Q \boldsymbol{Q}_h^{n, \ell}\right)\left(\boldsymbol{x}\right)=\boldsymbol{Q}_\sigma\left(\boldsymbol{x}, \Delta t_{n}\right)\\
		&=\sum_{\eta=0}^{D_{Q}^0-1}Q_{ij}^{\left(\eta\right)}
		\left(0\right)\boldsymbol{\varphi}_{ij}^{\left(\eta\right)}\left(\boldsymbol{x}\right)+
		\sum_{\mu=1}^{k}\exp\left( -\Delta t_{n} \sum_{m=0}^\mu \delta_{ij}^{m} \left(\boldsymbol{Q}_h^{n,\ell}\right)\right)
		\sum_{\eta=D_{Q}^{\mu-1}}^{D_{Q}^\mu-1}Q_{ij}^{
			\left(\eta\right)}\left(0\right)\boldsymbol{\varphi}_{ij}^{\left(\eta\right)}
		\left(\boldsymbol{x}\right).
	\end{split}
\end{equation}

\begin{remark}
The OE procedure \eqref{Equ:2DfinalQ} not only effectively suppresses potential spurious oscillations but also naturally maintains the LDF property of the magnetic field. Moreover, this LDF OE procedure inherits many notable advantages of the component-wise OE procedure (\ref{Equ:2DfinalR}) mentioned in Remark \ref{rem1:OE}, including:
being non-intrusive, which enables easy integration into existing (LDF) DG codes; 
straightforward extensibility to general meshes;
stability under normal CFL conditions; 
scale invariance and evolution invariance, which are crucial for effectively eliminating spurious oscillations across various scales and wave speeds; 
being free of characteristic decomposition; 
retaining  the conservation and optimal convergence rates of the original DG method. 
All these advantages underscore the strengths of the LDF OE procedure.

\end{remark}

\section{Positivity preservation for LDF OEDG method}\label{sec:PPLDFOEDGMHD}

For physical solutions, both the density and internal energy (or pressure) must be positive. We define the set of physically admissible states as follows:
\begin{equation}\label{Equ:phyadmiset}
	\mathcal{G}=\left\{\boldsymbol{U}=\left(\rho, \boldsymbol{m}, \boldsymbol{B}, E\right)^{\top}: \rho>0,~ \mathcal{E}\left(\boldsymbol{U}\right):=E-\frac{1}{2}\left(
	\frac{{\|\boldsymbol{m}\|}^{2}}{\rho} + {\|\boldsymbol{B}\|}^{2}
	\right)>0\right\},
\end{equation}
where $\mathcal{E}\left(\boldsymbol{U}\right)$ denotes the internal energy. 
When the polynomial degree $k \geq 1$, the OEDG method does not always ensure that $\boldsymbol{U}_{\sigma}^{n} \in \mathcal{G}$. Based on the Zhang--Shu framework \cite{ZHANG20103091,ZHANG20108918}, to achieve a bound-preserving DG scheme for hyperbolic conservation laws, one must first identify a DG scheme that preserves the updated cell averages within the bounds. Subsequently, a local scaling bound-preserving limiter is used to ensure that the DG solution satisfies these bounds at all the points of interest. 
According to discoveries in \cite{Wu2018Positivity,Wu2019Provably}, for the conservative MHD system \eqref{Equ:govequ2D}, the PP property of the updated cell averages of the DG solution is closely connected with a set of globally coupled discrete divergence-free conditions, which, however, are incompatible with the standard local scaling limiter. An effective way \cite{Wu2018A,Wu2019Provably} to address this issue is to consider the following symmetrizable MHD system with an additional source term: 
\begin{equation}\label{Equ:govequsource2D}
	\frac{\partial \boldsymbol{U}}{\partial t}+\frac{\partial\boldsymbol{F}_{1}\left(\boldsymbol{U}\right)}{\partial{x}}
	+\frac{\partial\boldsymbol{F}_{2}\left(\boldsymbol{U}\right)}{\partial{y}}=
	-\left(\nabla \cdot \boldsymbol{B}\right)\boldsymbol{S}\left(\boldsymbol{U}\right),
\end{equation}
where $\boldsymbol{S}\left(\boldsymbol{U}\right)=\left(0, \boldsymbol{B}, \boldsymbol{u}, \boldsymbol{u}\cdot\boldsymbol{B}\right)^{\top}$, and the right-hand side of (\ref{Equ:govequsource2D}) is termed the Godunov--Powell source term \cite{Godunov1972SymmetricFO,Powell1995AnUS}. As demonstrated in \cite{Wu2018A,Wu2019Provably}, this source term can help to eliminate the effect of divergence error on the PP property.

In this section, we will utilize the optimal convex decomposition approach proposed in \cite{Cui2024On,CUI2023111882} to derive an efficient, provably PP LDF OEDG method. This method will be built on a PP HLL flux, a properly upwind discretization of the Godunov--Powell source term, and a simple scaling PP limiter. Our method automatically ensures that the updated cell averages of the OEDG solution remain within $\mathcal{G}$, while the PP limiter enforces the values of the OEDG solution in $\mathcal{G}$ at any points of interest.

\subsection{Positivity-preserving HLL flux}\label{Sec:PPHLLflux}

The HLL fluxes, with $\ell=1$ for the $x$-direction and $\ell=2$ for the $y$-direction,  are given by
\begin{equation}\label{Equ:HLL2D}
	\widehat{\boldsymbol{F}_{\ell}}\left(\boldsymbol{U}^{-}, \boldsymbol{U}^{+} \right)=
	\left\{
	\begin{aligned}
		&\boldsymbol{F}_{\ell}\left(\boldsymbol{U}^{-}\right),  &0 \leq \mathscr{V}_{\ell,l}<\mathscr{V}_{\ell,r}, \\ &\frac{
			\mathscr{V}_{\ell,r}\boldsymbol{F}_{\ell}\left(\boldsymbol{U}^{-}\right)-
			\mathscr{V}_{\ell,l}\boldsymbol{F}_{\ell}\left(\boldsymbol{U}^{+}\right)+
			\mathscr{V}_{\ell,l} \mathscr{V}_{\ell,r}\left(\boldsymbol{U}^{+}-\boldsymbol{U}^{-}\right)}
		{\mathscr{V}_{\ell,r}-\mathscr{V}_{\ell,l}},   &\mathscr{V}_{\ell,l}<0<\mathscr{V}_{\ell,r}, \\ &\boldsymbol{F}_{\ell}\left(\boldsymbol{U}^{+}\right), &\mathscr{V}_{\ell,l}<\mathscr{V}_{\ell,r} \leq 0,
	\end{aligned}
	\right.
\end{equation}
where $\mathscr{V}_{\ell,l}$ and $\mathscr{V}_{\ell,r}$ are functions of $\boldsymbol{U}^{-}$ and $\boldsymbol{U}^{+}$, denoting the estimated minimum and maximum wave speeds, respectively.

The estimates of $\{\mathscr{V}_{\ell,l}\}$ and $\{\mathscr{V}_{\ell,r}\}$ are crucial for the PP property \cite{Wu2019Provably}. 
For any pair of admissible states $\boldsymbol{U}$ and $\tilde{\boldsymbol{U}}$, define
\begin{equation}\label{Equ:alphaleft}
	\alpha_{\ell,l}\left(\boldsymbol{U}, \tilde{\boldsymbol{U}}\right):=\min \left\{u_{\ell},\frac{\sqrt{\rho}u_{\ell}+\sqrt{\tilde{\rho}}\tilde{u}_{\ell}}
	{\sqrt{\rho}+\sqrt{\tilde{\rho}}}\right\}-\mathscr{C}_{\ell}\left(\boldsymbol{U}\right)-
	\frac{\|\boldsymbol{B}-\tilde{\boldsymbol{B}}\|}{\sqrt{\rho}+\sqrt{\tilde{\rho}}},
	\quad \ell=1, 2,
\end{equation}
\begin{equation}\label{Equ:alpharight}
	\alpha_{\ell,r}\left(\boldsymbol{U}, \tilde{\boldsymbol{U}}\right):=\max \left\{u_{\ell},\frac{\sqrt{\rho}u_{\ell}+\sqrt{\tilde{\rho}}\tilde{u}_{\ell}}
	{\sqrt{\rho}+\sqrt{\tilde{\rho}}}\right\}+\mathscr{C}_{\ell}\left(\boldsymbol{U}\right)+
	\frac{\|\boldsymbol{B}-\tilde{\boldsymbol{B}}\|}{\sqrt{\rho}+\sqrt{\tilde{\rho}}},
	\quad\ell=1, 2,
\end{equation}
where
\begin{equation}\label{Equ:speedc}
	\mathscr{C}_{\ell}\left(\boldsymbol{U}\right):=\frac{1}{\sqrt{2}}\left[\mathscr{C}_s^2+
	\frac{\|\boldsymbol{B}\|^2}{\rho}+\sqrt{\left({\mathscr{C}_s}^2+\frac{\|\boldsymbol{B}\|^2}{\rho}
		\right)^2-4 \frac{\mathscr{C}_s^2 {B_{\ell}}^2}{\rho}}\right]^{\frac{1}{2}} \quad \text{with} \quad \mathscr{C}_s=\sqrt{\frac{\left(\gamma-1\right)p}{2\rho}}, \quad \ell=1, 2.
\end{equation}
According to the rigorous PP analysis in \cite{Wu2019Provably},  $\{\mathscr{V}_{\ell,l}\}$ and $\{\mathscr{V}_{\ell,r}\}$ should satisfy
\begin{equation}\label{Equ:PPspeedcond}
	\mathscr{V}_{\ell,l}\leq\alpha_{\ell,l}\left(\boldsymbol{U}^{-}, \boldsymbol{U}^{+}\right), \quad \mathscr{V}_{\ell,r}\geq\alpha_{\ell,r}\left(\boldsymbol{U}^{+}, \boldsymbol{U}^{-}\right), \quad \ell=1, 2.
\end{equation}
We take $\mathscr{V}_{\ell,l}$ and $\mathscr{V}_{\ell,r}$ as
\begin{equation}\label{Equ:WuandShuspeed}
	\begin{aligned}
		\mathscr{V}_{\ell,l} &=\min \left\{\alpha_{\ell,l}\left(\boldsymbol{U}^{-}, \boldsymbol{U}^{+}\right), \lambda_{\ell,\min}\left(\boldsymbol{U}^{-}\right), \lambda_{\ell,\min}\left(\boldsymbol{U}^{\text{Roe}} \right), \lambda_{\ell,\min}\left(\boldsymbol{U}^{+} \right)  \right\}, \quad \ell=1, 2,\\
		\mathscr{V}_{\ell,r} &=\max \left\{\alpha_{\ell,r}\left(\boldsymbol{U}^{+}, \boldsymbol{U}^{-}\right), \lambda_{\ell,\max}\left(\boldsymbol{U}^{-}\right), \lambda_{\ell,\max}\left(\boldsymbol{U}^{\text{Roe}} \right), \lambda_{\ell,\max}\left(\boldsymbol{U}^{+} \right)    \right\}, \quad \ell=1, 2,
	\end{aligned}
\end{equation}
where $\lambda_{\ell,\min/\max}\left(\boldsymbol{U}\right)$ are the minimum and maximum eigenvalues of the Jacobian matrix $\frac{\partial \boldsymbol{F}_{\ell}}{\partial {\boldsymbol{U}}}$, respectively, and  $\lambda_{\ell,\min/\max}\left(\boldsymbol{U}^{\text{Roe}}\right)$ are those of the Roe matrix \cite{Powell1995AnUS}. Let
$\mathscr{V}_{\ell}^{-}=\min\left\{\mathscr{V}_{\ell,l}, 0\right\}$ and
$\mathscr{V}_{\ell}^{+}=\max\left\{\mathscr{V}_{\ell,r}, 0\right\}$,
then the HLL fluxes in (\ref{Equ:HLL2D}) can be rewritten as
\begin{equation}\label{Equ:1DHLLfluxWu}
	\widehat{\boldsymbol{F}_\ell}\left(\boldsymbol{U}^{-},\boldsymbol{U}^{+}\right)=
	\frac{\mathscr{V}_{\ell}^{+}\boldsymbol{F}_\ell\left(\boldsymbol{U}^{-}\right)-
		\mathscr{V}_{\ell}^{-}\boldsymbol{F}_\ell\left(\boldsymbol{U}^{+}\right)+
		\mathscr{V}_{\ell}^{-}\mathscr{V}_{\ell}^{+}\left(\boldsymbol{U}^{+}-\boldsymbol{U}^{-}\right)}
	{\mathscr{V}_{\ell}^{+}-\mathscr{V}_{\ell}^{-}}, \quad \ell=1, 2.
\end{equation}

\subsection{Upwind discretization of Godunov--Powell source term}

According to the theory in \cite{Wu2018A,Wu2019Provably}, incorporating a properly discretized Godunov--Powell source term into the DG schemes is important to ensure the PP property of the updated cell averages. The source terms can help eliminate the effect of the normal magnetic jump across the cell interface on the PP property. Since the PP property of the updated cell averages is related only to the discrete equations satisfied by the cell averages, we consider adding the Godunov--Powell source term only to the evolution equations of the cell averages. Meanwhile, the evolution of the higher-order ``moments" continues to use the OEDG schemes designed in section \ref{sec:2Dgovequcon}.

Let $\boldsymbol{U}_{\sigma}^{n}\left(x, y\right)$ denote the OEDG solution  at the $n$th time level after LDF OE procedure, and let $\bar{\boldsymbol{U}}_{ij}^{n}$ represent its cell average over $I_{ij}$. 
Since an SSP time discretization can be regarded as a convex combination of the forward Euler method, we only need to discuss the PP property of the OEDG method coupled with the forward Euler time discretization. By incorporating an ``upwind" discretized Godunov--Powell source term into the evolution equations for the cell averages, we design
\begin{equation*}\label{Equ:averageMHDsource}
	\bar{\boldsymbol{U}}_{ij}^{n+1}=\bar{\boldsymbol{U}}_{ij}^{n}+\Delta t_{n}
	\boldsymbol{L}_{ij}\left(\boldsymbol{U}_{\sigma}^{n}\right),
\end{equation*}
where
\begin{equation}\label{Equ:averageMHDsourceL}
	\begin{aligned}
		\boldsymbol{L}_{ij}\left(\boldsymbol{U}_{\sigma}^{n}\right)
		&:=-\frac{1}{\Delta x} \sum_{\mu=1}^{q}\omega_{\mu}^{G}\bigg[\left(
		\widehat{\boldsymbol{F}_{1}}_{, i+\frac{1}{2}}
		\left(y_{j}^{(\mu)}\right)-
		\widehat{\boldsymbol{F}_{1}}_{, i-\frac{1}{2}}
		\left(y_{j}^{(\mu)}\right)\right)\\
		&+\mathscr{B}_{1, i+\frac{1}{2}}^-
		\left(y_{j}^{(\mu)}\right)\boldsymbol{S}\left(\boldsymbol{U}_{\sigma}^{n}
		\left(x_{i+\frac{1}{2}}^{-}, y_{j}^{(\mu)}\right)\right)+
		\mathscr{B}_{1, i-\frac{1}{2}}^+
		\left(y_{j}^{(\mu)}\right)\boldsymbol{S}\left(\boldsymbol{U}_{\sigma}^{n}
		\left(x_{i-\frac{1}{2}}^{+}, y_{j}^{(\mu)}\right)\right)
		\bigg]
		\\
		&-\frac{1}{\Delta y} \sum_{\mu=1}^{q}\omega_{\mu}^{G}\bigg[\left(
		\widehat{\boldsymbol{F}_{2}}_{, j+\frac{1}{2}}
		\left(x_{i}^{(\mu)}\right)-
		\widehat{\boldsymbol{F}_{2}}_{, j-\frac{1}{2}}
		\left(x_{i}^{(\mu)}\right)\right)\\
		&+\mathscr{B}_{2, j+\frac{1}{2}}^-
		\left(x_{i}^{(\mu)}\right)\boldsymbol{S}\left(\boldsymbol{U}_{\sigma}^{n}
		\left(x_{i}^{(\mu)}, y_{j+\frac{1}{2}}^{-}\right)\right)+
		\mathscr{B}_{2, j-\frac{1}{2}}^+
		\left(x_{i}^{(\mu)}\right)\boldsymbol{S}\left(\boldsymbol{U}_{\sigma}^{n}
		\left(x_{i}^{(\mu)}, y_{j-\frac{1}{2}}^{+}\right)\right)
		\bigg].
	\end{aligned}
\end{equation}
Here the Gauss quadrature nodes 
$\left\{x_{i}^{(\mu)}\right\}_{\mu=1}^{q}$ and $\left\{y_{j}^{(\mu)}\right\}_{\mu=1}^{q}$ as well as the weights 
$\left\{\omega_{\mu}^{G}\right\}_{\mu=1}^{q}$ are the same as in \eqref{Equ:2DsemidiscreR}. 
The numerical fluxes read 
\begin{align*}
	&\widehat{\boldsymbol{F}_{1}}_{, i+\frac{1}{2}}\left(y\right):=
	\widehat{\boldsymbol{F}_{1}}\left(
	\boldsymbol{U}_{\sigma}^{n}\left(x_{i+\frac{1}{2}}^{-}, y\right),
	\boldsymbol{U}_{\sigma}^{n}\left(x_{i+\frac{1}{2}}^{+}, y\right)\right), \\
	&\widehat{\boldsymbol{F}_{2}}_{, j+\frac{1}{2}}\left(x\right):=
	\widehat{\boldsymbol{F}_{2}}\left(
	\boldsymbol{U}_{\sigma}^{n}\left(x, y_{j+\frac{1}{2}}^{-}\right),
	\boldsymbol{U}_{\sigma}^{n}\left(x, y_{j+\frac{1}{2}}^{+}\right)\right), 
\end{align*}
In (\ref{Equ:averageMHDsourceL}), we take
\begin{align*}
	\mathscr{B}_{1, i+\frac{1}{2}}^\pm\left(y\right) :=
	\frac{\pm \mathscr{V}_{1,i+\frac{1}{2}}^{\pm}\left(y\right)}
	{\mathscr{V}_{1,i+\frac{1}{2}}^{+}\left(y\right)-
		\mathscr{V}_{1,i+\frac{1}{2}}^{-}\left(y\right)}
	\left[{\left(B_{1}\right)}_{\sigma}^{n}\left(x_{i+\frac{1}{2}}^{+}, y\right)-
	{\left(B_{1}\right)}_{\sigma}^{n}\left(x_{i+\frac{1}{2}}^{-}, y\right)\right],\\
	\mathscr{B}_{2, j+\frac{1}{2}}^{\pm}\left(x\right) := \frac{\pm\mathscr{V}_{2,j+\frac{1}{2}}^{\pm}\left(x\right)}
	{\mathscr{V}_{2,j+\frac{1}{2}}^{+}\left(x\right)-
		\mathscr{V}_{2,j+\frac{1}{2}}^{-}\left(x\right)}
	\left[{\left(B_{2}\right)}_{\sigma}^{n}\left(x, y_{j+\frac{1}{2}}^{+}\right)-
	{\left(B_{2}\right)}_{\sigma}^{n}\left(x, y_{j+\frac{1}{2}}^{-}\right)\right]
\end{align*}
with 
\begin{align*}
	\mathscr{V}_{1,i+\frac{1}{2}}^{\pm}\left(y\right):=\mathscr{V}_{1}^{\pm}\left(
	\boldsymbol{U}_{\sigma}^{n}\left(x_{i+\frac{1}{2}}^{-}, y\right),
	\boldsymbol{U}_{\sigma}^{n}\left(x_{i+\frac{1}{2}}^{+}, y\right)\right),\\
	\mathscr{V}_{2,j+\frac{1}{2}}^{\pm}\left(x\right):=\mathscr{V}_{2}^{\pm}\left(
	\boldsymbol{U}_{\sigma}^{n}\left(x, y_{j+\frac{1}{2}}^{-}\right),
	\boldsymbol{U}_{\sigma}^{n}\left(x, y_{j+\frac{1}{2}}^{+}\right)\right),
\end{align*}
where $\mathscr{V}_{1}^{\pm}$ and $\mathscr{V}_{2}^{\pm}$ are the wave speeds in the HLL fluxes discussed in section \ref{Sec:PPHLLflux}. Note that 
\begin{align}\label{eq:split1}
	\mathscr{B}_{1, i+\frac{1}{2}}^+\left(y\right) + \mathscr{B}_{1, i+\frac{1}{2}}^-\left(y\right) &= {\left(B_{1}\right)}_{\sigma}^{n}\left(x_{i+\frac{1}{2}}^{+}, y\right)-
{\left(B_{1}\right)}_{\sigma}^{n}\left(x_{i+\frac{1}{2}}^{-}, y\right)=: \left[\!\left[\left(B_{1}\right)_{\sigma}^{n}
\right]\!\right]_{i+\frac12}(y),
\\ \label{eq:split2}
\mathscr{B}_{2, j+\frac{1}{2}}^+ \left(x\right) + \mathscr{B}_{2, j+\frac{1}{2}}^- \left(x\right) &= {\left(B_{2}\right)}_{\sigma}^{n}\left(x, y_{j+\frac{1}{2}}^{+}\right)-
{\left(B_{2}\right)}_{\sigma}^{n}\left(x, y_{j+\frac{1}{2}}^{-}\right)=: \left[\!\left[\left(B_{2}\right)_{\sigma}^{n}
\right]\!\right]_{j+\frac12}(x)
\end{align}
are actually wave-dependent splittings of the jumps in the normal magnetic component across the cell interfaces $x=x_{i+\frac12}$ and $y=y_{j+\frac12}$, respectively. 
Such ``upwind'' splittings are crucial for the PP property of the updated cell averages; see section \ref{sec:PPproof}. 

More specifically, in \eqref{Equ:averageMHDsourceL}, 
the terms involving $\mathscr{B}_{1, i+\frac{1}{2}}^\pm\left(y\right)$ and  $\mathscr{B}_{2, j+\frac{1}{2}}^{\pm}\left(x\right)$  are actually an ``unwind'' discretization of the cell average of the Godunov--Powell source term 
\begin{equation}\label{key000}
	\frac{1}{\Delta x \Delta y} \int_{y_{j-\frac12}}^{y_{j+\frac12}} \int_{x_{i-\frac12}}^{x_{i+\frac12}} \left(-\nabla \cdot \boldsymbol{B}_{\sigma}^n\right)\boldsymbol{S}\left(\boldsymbol{U}_{\sigma}^n\right)~{\rm d} x {\rm d}y. 
\end{equation}
Thanks to the LDF property of $\boldsymbol{B}_{\sigma}^n$, the local divergence $\nabla \cdot \boldsymbol{B}_{\sigma}^n$ vanishes within each cell $I_{ij}$. Therefore, to approximate \eqref{key000}, we only need to measure the ``weak" divergence/derivatives of the normal magnetic component at the cell interface.

\begin{figure}[!htb]
	\centering
	\includegraphics[scale=0.5]{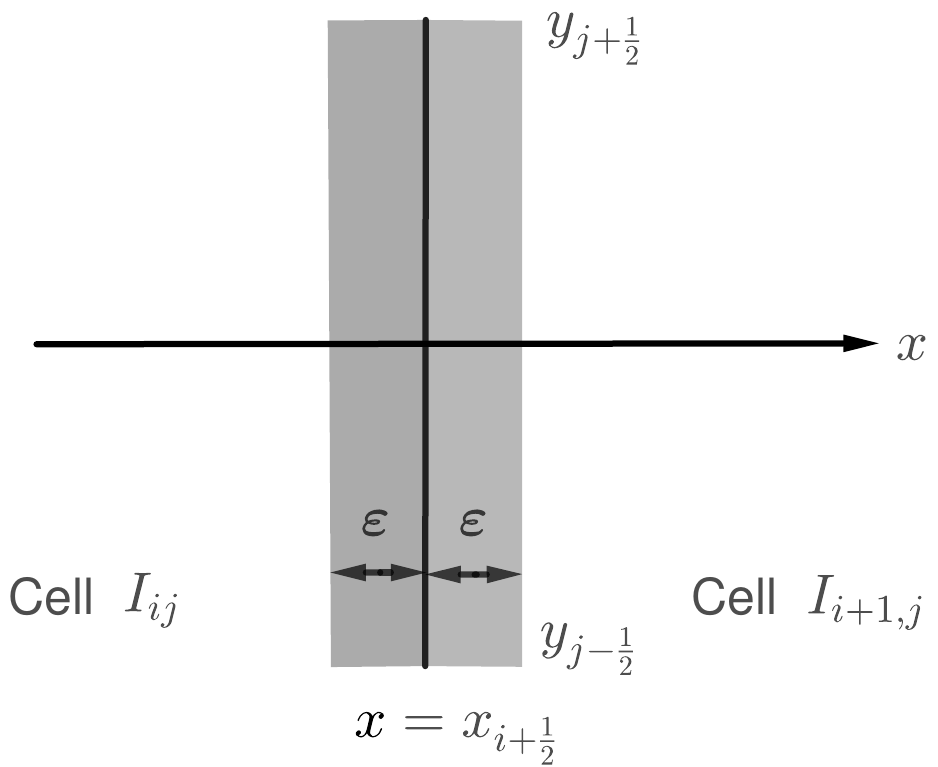}
	\caption{The sketch near the interface $x=x_{i+\frac{1}{2}}$: the left/right part of the gray area represents the local region of cell $I_{ij}/I_{i+1,j}$ near the interface with a thickness of $\varepsilon\to0^{+}$.
	}
	\label{Fig:interface}
\end{figure}

\begin{figure}[!htb]
	\centering
    \subfigcapskip=4pt
    \begin{tabular}{ccc}
    \subfigure[]{
		\label{[level.sub.1]}
		\includegraphics[width=0.32\linewidth]{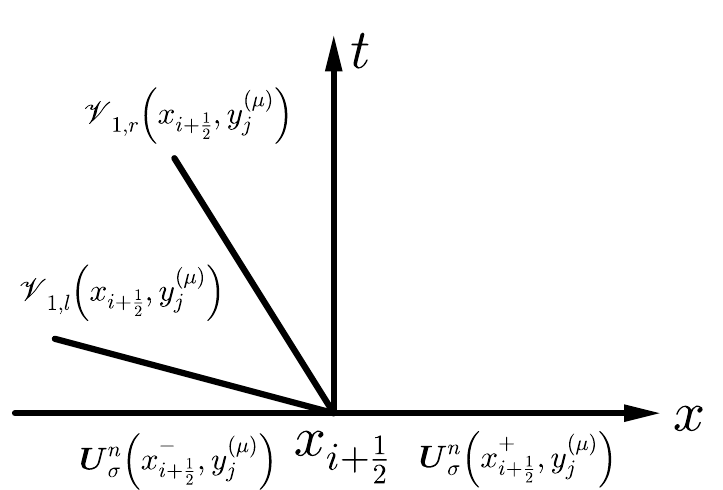}}
	\subfigure[]{
		\label{level.sub.2}
		\includegraphics[width=0.32\linewidth]{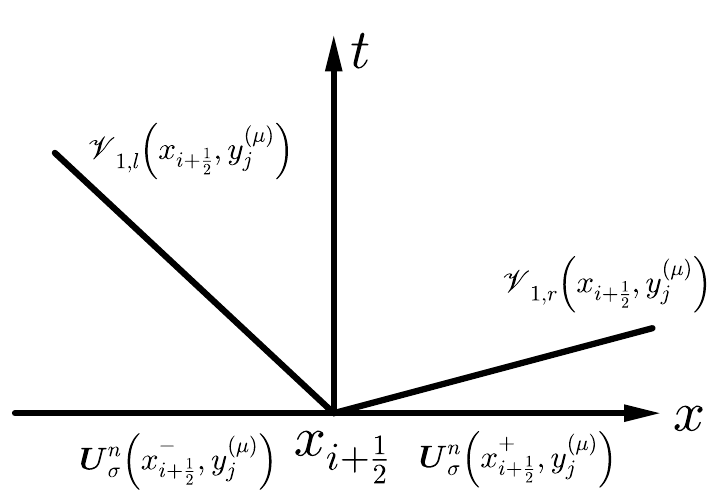}}
    \subfigure[]{
		\label{level.sub.3}
		\includegraphics[width=0.32\linewidth]{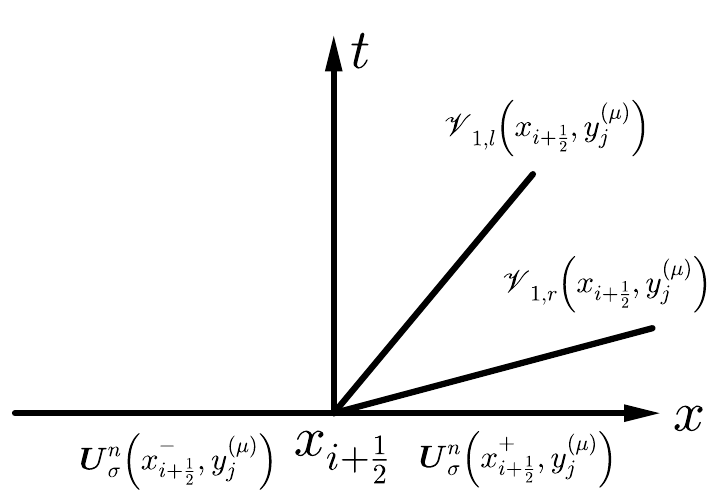}}
    \end{tabular}
	\caption{The three HLL wave patterns at the interface $x=x_{i+\frac{1}{2}}$ with the (estimated) minimum wave speed $\mathscr{V}_{1,l}\left(x_{i+\frac{1}{2}}, y_{j}^{(\mu)}\right)$ and the (estimated) maximum wave speed $\mathscr{V}_{1,r}\left(x_{i+\frac{1}{2}}, y_{j}^{(\mu)}\right)$.
    \label{Fig:wavepattern}
	}
\end{figure}

Take the interface $x=x_{i+\frac{1}{2}}$ as an example. See 
Figure \ref{Fig:interface}, where the left/right part of the gray area represents the local region of cell $I_{ij}/I_{i+1,j}$ near the interface $x=x_{i+\frac{1}{2}}$ with a thickness of $\varepsilon\to0^{+}$. 
The contribution of the Godunov–Powell source term due to the weak divergence at the interface $x=x_{i+\frac{1}{2}}$ can be approximated by  
\begin{equation}\label{Equ:sourceface}
	\begin{aligned}
	&	\quad\frac{1}{\Delta x \Delta y} \lim_{\varepsilon \to 0^{+}}
     \int_{x_{i+\frac{1}{2}}-\varepsilon}^{x_{i+\frac{1}{2}}+\varepsilon}
      \int_{y_{j-\frac{1}{2}}}^{y_{j+\frac{1}{2}}}
      \left(\frac{\partial \left(B_{1}\right)_{\sigma}^{n}}{\partial x}+\frac{\partial \left(B_{2}\right)_{\sigma}^{n}}{\partial y}\right)
      \boldsymbol{S}\left(\boldsymbol{U}_{\sigma}^{n}\right)
      \mathrm{~d}x\mathrm{d}y\\
      &=\frac{1}{\Delta x \Delta y} \lim_{\varepsilon \to 0^{+}}
     \int_{x_{i+\frac{1}{2}}-\varepsilon}^{x_{i+\frac{1}{2}}+\varepsilon}
      \int_{y_{j-\frac{1}{2}}}^{y_{j+\frac{1}{2}}}
      \frac{\partial \left(B_{1}\right)_{\sigma}^{n}}{\partial x}\boldsymbol{S}\left(\boldsymbol{U}_{\sigma}^{n}\right)
      \mathrm{~d}x\mathrm{d}y \\
      &\approx \frac{1}{\Delta x}\sum_{\mu=1}^{q}\omega_{\mu}^{G}\lim_{\varepsilon \to 0^{+}}\int_{x_{i+\frac{1}{2}}-\varepsilon}^{x_{i+\frac{1}{2}}+\varepsilon}
      \frac{\partial \left(B_{1}\right)_{\sigma}^{n}\left(x, y_{j}^{(\mu)}\right)}{\partial x}\boldsymbol{S}\left(\boldsymbol{U}_{\sigma}^{n}\left(x, y_{j}^{(\mu)}\right)\right)
      \mathrm{~d}x \\
      &\approx \frac{1}{\Delta x}\sum_{\mu=1}^{q}\omega_{\mu}^{G}\bigg[
      \boldsymbol{S}\left(\boldsymbol{U}_{\sigma}^{n}\big(x_{i+\frac{1}{2}}^-, y_{j}^{(\mu)}\big)\right)\omega_{\text{jump}}^{-}
      \left[\!\left[\left(B_{1}\right)_{\sigma}^{n}
      \right]\!\right]_{i+\frac12}(y_{j}^{(\mu)})+
      \\
      & \qquad \qquad \qquad ~~
       \boldsymbol{S}\left(\boldsymbol{U}_{\sigma}^{n}\big(x_{i+\frac{1}{2}}^+, y_{j}^{(\mu)}\big)\right)\omega_{\text{jump}}^{+}
      \left[\!\left[\left(B_{1}\right)_{\sigma}^{n}
      \right]\!\right]_{i+\frac12}(y_{j}^{(\mu)})
      \bigg],
	\end{aligned}
\end{equation}
where $\left[\!\left[\left(B_{1}\right)_{\sigma}^{n}
\right]\!\right]_{i+\frac12}(y)$ 
 is the jump across the interface $x=x_{i+\frac{1}{2}}$, and $\omega_{\text{jump}}^{-}$ (resp.~$\omega_{\text{jump}}^{+}$) represents the  weight assigned to the cell $I_{ij}$ (resp.~$I_{i+1,j}$). These two weights satisfy $\omega_{\text{jump}}^{-}+\omega_{\text{jump}}^{+}=1$, as we have discussed in \eqref{eq:split1}. 
In this work, we adopt the following wave-dependent weights $\omega_{\text{jump}}^{-}$ and $\omega_{\text{jump}}^{+}$: 
\begin{equation*}
	\omega_{\text{jump}}^{-}=
	\frac{- \mathscr{V}_{1,i+\frac{1}{2}}^{-}\left(y_{j}^{(\mu)}\right)}
	{\mathscr{V}_{1,i+\frac{1}{2}}^{+}\left(y_{j}^{(\mu)}\right)-
		\mathscr{V}_{1,i+\frac{1}{2}}^{-}\left(y_{j}^{(\mu)}\right)}, \qquad 
	\omega_{\text{jump}}^{+}=
	\frac{\mathscr{V}_{1,i+\frac{1}{2}}^{+}\left(y_{j}^{(\mu)}\right)}
	{\mathscr{V}_{1,i+\frac{1}{2}}^{+}\left(y_{j}^{(\mu)}\right)-
		\mathscr{V}_{1,i+\frac{1}{2}}^{-}\left(y_{j}^{(\mu)}\right)}.
\end{equation*}
Let $\mathscr{V}_{1,l}\left(x_{i+\frac{1}{2}}, y_{j}^{(\mu)}\right)$ and $\mathscr{V}_{1,r}\left(x_{i+\frac{1}{2}}, y_{j}^{(\mu)}\right)$ denote the (estimated) minimum and maximum wave speeds at the interface $x=x_{i+\frac{1}{2}}$, respectively. Recall that  
\begin{equation*}
	\mathscr{V}_{1,i+\frac{1}{2}}^{-}\left(y_{j}^{(\mu)}\right)=
\min\left\{0,\mathscr{V}_{1,l}\left(x_{i+\frac{1}{2}}, y_{j}^{(\mu)}\right)\right\}, \quad \mathscr{V}_{1,i+\frac{1}{2}}^{+}\left(y_{j}^{(\mu)}\right)=
\max\left\{0,\mathscr{V}_{1,r}\left(x_{i+\frac{1}{2}}, y_{j}^{(\mu)}\right)\right\}.
\end{equation*}
Figure \ref{Fig:wavepattern} shows the three HLL wave patterns at the interface $x=x_{i+\frac{1}{2}}$.
The corresponding weights are listed below:
\begin{itemize}
  \item For case (a), $\omega_{\text{jump}}^{-}=
	\frac{- \mathscr{V}_{1,i+\frac{1}{2}}^{-}\left(y_{j}^{(\mu)}\right)}
	{\mathscr{V}_{1,i+\frac{1}{2}}^{+}\left(y_{j}^{(\mu)}\right)-
		\mathscr{V}_{1,i+\frac{1}{2}}^{-}\left(y_{j}^{(\mu)}\right)}$, $\omega_{\text{jump}}^{+}=0$;
  \item For case (b), $\omega_{\text{jump}}^{-}=
	\frac{- \mathscr{V}_{1,i+\frac{1}{2}}^{-}\left(y_{j}^{(\mu)}\right)}
	{\mathscr{V}_{1,i+\frac{1}{2}}^{+}\left(y_{j}^{(\mu)}\right)-
		\mathscr{V}_{1,i+\frac{1}{2}}^{-}\left(y_{j}^{(\mu)}\right)}$, $\omega_{\text{jump}}^{+}=
	\frac{ \mathscr{V}_{1,i+\frac{1}{2}}^{+}\left(y_{j}^{(\mu)}\right)}
	{\mathscr{V}_{1,i+\frac{1}{2}}^{+}\left(y_{j}^{(\mu)}\right)-
		\mathscr{V}_{1,i+\frac{1}{2}}^{-}\left(y_{j}^{(\mu)}\right)}$;
  \item For case (c), $\omega_{\text{jump}}^{-}=0$, $\omega_{\text{jump}}^{+}=
	\frac{ \mathscr{V}_{1,i+\frac{1}{2}}^{+}\left(y_{j}^{(\mu)}\right)}
	{\mathscr{V}_{1,i+\frac{1}{2}}^{+}\left(y_{j}^{(\mu)}\right)-
		\mathscr{V}_{1,i+\frac{1}{2}}^{-}\left(y_{j}^{(\mu)}\right)}$.
\end{itemize}
Such choices reflect an ``upwind'' discretization of the Godunov--Powell source term. It is worth mentioning that this discretization is actually motivated by rigorous PP analysis, which will become clear in the proof of Theorem \ref{Thm:2DPPpropertytheorem} in section \ref{sec:PPproof}.

\subsection{Convex decomposition of cell averages}

In the framework \cite{ZHANG20103091,CUI2023111882,Cui2024On} of constructing provably PP high-order schemes, a key ingredient 
is the decomposition of the cell averages of the numerical solution into
a convex combination of the point values at certain quadrature points.
The convex decomposition determines not only the theoretical PP CFL condition of the resulting scheme, but also identifies the points at which to perform the PP limiter, as will be shown in section \ref{Sec:PPlimiter}.

Recall that $\{x_{i}^{(\mu)}\}_{\mu=1}^{q}$ and $\{y_{j}^{(\mu)}\}_{\mu=1}^{q}$ denote the Gauss quadrature nodes, and $\{\omega_{\mu}^{G}\}_{\mu=1}^{q}$ are the associated weights in the DG schemes in \eqref{Equ:2DsemidiscreR}. Let $a_{1}$ and $a_{2}$ represent the (estimated) maximum wave speeds in the $x$- and $y$-directions, respectively. According to \cite{CUI2023111882,Cui2024On}, a general feasible convex decomposition for polynomial space $\mathbb{P}^k(I_{ij})$ on a rectangular cell can be expressed as 
\begin{equation}\label{Equ:2DgeneralCD}
	\begin{aligned}
		\bar{\boldsymbol{U}}_{ij}=&\sum_{\mu=1}^q \omega_\mu^{G}\left[\omega_1^{+} \boldsymbol{U}_{ij}\left(x_{i+\frac{1}{2}}^{-}, y_{j}^{(\mu)}\right)+\omega_1^{-}
		\boldsymbol{U}_{ij}\left(x_{i-\frac{1}{2}}^{+}, y_{j}^{(\mu)}\right)+\omega_2^{+} \boldsymbol{U}_{ij}\left(x_{i}^{(\mu)}, y_{j+\frac{1}{2}}^{-}\right)+\omega_2^{-}
		\boldsymbol{U}_{ij}\left(x_{i}^{(\mu)}, y_{j-\frac{1}{2}}^{+}\right)
		\right]\\
		&+\sum_{s=1}^S \omega_s \boldsymbol{U}_{ij}\left(x_{ij}^{(s)}, y_{ij}^{(s)}\right),
	\end{aligned}
\end{equation}
which should simultaneously satisfy the following three conditions:
\begin{description}
	\item[(i)] the convex decomposition holds exactly for all $\boldsymbol{U}_{ij}\left(x, y\right)\in [\mathbb{P}^k\left(I_{ij}\right)]^8$;
	\item[(ii)] the weights $\left\{\omega_{1}^{\pm}, \omega_{2}^{\pm}, \omega_{s}\right\}$ are all positive and their summation equals one;
	\item[(iii)] the internal node set $\mathbb{S}_{ij}^{\text{internal}}:={\left\{\left(x_{ij}^{(s)}, y_{ij}^{(s)}\right)\right\}}_{s=1}^{S}\subset I_{ij}$.
\end{description}
In convex decomposition (\ref{Equ:2DgeneralCD}), the boundary node set is given by
\begin{equation*}\label{Equ:boundarynodeset}
	\mathbb{S}_{ij}^{\text{boundary}}:={\left\{\left(x_{i+\frac{1}{2}}^{-}, y_{j}^{(\mu)}\right)\right\}}_{\mu=1}^{q}\bigcup{\left\{\left(x_{i-\frac{1}{2}}^{+}, y_{j}^{(\mu)}\right)\right\}}_{\mu=1}^{q}\bigcup{\left\{\left(x_{i}^{(\mu)},
		y_{j+\frac{1}{2}}^{-} \right)\right\}}_{\mu=1}^{q}\bigcup{\left\{\left(x_{i}^{(\mu)},
		y_{j-\frac{1}{2}}^{+} \right)\right\}}_{\mu=1}^{q}. 
\end{equation*}
Based on \eqref{Equ:2DgeneralCD}, the cell average $\bar{\boldsymbol{U}}_{ij}$ is decomposed into a convex combination of solution values at the
quadrature points belonging to  $$\mathbb{S}_{ij}:=\mathbb{S}_{ij}^{\text{boundary}}\bigcup \mathbb{S}_{ij}^{\text{internal}},$$
at which the PP limiter will be performed; see section \ref{Sec:PPlimiter}.

Denote $\left\{\hat{x}_{i}^{(\mu)}\right\}_{\mu=1}^{L}$ and $\left\{\hat{y}_{j}^{(\mu)}\right\}_{\mu=1}^{L}$
as the $L$-point Gauss--Lobatto quadrature nodes with $L=\lceil\frac{k+3}{2}\rceil$ in $\left[x_{i-\frac{1}{2}}, x_{i+\frac{1}{2}}\right]$ and $\left[y_{j-\frac{1}{2}}, y_{j+\frac{1}{2}}\right]$, respectively, and $\left\{\hat{\omega}_{\mu}^{GL}\right\}_{\mu=1}^{L}$ are the associated weights for the interval $[-\frac12,\frac12]$ with  $\hat{\omega}_{1}^{GL}=\hat{\omega}_{L}^{GL}=\frac{1}{L(L-1)}$. 
Based on the tensor product of the $L$-point Gauss--Lobatto quadrature 
and the $q$-point Gauss quadrature,
Zhang and Shu \cite{ZHANG20103091} proposed  the following classic convex decomposition. 

\noindent
\textbf{\emph{Zhang--Shu convex decomposition}:}
\begin{equation}\label{Equ:2DZhangShuCD}
	\begin{aligned}
		\bar{\boldsymbol{U}}_{ij}=&\sum_{\mu=1}^q \omega_\mu^{G}\hat{\omega}_1^{GL}\left[\kappa_1 \boldsymbol{U}_{ij}\left(x_{i+\frac{1}{2}}^{-}, y_{j}^{(\mu)}\right)+\kappa_1
		\boldsymbol{U}_{ij}\left(x_{i-\frac{1}{2}}^{+}, y_{j}^{(\mu)}\right)\right.\\
		&\left.+\kappa_2 \boldsymbol{U}_{ij}\left(x_{i}^{(\mu)}, y_{j+\frac{1}{2}}^{-}, \right)+\kappa_2
		\boldsymbol{U}_{ij}\left(x_{i}^{(\mu)}, y_{j-\frac{1}{2}}^{+}\right)
		\right]\\
		&+\sum_{s=2}^{L-1} \sum_{\mu=1}^{q}\hat{\omega}_s^{GL} \omega_\mu^{G}\left[\kappa_1\boldsymbol{U}_{ij}\left(\hat{x}_{i}^{(s)}, y_{j}^{(\mu)}\right)+\kappa_2\boldsymbol{U}_{ij}\left(x_{i}^{(\mu)}, \hat{y}_{j}^{(s)}\right)\right]
	\end{aligned}
\end{equation}
with
\begin{equation*}\label{Equ:2DZhangShuk1k2CD}
	\kappa_1:=\frac{\frac{a_1}{\Delta x}}{\frac{a_1}{\Delta x}+\frac{a_2}{\Delta y}}, \quad \kappa_2:=\frac{\frac{a_2}{\Delta y}}{\frac{a_1}{\Delta x}+\frac{a_2}{\Delta y}}.
\end{equation*}
Cui, Ding, and Wu \cite{CUI2023111882,Cui2024On} discovered an optimal convex decomposition involving much fewer nodes and allowing the mildest PP time step-size, as shown below for the cases of $\mathbb P^2$ and $\mathbb P^3$ spaces. For the optimal convex decomposition for general $\mathbb P^k$ spaces with $k\ge 4$, see \cite{Cui2024On}. 

\noindent
\textbf{\emph{Optimal convex decomposition for $\mathbb P^2$ and $\mathbb P^3$}:}
\begin{equation}\label{Equ:2DoptimalCD}
	\begin{aligned}
		\bar{\boldsymbol{U}}_{ij}=&\frac{\kappa_{1}}{2}\sum_{\mu=1}^q \omega_\mu^{G}\left[ \boldsymbol{U}_{ij}\left(x_{i+\frac{1}{2}}^{-}, y_{j}^{(\mu)}\right)+
		\boldsymbol{U}_{ij}\left(x_{i-\frac{1}{2}}^{+}, y_{j}^{(\mu)}\right)\right]\\
		&+\frac{\kappa_{2}}{2}\sum_{\mu=1}^q\omega_\mu^{G}\left[ \boldsymbol{U}_{ij}\left(x_{i}^{(\mu)}, y_{j+\frac{1}{2}}^{-} \right)+
		\boldsymbol{U}_{ij}\left(x_{i}^{(\mu)}, y_{j-\frac{1}{2}}^{+}\right)
		\right]\\
		&+\omega\sum_{s} \boldsymbol{U}_{ij}\left(\hat{x}_{s}, \hat{y}_{s}\right)
	\end{aligned}
\end{equation}
with
\begin{equation*}\label{Equ:2DoptinternodeCD}
	\left\{\left(\hat{x}_s, \hat{y}_s\right)\right\}=
	\left\{
	\begin{aligned}
		&\left(x_i, y_j \pm \frac{\Delta y}{2 \sqrt{3}} \sqrt{\frac{\phi_*-\phi_2}{\phi_*}}\right), &\text { if } \phi_1 \geq \phi_2, \\
		&\left(x_i \pm \frac{\Delta x}{2 \sqrt{3}} \sqrt{\frac{\phi_*-\phi_1}{\phi_*}}, y_j\right), &\text { if } \phi_1<\phi_2,
	\end{aligned}
	\right.
\end{equation*}
where
\begin{equation*}
	\phi_1=\frac{a_1}{\Delta x}, \quad \phi_2=\frac{a_2}{\Delta y}, \quad \phi_*=\max \left\{\phi_1, \phi_2\right\}, \quad \psi=\phi_1+\phi_2+2 \phi_*, \quad \kappa_1=\frac{\phi_1}{\psi}, \quad \kappa_2=\frac{\phi_2}{\psi}, \quad \omega=\frac{\phi_*}{\psi} .
\end{equation*}

\subsection{Positivity--preserving limiter}\label{Sec:PPlimiter}

Let $\boldsymbol{U}_{ij}^{n}\left(\boldsymbol{x}\right)$ denote the OEDG solution  in $I_{ij}$ at the $n$th time level after the LDF OE procedure, and let  $\bar{\boldsymbol{U}}_{ij}^{n}$ be its cell average over $I_{ij}$. 
As it will be shown in section \ref{sec:PPproof}, the OEDG method with the PP HLL flux and upwind discrete Godunov--Powell source term preserves the positivity for the updated cell averages under a CFL condition, if the OEDG solution satisfies 
\begin{equation}\label{Equ:PPGL1D}
	{\boldsymbol{U}}_{ij}^n\left(\boldsymbol{x}\right) \in \mathcal{G} \quad \forall \boldsymbol{x} \in \mathbb{S}_{ij}
\end{equation}
where $\mathbb{S}_{ij}$ is the set of all nodes involved in the adopted convex decomposition (\ref{Equ:2DgeneralCD}).

The LDF OE procedure, although enhancing the numerical stability for strong shocks, does not necessarily guarantee the condition \eqref{Equ:PPGL1D}, which should be enforced by a simple PP limiter. 
For readers' convenience, here we briefly review the PP limiter \cite{ZHANG20108918,CHENG2013255,Wu2019Provably}. 
The PP limiting procedure mainly consists of two steps.
First, modify the density to enforce the positivity by
\begin{equation}\label{Equ:PPlimistep1}
	\widehat{\rho}_{ij}\left(\boldsymbol{x}\right)=\theta_1\left(\rho_{ij}^n\left(\boldsymbol{x}\right)-\bar{\rho}_{ij}^n\right)+
	\bar{\rho}_{ij}^n  \quad \text{with} \quad \theta_1=\min \left\{1, \frac{\bar{\rho}_{ij}^n-\epsilon_1}{\bar{\rho}_{ij}^n-\min _{\boldsymbol{x} \in \mathbb{S}_{ij}} \rho_{ij}^n\left(\boldsymbol{x}\right)}\right\} .
\end{equation}
Then, modify $\widehat{\boldsymbol{U}}_{ij}\left(\boldsymbol{x}\right):=
\left(\widehat{\rho}_{ij}\left(\boldsymbol{x}\right), \boldsymbol{m}_{ij}^n\left(\boldsymbol{x}\right), \boldsymbol{B}_{ij}^n\left(\boldsymbol{x}\right), E_{ij}^n\left(\boldsymbol{x}\right)\right)^{\top}$ to ensure positive internal energy by
\begin{equation}\label{Equ:PPlimistep2}
	\widetilde{\boldsymbol{U}}_{ij}^n\left(\boldsymbol{x}\right)=\theta_2
	\left(\widehat{\boldsymbol{U}}_{ij}\left(\boldsymbol{x}\right)-\bar{\boldsymbol{U}}_{ij}^n\right)+
	\bar{\boldsymbol{U}}_{ij}^n \quad \text{with} \quad \theta_2=\min \left\{1, \frac{\mathcal{E}\left(\bar{\boldsymbol{U}}_{ij}^n\right)-\epsilon_2}
	{\mathcal{E}\left(\bar{\boldsymbol{U}}_{ij}^n\right)-\min _{\boldsymbol{x} \in \mathbb{S}_{ij}} \mathcal{E}\left(\widehat{\boldsymbol{U}}_{ij}\left(\boldsymbol{x}\right)\right)}\right\}.
\end{equation}
Here $\epsilon_1=\min \left\{10^{-13}, \bar{\rho}_{ij}^n\right\}$ and $\epsilon_2=\min \left\{10^{-13}, \mathcal{E}\left(\bar{\boldsymbol{U}}_{ij}^n\right)\right\}$ are introduced to avoid the impact of rounding errors on PP property \cite{Wu2019Provably}. 
Note that the modified solution has the same cell average $\bar{\boldsymbol{U}}_{ij}^{n}$ on $I_{ij}$ and holds
\begin{equation}\label{Equ:PPlimihold}
	\widetilde{\boldsymbol{U}}_{ij}^n\left(\boldsymbol{x}\right) \in \mathcal{G}_\epsilon=
	\left\{\boldsymbol{U}: \rho \geq \epsilon_1, \mathcal{E}\left(\boldsymbol{U}\right) \geq \epsilon_2\right\} \quad \forall \boldsymbol{x} \in \mathbb{S}_{ij}.
\end{equation}
It is worth mentioning that such a limiter does not reduce the high-order accuracy of the DG schemes \cite{ZHANG20103091,ZHANG20108918}. The PP limiter is performed right after each LDF OE procedure, and it maintains the LDF property for the magnetic field.

\subsection{Algorithmic flow of PP LDF OEDG schemes}\label{Sec:PPscheme}

Based on the LDF OEDG method, the proper upwind discretization of Godunov--Powell source term, and the PP limiter, we construct the PP LDF OEDG schemes as follows. 
For clarity, we focus on the detailed descriptions of the schemes with the forward Euler time discretization. 

\begin{description}
	\item[Step 0] 
	Initialize the approximate solution  $\boldsymbol{U}_{h}^{0}\left(\boldsymbol{x}\right)$ 
	by the $L^{2}$-projection of initial data of $\boldsymbol{R}$ onto the space $\mathbb{V}_{R}^{k}$ and the $L^{2}$-projection of initial data of $\boldsymbol{Q}$ onto the LDF space $\mathbb{V}_{Q}^{k}$. 
	The cell average $\bar{\boldsymbol{U}}_{ij}^{0}\in \mathcal{G}$, because of the convexity of $\mathcal{G}$.

	\item[Step 1] Given the approximate solution $\boldsymbol{U}_{h}^{n}\left(\boldsymbol{x}\right)$ with 
	admissible cell averages $\left\{\bar{\boldsymbol{U}}_{ij}^{n}\right\}$, 
	perform LDF OE procedure 
	\begin{equation}\label{Equ:LDFOEforQandR}
		\boldsymbol{R}_{\sigma}^{n}\left(\boldsymbol{x}\right)=	\mathcal{F}_{\text{OE}}^R \boldsymbol{R}_{h}^{n}\left(\boldsymbol{x}\right), \quad
		\boldsymbol{Q}_{\sigma}^{n}\left(\boldsymbol{x}\right)=	\mathcal{F}_{\text{OE}}^Q \boldsymbol{Q}_{h}^{n}\left(\boldsymbol{x}\right),
	\end{equation}
	to obtain the OEDG solution $\boldsymbol{U}_{\sigma}^{n}\left(\boldsymbol{x}\right)$. 
	
	\item[Step 2] Given the OEDG solution $\boldsymbol{U}_{\sigma}^{n}\left(\boldsymbol{x}\right)$
	with the admissible cell averages $\left\{\bar{\boldsymbol{U}}_{ij}^{n}\right\}$, define the polynomial vector of $\boldsymbol{U}_{\sigma}^{n}\left(\boldsymbol{x}\right)$ on $I_{ij}$ as 
	$\boldsymbol{U}_{ij}^{n}\left(\boldsymbol{x}\right)$, and perform PP limiter (section \ref{Sec:PPlimiter}) to modify $
	\left\{\boldsymbol{U}_{ij}^{n}\left(\boldsymbol{x}\right)\right\}$ such that
	the modified polynomial vector  $\left\{\widetilde{\boldsymbol{U}}_{ij}^{n}\left(\boldsymbol{x}\right)\right\}$ satisfies 
	\begin{equation}\label{Equ:2DPPlimitercond}
		\widetilde{\boldsymbol{U}}_{ij}^{n}\left(\boldsymbol{x}\right)\in\mathcal{G}
		\quad \forall \boldsymbol{x}\in \mathbb{S}_{ij},
	\end{equation}
	where $\mathbb{S}_{ij}$ is the set of all the nodes of the adopted convex decomposition (\ref{Equ:2DgeneralCD}). We opt to use the optimal convex decomposition (\ref{Equ:2DoptimalCD}) in our computations. 
	\item[Step 3] Let $\widetilde{\boldsymbol{U}}_{\sigma}^{n}\left(\boldsymbol{x}\right)$ denote the piecewise polynomial solution defined by $\widetilde{\boldsymbol{U}}_{ij}^{n}\left(\boldsymbol{x}\right)$. 
	Perform an SSP Runge--Kutta time discretization to evolve the OEDG solution. 
	For clarity, we only provide specific details for the forward Euler method. 
	The cell averages $\{\bar{\boldsymbol{U}}_{ij}^{n+1}\}$ are updated by utilizing the following scheme
	\begin{equation}\label{Equ:updatecellave}
		\bar{\boldsymbol{U}}_{ij}^{n+1}=\bar{\boldsymbol{U}}_{ij}^{n}+\Delta t_{n}
		\boldsymbol{L}_{ij}\left(\widetilde{\boldsymbol{U}}_{\sigma}^{n}\right),
	\end{equation}
	where the operator $\boldsymbol{L}_{ij}$ is defined by \eqref{Equ:averageMHDsourceL} with the properly discretized Godunov--Powell source term. 
     The high-order ``moments" of the DG solution are updated by evolving the semi-discrete schemes (\ref{Equ:2DsemidiscreR}) and (\ref{Equ:2DsemidiscreQ}) with the forward Euler method. 
	Note that the scheme (\ref{Equ:updatecellave})  ensures  $\bar{\boldsymbol{U}}_{ij}^{n+1}\in\mathcal{G}$, as shown in Theorem \ref{Thm:2DPPpropertytheorem} later.
	\item[Step 4] If $t_{n+1}<T$, assign $n\leftarrow n+1$ and go to Step 1; otherwise, output numerical results.
\end{description}

\subsection{Rigorous PP analysis}\label{sec:PPproof}

To prove $\bar{\boldsymbol{U}}_{ij}^{n+1}\in\mathcal{G}$ in scheme (\ref{Equ:updatecellave}), we first introduce several key auxiliary lemmas related to the GQL approach \cite{Wu2018Positivity,Wu2023Geometric}, which will be used in the proof of $\bar{\boldsymbol{U}}_{ij}^{n+1}\in\mathcal{G}$ in Theorem \ref{Thm:2DPPpropertytheorem}.

\begin{lemma}[GQL representation \cite{Wu2018Positivity}]\label{Lem:admistatequiv}
	The admissible state set $\mathcal{G}$ is exactly equivalent to
	\begin{equation}\label{Equ:admistatequiv}
		\mathcal{G}_*=\left\{\boldsymbol{U}=\left(\rho, \boldsymbol{m}, \boldsymbol{B}, E\right)^{\top}:~ \rho>0,~ \boldsymbol{U} \cdot \boldsymbol{n}^*+\frac{\left\|\boldsymbol{B}^*\right\|^2}{2}>0~~ \forall \boldsymbol{v}^*, \boldsymbol{B}^* \in \mathbb{R}^3\right\},
	\end{equation}
	where $\boldsymbol{v}^{*}$ and $\boldsymbol{B}^{*}$ are the extra free auxiliary variables independent of $\boldsymbol{U}$, and 
	$$
	\boldsymbol{n}^*=\left(\frac{\left\|\boldsymbol{v}^*\right\|^2}{2},~-\boldsymbol{v}^*,~
	-\boldsymbol{B}^*,~ 1\right)^{\top}.
	$$
\end{lemma}
In the original definition \eqref{Equ:phyadmiset} of $\mathcal{G}$, the second constraint is nonlinear respect to $\boldsymbol{U}$, leading to challenges in analyzing the positivity of $\mathcal{E}(\boldsymbol{U})$. Thanks to the GQL representation \eqref{Equ:admistatequiv}, all equivalent constraints in $\mathcal{G}_*$ become linear with respect to $\boldsymbol{U}$, thereby greatly facilitating the PP analysis. The first proof of Lemma \ref{Lem:admistatequiv} can be found in \cite{Wu2018Positivity}. For more details about its geometric interpretation, refer to the GQL framework in \cite{Wu2023Geometric}.

\begin{lemma}[\cite{Wu2018A}]\label{Lem:sourcetermequiv}
	For any $\boldsymbol{U}\in\mathcal{G}$ and any $\boldsymbol{v}^{*}, \boldsymbol{B}^{*}\in \mathbb{R}^{3}$, we have
	\begin{equation*}\label{Equ:sourceequiv1}
		\boldsymbol{S}\left(\boldsymbol{U}\right) \cdot \boldsymbol{n}^*=\left(\boldsymbol{v}-\boldsymbol{v}^*\right) \cdot\left(\boldsymbol{B}-\boldsymbol{B}^*\right)-\boldsymbol{v}^* \cdot \boldsymbol{B}^*,
	\end{equation*}
	\begin{equation*}\label{Equ:sourceequiv2}
		\left|\sqrt{\rho}\left(\boldsymbol{v}-\boldsymbol{v}^*\right) \cdot\left(\boldsymbol{B}-\boldsymbol{B}^*\right)\right|<\boldsymbol{U} \cdot \boldsymbol{n}^*+\frac{\left\|\boldsymbol{B}^*\right\|^2}{2}.
	\end{equation*}
	Furthermore, for any $b\in\mathbb{R}$, it holds
	\begin{equation}\label{Equ:sourceinequiv1}
		-b\left(\boldsymbol{S}\left(\boldsymbol{U}\right) \cdot \boldsymbol{n}^*\right) \geq b\left(\boldsymbol{v}^* \cdot \boldsymbol{B}^*\right)-\frac{|b|}{\sqrt{\rho}}\left(\boldsymbol{U} \cdot \boldsymbol{n}^*+\frac{\left\|\boldsymbol{B}^*\right\|^2}{2}\right).
	\end{equation}
\end{lemma}

The proof of Lemma \ref{Lem:sourcetermequiv} can be found in \cite{Wu2018A}.

\begin{lemma}\label{Thm:2DPPstatetheorem}
	Given two admissible states $\boldsymbol{U}^{(1)}$ and $\boldsymbol{U}^{(2)}$, we define for $\ell\in \{1,2\}$ that 
	\begin{align*}
		&\hat{\alpha}_{\ell}^{(1)}:=\max\left\{
		\left(u_{\ell}\right)^{(1)},
		\frac{\sqrt{\rho^{(1)}}\left(u_{\ell}\right)^{(1)}+\sqrt{\rho^{(2)}}\left(u_{\ell}\right)^{(2)}}
		{\sqrt{\rho^{(1)}}+\sqrt{\rho^{(2)}}}\right\}+
		\mathscr{C}_{\ell}\left(\boldsymbol{U}^{(1)}\right)+
		\frac{\left\|\boldsymbol{B}^{(1)}-\boldsymbol{B}^{(2)}\right\|}
		{\sqrt{\rho^{(1)}}+\sqrt{\rho^{(2)}}},
		\\
		&\hat{\alpha}_{\ell}^{(2)}:=\max\left\{
		-\left(u_{\ell}\right)^{(2)},
		-\left(\frac{\sqrt{\rho^{(2)}}\left(u_{\ell}\right)^{(2)}+\sqrt{\rho^{(1)}}\left(u_{\ell}\right)^{(1)}}
		{\sqrt{\rho^{(2)}}+\sqrt{\rho^{(1)}}}\right)\right\}+
		\mathscr{C}_{\ell}\left(\boldsymbol{U}^{(2)}\right)+
		\frac{\left\|\boldsymbol{B}^{(2)}-\boldsymbol{B}^{(1)}\right\|}
		{\sqrt{\rho^{(2)}}+\sqrt{\rho^{(1)}}}.
	\end{align*}
	Then for any $\alpha_{\ell}^{(j)}\geq\hat{\alpha}_{\ell}^{(j)}$, $j=1, 2$, the states
	$$
	\begin{aligned}
		&\bar{\boldsymbol{U}}_{\ell}:=\frac{1}{\alpha_{\ell}^{(1)}+\alpha_{\ell}^{(2)}}
		\left[\left(\alpha_{\ell}^{(1)}\boldsymbol{U}^{(1)}-\boldsymbol{F}_{\ell}
		\left(\boldsymbol{U}^{(1)}\right)\right)+\left(\alpha_{\ell}^{(2)}
		\boldsymbol{U}^{(2)}+\boldsymbol{F}_{\ell}\left(\boldsymbol{U}^{(2)} \right)\right)\right], \quad \ell =1,2, 
	\end{aligned}
	$$
	belong to $\mathcal{G}_{\rho}:=\left\{\boldsymbol{U}=\left(\rho, \boldsymbol{m}, \boldsymbol{B}, E\right)^{\top}: \rho>0\right\}$, and satisfy
	$$
	\bar{\boldsymbol{U}}_{\ell}\cdot\boldsymbol{n}^{*}+\frac{\left|
		\boldsymbol{B}^{*}\right|^{2}}{2}\geq
	\frac{-\boldsymbol{v}^{*}\cdot\boldsymbol{B}^{*}}
	{\alpha_{\ell}^{(1)}+\alpha_{\ell}^{(2)}}\left(
	\left(B_{\ell}\right)^{(1)}-\left(B_{\ell}\right)^{(2)}\right)\quad
	\forall\boldsymbol{v}^{*},\boldsymbol{B}^{*}\in \mathbb{R}^{3}.
	$$
\end{lemma}


\noindent \textbf{\emph{Proof }} 
This directly follows from Theorem 1 in \cite{Wu2019Provably} with $N=2$, by taking
$s_{1}=s_{2}=1$, 
$\boldsymbol{\xi}^{(1)}=-\boldsymbol{\xi}^{(2)}=(1,0)$ for $\ell=1$, and 
$\boldsymbol{\xi}^{(1)}=-\boldsymbol{\xi}^{(2)}=(0,1)$ for $\ell=2.$   \qed

\begin{lemma}\label{Thm:2DPPfluxtheorem}
	Assume $\boldsymbol{U}^{-}, \boldsymbol{U}^{+} \in \mathcal{G}$.
	If the approximate wave speeds in the HLL flux satisfy
	\begin{equation}\label{Equ:2DPPwavespeedcon}
		\mathscr{V}_{\ell,r}\geq\alpha_{\ell,r}\left(\boldsymbol{U}^{+}, \boldsymbol{U}^{-}\right), \quad
		\mathscr{V}_{\ell,l}\leq\alpha_{\ell,l}\left(\boldsymbol{U}^{-}, \boldsymbol{U}^{+}\right), \quad
		\ell=1,2,
	\end{equation}
	then
	\begin{equation}\label{Equ:2DPPfluxequal1}
		\widehat{\boldsymbol{F}_{\ell}}\left(\boldsymbol{U}^{-}, \boldsymbol{U}^{+}\right)=
		\mathscr{V}_{\ell}^{-}\boldsymbol{H}_{\ell}\left(\boldsymbol{U}^{-}, \boldsymbol{U}^{+}\right)+\boldsymbol{F}_{\ell}\left(\boldsymbol{U}^{-}\right)-
		\mathscr{V}_{\ell}^{-}\boldsymbol{U}^{-},
	\end{equation}
	\begin{equation}\label{Equ:2DPPfluxequal2}
		\widehat{\boldsymbol{F}_{\ell}}\left(\boldsymbol{U}^{-}, \boldsymbol{U}^{+}\right)=
		\mathscr{V}_{\ell}^{+}\boldsymbol{H}_{\ell}\left(\boldsymbol{U}^{-}, \boldsymbol{U}^{+}\right)+\boldsymbol{F}_{\ell}\left(\boldsymbol{U}^{+}\right)-
		\mathscr{V}_{\ell}^{+}\boldsymbol{U}^{+}
	\end{equation}
	with the intermediate state
	\begin{equation*}\label{Equ:2DPPinterme}
		\boldsymbol{H}_{\ell}\left(\boldsymbol{U}^{-}, \boldsymbol{U}^{+}\right):=
		\frac{1}{\mathscr{V}_{\ell}^{+}-\mathscr{V}_{\ell}^{-}}\left(
		\mathscr{V}_{\ell}^{+}\boldsymbol{U}^{+}-\boldsymbol{F}_{\ell}\left(\boldsymbol{U}^{+}\right)
		-\mathscr{V}_{\ell}^{-}\boldsymbol{U}^{-}+\boldsymbol{F}_{\ell}\left(\boldsymbol{U}^{-}\right)
		\right)
	\end{equation*}
	belonging to $\mathcal{G}$ and satisfying 
	\begin{equation}\label{Equ:2DPPHinequal}
		\boldsymbol{H}_{\ell}\left(\boldsymbol{U}^{-}, \boldsymbol{U}^{+}\right) \cdot \boldsymbol{n}^*+\frac{\left\|\boldsymbol{B}^*\right\|^2}{2}\geq
		-\frac{\boldsymbol{v}^* \cdot \boldsymbol{B}^*}{\mathscr{V}_{\ell}^{+}-\mathscr{V}_{\ell}^{-}}\left( \boldsymbol{B}_{\ell}^{+}- \boldsymbol{B}_{\ell}^{-}\right) \quad \forall \boldsymbol{v}^*, \boldsymbol{B}^* \in \mathbb{R}^3.
	\end{equation}
\end{lemma}

\noindent \textbf{\emph{Proof }}
This directly follows from Theorem 2 in \cite{Wu2019Provably} by taking
$\boldsymbol{\xi}=(1,0)$ for $\ell=1$ and $\boldsymbol{\xi}=(0,1)$ for $\ell=2$. \qed

With the above auxiliary lemmas as the key cornerstone, we are now ready to prove $\bar{\boldsymbol{U}}_{ij}^{n+1}\in\mathcal{G}$ in scheme (\ref{Equ:updatecellave}).
Recall that $\widetilde{\boldsymbol{U}}_\sigma^{n}$ denotes the OEDG solution after the LDF OE and PP limiting procedures. 
To shorten our notations, define 
$$
\boldsymbol{U}_{i+\frac{1}{2},j}^{\pm, \mu}:=\widetilde{\boldsymbol{U}}_\sigma^{n}\left(x_{i+\frac{1}{2}}^{\pm}, y_{j}^{(\mu)}\right), \quad \boldsymbol{U}_{i, j+\frac{1}{2}}^{ \mu, \pm}:=\widetilde{\boldsymbol{U}}_\sigma^{n}\left(x_{i}^{(\mu)}, y_{j+\frac{1}{2}}^{\pm}\right),
$$
$$
\mathscr{V}_{1, i+\frac{1}{2}}^{\mu, \pm}:=
\mathscr{V}_{1}^{\pm}\left(
\widetilde{\boldsymbol{U}}_\sigma^{n}\left(x_{i+\frac{1}{2}}^{-}, y_{j}^{(\mu)}\right),
\widetilde{\boldsymbol{U}}_\sigma^{n}\left(x_{i+\frac{1}{2}}^{+}, y_{j}^{(\mu)}\right)
\right),
$$
$$
\mathscr{V}_{2, j+\frac{1}{2}}^{\mu, \pm}:=
\mathscr{V}_{2}^{\pm}\left(
\widetilde{\boldsymbol{U}}_\sigma^{n}\left(x_{i}^{(\mu)}, y_{j+\frac{1}{2}}^{-}\right),
\widetilde{\boldsymbol{U}}_\sigma^{n}\left(x_{i}^{(\mu)}, y_{j+\frac{1}{2}}^{+}\right)
\right),
$$
$$
\boldsymbol{H}_{1, i+\frac{1}{2}}^{\mu}:=
\boldsymbol{H}_{1}\left(
\widetilde{\boldsymbol{U}}_\sigma^{n}\left(x_{i+\frac{1}{2}}^{-}, y_{j}^{(\mu)}\right),
\widetilde{\boldsymbol{U}}_\sigma^{n}\left(x_{i+\frac{1}{2}}^{+}, y_{j}^{(\mu)}\right)
\right),
$$
$$
\boldsymbol{H}_{2, j+\frac{1}{2}}^{\mu}:=
\boldsymbol{H}_{2}\left(
\widetilde{\boldsymbol{U}}_\sigma^{n}\left(x_{i}^{(\mu)}, y_{j+\frac{1}{2}}^{-}\right),
\widetilde{\boldsymbol{U}}_\sigma^{n}\left(x_{i}^{(\mu)}, y_{j+\frac{1}{2}}^{+}\right)
\right),
$$
$$
\begin{aligned}
	\hat{\alpha}_{1, i\pm\frac{1}{2}}^{\mu, \mp}=
	&\max\left\{
	\pm\left(u_{1}\right)_{i\pm\frac{1}{2},j}^{\mp,\mu},
	\pm\left(
	\frac{\sqrt{\rho_{i\pm\frac{1}{2},j}^{\mp, \mu}} \left(u_{1}\right)_{i\pm\frac{1}{2},j}^{\mp,\mu}+
		\sqrt{\rho_{i\mp\frac{1}{2},j}^{\pm, \mu}}
		\left(u_{1}\right)_{i\mp\frac{1}{2},j}^{\pm,\mu}
	}{\sqrt{\rho_{i\pm\frac{1}{2},j}^{\mp, \mu}}+\sqrt{\rho_{i\mp\frac{1}{2},j}^{\pm, \mu}}}\right)\right\}\\
	&+\mathscr{C}_{1}\left(\boldsymbol{U}_{i\pm\frac{1}{2},j}^{\mp, \mu}\right)+
	\frac{\left\|\boldsymbol{B}_{i\pm\frac{1}{2},j}^{\mp, \mu} - \boldsymbol{B}_{i\mp\frac{1}{2},j}^{\pm, \mu}\right\|}{\sqrt{\rho_{i\pm\frac{1}{2},j}^{\mp, \mu}}+\sqrt{\rho_{i\mp\frac{1}{2},j}^{\pm, \mu}}},\\
\end{aligned}
$$
$$
\begin{aligned}
	\hat{\alpha}_{2, j\pm\frac{1}{2}}^{\mu, \mp}=
	&\max\left\{
	\pm\left(u_{2}\right)_{i,j\pm\frac{1}{2}}^{\mu,\mp},
	\pm\left(
	\frac{ \sqrt{\rho_{i,j\pm\frac{1}{2}}^{\mu,\mp}} \left(u_{2}\right)_{i,j\pm\frac{1}{2}}^{\mu,\mp}+
		\sqrt{\rho_{i,j\mp\frac{1}{2}}^{\mu,\pm}}
		\left(u_{2}\right)_{i,j\mp\frac{1}{2}}^{\mu,\pm}
	}{\sqrt{\rho_{i,j\pm\frac{1}{2}}^{\mu,\mp}}+\sqrt{\rho_{i,j\mp\frac{1}{2}}^{ \mu,\pm}}}\right)\right\}\\
	&+\mathscr{C}_{2}\left(\boldsymbol{U}_{i,j\pm\frac{1}{2}}^{\mu,\mp}\right)+
	\frac{\left\|\boldsymbol{B}_{i,j\pm\frac{1}{2}}^{\mu,\mp} - \boldsymbol{B}_{i,j\mp\frac{1}{2}}^{\mu,\pm}\right\|}{\sqrt{\rho_{i,j\pm\frac{1}{2}}^{ \mu,\mp}}+\sqrt{\rho_{i,j\mp\frac{1}{2}}^{\mu,\pm}}}.\\
\end{aligned}
$$

The PP property of the updated cell averages in scheme (\ref{Equ:updatecellave}) is shown as follows. 

\begin{theorem}\label{Thm:2DPPpropertytheorem}
	Let $\left\{ \widetilde{\boldsymbol{U}}_{ij}^{n}\left( \boldsymbol{x}\right)\right\}$ denote the OEDG solution polynomial vector in cell $I_{ij}$ after the OE and PP limiting procedures, so that $\widetilde{\boldsymbol{B}}_{ij}^{n}\left( \boldsymbol{x}\right)$ is LDF and 
	$\widetilde{\boldsymbol{U}}_{ij}^{n}\left( \boldsymbol{x}\right)$ satisfies condition (\ref{Equ:2DPPlimitercond}). 
	If the wave speeds in the HLL flux satisfy the condition (\ref{Equ:2DPPwavespeedcon}), then the scheme (\ref{Equ:updatecellave}) preserves  $\bar{\boldsymbol{U}}_{ij}^{n+1}\in\mathcal{G}$ under the CFL-type condition
	\begin{equation}\label{Equ:2DCFLcondi}
		\Delta t_{n} \leq \min\left\{
		\frac{\Delta x \omega_{1}^{+}}{\alpha_{1,i+\frac{1}{2}}^{\mu}},
		\frac{\Delta x \omega_{1}^{-}}{\alpha_{1,i-\frac{1}{2}}^{\mu}},
		\frac{\Delta y \omega_{2}^{+}}{\alpha_{2,j+\frac{1}{2}}^{\mu}},
		\frac{\Delta y \omega_{2}^{-}}{\alpha_{2,j-\frac{1}{2}}^{\mu}}
		\right\}  \quad \forall  i,j,  1\leq \mu \leq q,
	\end{equation}
	where 
	$$
	\alpha_{1,i+\frac{1}{2}}^{\mu}=\hat{\alpha}_{1,i+\frac{1}{2}}^{\mu,-}-
	\mathscr{V}_{1,i+\frac{1}{2}}^{\mu,-}+
	\frac{\left|\mathscr{B}_{1, i+\frac{1}{2}}^{-}
		\left(y_{j}^{(\mu)}\right) \right|}{\sqrt{\rho_{i+\frac{1}{2},j}^{-,\mu}}},
	\quad
	\alpha_{1,i-\frac{1}{2}}^{\mu}=\hat{\alpha}_{1,i-\frac{1}{2}}^{\mu,+}+
	\mathscr{V}_{1,i-\frac{1}{2}}^{\mu,+}+
	\frac{\left|\mathscr{B}_{1, i-\frac{1}{2}}^{+}
		\left(y_{j}^{(\mu)}\right) \right|}{\sqrt{\rho_{i-\frac{1}{2},j}^{+,\mu}}},
	$$
	$$
	\alpha_{2,j+\frac{1}{2}}^{\mu}=\hat{\alpha}_{2,j+\frac{1}{2}}^{\mu,-}-
	\mathscr{V}_{2,j+\frac{1}{2}}^{\mu,-}+
	\frac{\left|\mathscr{B}_{2, j+\frac{1}{2}}^{-}
		\left(x_{i}^{(\mu)}\right) \right|}{\sqrt{\rho_{i,j+\frac{1}{2}}^{\mu,-}}},
	\quad
	\alpha_{2,j-\frac{1}{2}}^{\mu}=\hat{\alpha}_{2,j-\frac{1}{2}}^{\mu,+}+
	\mathscr{V}_{2,j-\frac{1}{2}}^{\mu,+}+
	\frac{\left|\mathscr{B}_{2, j-\frac{1}{2}}^{+}
		\left(x_{i}^{(\mu)}\right) \right|}{\sqrt{\rho_{i,j-\frac{1}{2}}^{\mu,+}}}.
	$$
\end{theorem}

\noindent \textbf{\emph{Proof }}
According to the identities (\ref{Equ:2DPPfluxequal1})--(\ref{Equ:2DPPfluxequal2}) in Lemma \ref{Thm:2DPPfluxtheorem}, the scheme (\ref{Equ:updatecellave}) can be rewritten as
\begin{equation*}
	\begin{aligned}
		\bar{\boldsymbol{U}}_{ij}^{n+1}=\bar{\boldsymbol{U}}_{ij}^{n}-&\frac{\Delta t_{n}}{\Delta x}
		\sum_{\mu=1}^{q}\omega_{\mu}^{G}\left[\left( \mathscr{V}_{1,i+\frac{1}{2}}^{\mu,-}\boldsymbol{H}_{1,i+\frac{1}{2}}^{\mu}+\boldsymbol{F}_{1}
		\left(\boldsymbol{U}_{i+\frac{1}{2},j}^{-,\mu}\right)-\mathscr{V}_{1,i+\frac{1}{2}}^{\mu,-}
		\boldsymbol{U}_{i+\frac{1}{2},j}^{-,\mu}
		\right)\right.\\
		& -\left( \mathscr{V}_{1,i-\frac{1}{2}}^{\mu,+}\boldsymbol{H}_{1,i-\frac{1}{2}}^{\mu}+\boldsymbol{F}_{1}
		\left(\boldsymbol{U}_{i-\frac{1}{2},j}^{+,\mu}\right)-\mathscr{V}_{1,i-\frac{1}{2}}^{\mu,+}
		\boldsymbol{U}_{i-\frac{1}{2},j}^{+,\mu}
		\right)\\
		&\left.+\mathscr{B}_{1, i+\frac{1}{2}}^{-}
		\left(y_{j}^{(\mu)}\right)\boldsymbol{S}\left(\boldsymbol{U}_{i+\frac{1}{2},j}^{-,\mu}\right)+
		\mathscr{B}_{1, i-\frac{1}{2}}^{+}
		\left(y_{j}^{(\mu)}\right)\boldsymbol{S}\left(\boldsymbol{U}_{i-\frac{1}{2},j}^{+,\mu}\right)\right]\\
		-&\frac{\Delta t_{n}}{\Delta y}
		\sum_{\mu=1}^{q}\omega_{\mu}^{G}\left[\left( \mathscr{V}_{2,j+\frac{1}{2}}^{\mu,-}\boldsymbol{H}_{2,j+\frac{1}{2}}^{\mu}+\boldsymbol{F}_{2}
		\left(\boldsymbol{U}_{i,j+\frac{1}{2}}^{\mu,-}\right)-\mathscr{V}_{2,j+\frac{1}{2}}^{\mu,-}
		\boldsymbol{U}_{i,j+\frac{1}{2}}^{\mu,-}
		\right)\right.\\
		& -\left( \mathscr{V}_{2,j-\frac{1}{2}}^{\mu,+}\boldsymbol{H}_{2,j-\frac{1}{2}}^{\mu}+\boldsymbol{F}_{2}
		\left(\boldsymbol{U}_{i,j-\frac{1}{2}}^{\mu,+}\right)-\mathscr{V}_{2,j-\frac{1}{2}}^{\mu,+}
		\boldsymbol{U}_{i,j-\frac{1}{2}}^{\mu,+}
		\right)\\
		&\left.+\mathscr{B}_{2, j+\frac{1}{2}}^{-}
		\left(x_{i}^{(\mu)}\right)\boldsymbol{S}\left(\boldsymbol{U}_{i,j+\frac{1}{2}}^{\mu,-}\right)+
		\mathscr{B}_{2, j-\frac{1}{2}}^{+}
		\left(x_{i}^{(\mu)}\right)\boldsymbol{S}\left(\boldsymbol{U}_{i,j-\frac{1}{2}}^{\mu,+}\right)\right].
	\end{aligned}
\end{equation*}
It can be further split into 
\begin{equation}\label{Equ:rewrite2Davescheme}
	\bar{\boldsymbol{U}}_{ij}^{n+1}=\bar{\boldsymbol{U}}_{ij}^n+
	\boldsymbol{\Xi}_1+\boldsymbol{\Xi}_2+\boldsymbol{\Xi}_3+\boldsymbol{\Xi}_4
\end{equation}
with
$$
\begin{aligned}
	\boldsymbol{\Xi}_1:=&\frac{\Delta t_{n}}{\Delta x}\sum_{\mu=1}^{q}\omega_{\mu}^{G}\left[
	\left(\mathscr{V}_{1,i+\frac{1}{2}}^{\mu,-}-\hat{\alpha}_{1,i+\frac{1}{2}}^{\mu,-}\right)
	\boldsymbol{U}_{i+\frac{1}{2},j}^{-,\mu}-\left(\mathscr{V}_{1,i-\frac{1}{2}}^{\mu,+}+
	\hat{\alpha}_{1,i-\frac{1}{2}}^{\mu,+}\right)
	\boldsymbol{U}_{i-\frac{1}{2},j}^{+,\mu}
	\right]\\
	&+\frac{\Delta t_{n}}{\Delta y}\sum_{\mu=1}^{q}\omega_{\mu}^{G}\left[
	\left(\mathscr{V}_{2,j+\frac{1}{2}}^{\mu,-}-\hat{\alpha}_{2,j+\frac{1}{2}}^{\mu,-}\right)
	\boldsymbol{U}_{i,j+\frac{1}{2}}^{\mu,-}-\left(\mathscr{V}_{2,j-\frac{1}{2}}^{\mu,+}+
	\hat{\alpha}_{2,j-\frac{1}{2}}^{\mu,+}\right)
	\boldsymbol{U}_{i,j-\frac{1}{2}}^{\mu,+}
	\right],\\
\end{aligned}
$$
$$
\begin{aligned}
	\boldsymbol{\Xi}_2:=&\frac{\Delta t_{n}}{\Delta x}\sum_{\mu=1}^{q}\omega_{\mu}^{G}\left[
	\left(\hat{\alpha}_{1,i+\frac{1}{2}}^{\mu,-}\boldsymbol{U}_{i+\frac{1}{2},j}^{-,\mu}-
	\boldsymbol{F}_{1}\left(\boldsymbol{U}_{i+\frac{1}{2},j}^{-,\mu}\right)\right)+
	\left(\hat{\alpha}_{1,i-\frac{1}{2}}^{\mu,+}\boldsymbol{U}_{i-\frac{1}{2},j}^{+,\mu}+
	\boldsymbol{F}_{1}\left(\boldsymbol{U}_{i-\frac{1}{2},j}^{+,\mu}\right)\right)\right]\\
	&+\frac{\Delta t_{n}}{\Delta y}\sum_{\mu=1}^{q}\omega_{\mu}^{G}\left[
	\left(\hat{\alpha}_{2,j+\frac{1}{2}}^{\mu,-}\boldsymbol{U}_{i,j+\frac{1}{2}}^{\mu,-}-
	\boldsymbol{F}_{2}\left(\boldsymbol{U}_{i,j+\frac{1}{2}}^{\mu,-}\right)\right)+
	\left(\hat{\alpha}_{2,j-\frac{1}{2}}^{\mu,+}\boldsymbol{U}_{i,j-\frac{1}{2}}^{\mu,+}+
	\boldsymbol{F}_{2}\left(\boldsymbol{U}_{i,j-\frac{1}{2}}^{\mu,+}\right)\right)\right],\\
\end{aligned}
$$
$$
\begin{aligned}
	\boldsymbol{\Xi}_3:=&-\frac{\Delta t_{n}}{\Delta x}\sum_{\mu=1}^{q}\omega_{\mu}^{G}\left(
	\mathscr{V}_{1,i+\frac{1}{2}}^{\mu,-}\boldsymbol{H}_{1,i+\frac{1}{2}}^{\mu}-
	\mathscr{V}_{1,i-\frac{1}{2}}^{\mu,+}\boldsymbol{H}_{1,i-\frac{1}{2}}^{\mu}\right)
	-\frac{\Delta t_{n}}{\Delta y}\sum_{\mu=1}^{q}\omega_{\mu}^{G}\left(
	\mathscr{V}_{2,j+\frac{1}{2}}^{\mu,-}\boldsymbol{H}_{2,j+\frac{1}{2}}^{\mu}-
	\mathscr{V}_{2,j-\frac{1}{2}}^{\mu,+}\boldsymbol{H}_{2,j-\frac{1}{2}}^{\mu}\right),\\
\end{aligned}
$$
$$
\begin{aligned}
	\boldsymbol{\Xi}_4:=&-\frac{\Delta t_{n}}{\Delta x}\sum_{\mu=1}^{q}\omega_{\mu}^{G}
	\left(\mathscr{B}_{1, i+\frac{1}{2}}^{-}
	\left(y_{j}^{(\mu)}\right)\boldsymbol{S}\left(\boldsymbol{U}_{i+\frac{1}{2},j}^{-,\mu}\right)+
	\mathscr{B}_{1, i-\frac{1}{2}}^{+}
	\left(y_{j}^{(\mu)}\right)\boldsymbol{S}\left(\boldsymbol{U}_{i-\frac{1}{2},j}^{+,\mu}\right)\right)\\
	&-\frac{\Delta t_{n}}{\Delta y}\sum_{\mu=1}^{q}\omega_{\mu}^{G}
	\left(\mathscr{B}_{2, j+\frac{1}{2}}^{-}
	\left(x_{i}^{(\mu)}\right)\boldsymbol{S}\left(\boldsymbol{U}_{i,j+\frac{1}{2}}^{\mu,-}\right)+
	\mathscr{B}_{2, j-\frac{1}{2}}^{+}
	\left(x_{i}^{(\mu)}\right)\boldsymbol{S}\left(\boldsymbol{U}_{i,j-\frac{1}{2}}^{\mu,+}\right)\right).
\end{aligned}
$$
We define, for any $1\leq\mu\leq q$,
$$
\begin{aligned}
	&\bar{\boldsymbol{U}}_{1,ij}^{(\mu)}:=\frac{1}{\hat{\alpha}_{1,i+\frac{1}{2}}^{\mu,-}
		+\hat{\alpha}_{1,i-\frac{1}{2}}^{\mu,+}}\left[ \left(\hat{\alpha}_{1,i+\frac{1}{2}}^{\mu,-}\boldsymbol{U}_{i+\frac{1}{2},j}^{-,\mu}-
	\boldsymbol{F}_{1}\left(\boldsymbol{U}_{i+\frac{1}{2},j}^{-,\mu}\right)\right)\right.
	+\left.\left(\hat{\alpha}_{1,i-\frac{1}{2}}^{\mu,+}\boldsymbol{U}_{i-\frac{1}{2},j}^{+,\mu}+
	\boldsymbol{F}_{1}\left(\boldsymbol{U}_{i-\frac{1}{2},j}^{+,\mu}\right)\right)\right],\\
\end{aligned}
$$
$$
\begin{aligned}
	&\bar{\boldsymbol{U}}_{2,ij}^{(\mu)}:=\frac{1}{\hat{\alpha}_{2,j+\frac{1}{2}}^{
			\mu,-}+\hat{\alpha}_{2,j-\frac{1}{2}}^{\mu,+}}\left[\left(\hat{\alpha}_{2,j+\frac{1}{2}}^{\mu,-}\boldsymbol{U}_{i,j+\frac{1}{2}}^{\mu,-}-
	\boldsymbol{F}_{2}\left(\boldsymbol{U}_{i,j+\frac{1}{2}}^{\mu,-}\right)\right)\right.
	+\left.\left(\hat{\alpha}_{2,j-\frac{1}{2}}^{\mu,+}\boldsymbol{U}_{i,j-\frac{1}{2}}^{\mu,+}+
	\boldsymbol{F}_{2}\left(\boldsymbol{U}_{i,j-\frac{1}{2}}^{\mu,+}\right)\right)\right],
\end{aligned}
$$
then $\boldsymbol{\Xi}_2$ is equal to
$$
\begin{aligned}
	\boldsymbol{\Xi}_2=\frac{\Delta t_{n}}{\Delta x}\sum_{\mu=1}^{q}\omega_{\mu}^{G}
	\left(\hat{\alpha}_{1,i+\frac{1}{2}}^{\mu,-}
	+\hat{\alpha}_{1,i-\frac{1}{2}}^{\mu,+}\right)\bar{\boldsymbol{U}}_{1,ij}^{(\mu)}
	+\frac{\Delta t_{n}}{\Delta y}\sum_{\mu=1}^{q}\omega_{\mu}^{G}
	\left(\hat{\alpha}_{2,j+\frac{1}{2}}^{\mu,-}+\hat{\alpha}_{2,j-\frac{1}{2}}^{\mu,+}
	\right)\bar{\boldsymbol{U}}_{2,ij}^{(\mu)}.
\end{aligned}
$$

Now, we first show $\bar{\rho}_{ij}^{n+1}>0$. 
Thanks to Lemma \ref{Thm:2DPPstatetheorem}, we have, for all $1\leq\mu\leq q$, $\bar{\boldsymbol{U}}_{1,ij}^{(\mu)}\in\mathcal{G}_{\rho}$ and 
$\bar{\boldsymbol{U}}_{2,ij}^{(\mu)}\in\mathcal{G}_{\rho}$. Since 
$\hat{\alpha}_{1,i+\frac{1}{2}}^{\mu,-}
+\hat{\alpha}_{1,i-\frac{1}{2}}^{\mu,+}>0$ 
and $\hat{\alpha}_{2,j+\frac{1}{2}}^{\mu,-}+\hat{\alpha}_{2,j-\frac{1}{2}}^{\mu,+}
>0$,  we obtain $\boldsymbol{\Xi}_2\in\mathcal{G}_{\rho}$. 
According to Lemma \ref{Thm:2DPPfluxtheorem}, we have $\boldsymbol{H}_{1,i+\frac{1}{2}}^{\mu}$, $\boldsymbol{H}_{1,i-\frac{1}{2}}^{\mu}$, $\boldsymbol{H}_{2,j+\frac{1}{2}}^{\mu}$, and $\boldsymbol{H}_{2,j-\frac{1}{2}}^{\mu}$ all belong to $\mathcal{G}_{\rho}$,
which means $\boldsymbol{\Xi}_3\in\mathcal{G}_{\rho}$.
Since the first component of $\boldsymbol{\Xi}_4$ is zero, we deduce from (\ref{Equ:rewrite2Davescheme}) that
\begin{align*}
	\bar{\rho}_{ij}^{n+1}&>\bar{\rho}_{ij}^{n}+\frac{\Delta t_{n}}{\Delta x}\sum_{\mu=1}^{q}\omega_{\mu}^{G}\left[
	\left(\mathscr{V}_{1,i+\frac{1}{2}}^{\mu,-}-\hat{\alpha}_{1,i+\frac{1}{2}}^{\mu,-}\right)
	\rho_{i+\frac{1}{2},j}^{-,\mu}-\left(\mathscr{V}_{1,i-\frac{1}{2}}^{\mu,+}+
	\hat{\alpha}_{1,i-\frac{1}{2}}^{\mu,+}\right)
	\rho_{i-\frac{1}{2},j}^{+,\mu}
	\right]\\
	&\quad +\frac{\Delta t_{n}}{\Delta y}\sum_{\mu=1}^{q}\omega_{\mu}^{G}\left[
	\left(\mathscr{V}_{2,j+\frac{1}{2}}^{\mu,-}-\hat{\alpha}_{2,j+\frac{1}{2}}^{\mu,-}\right)
	\rho_{i,j+\frac{1}{2}}^{\mu,-}-\left(\mathscr{V}_{2,j-\frac{1}{2}}^{\mu,+}+
	\hat{\alpha}_{2,j-\frac{1}{2}}^{\mu,+}\right)
	\rho_{i,j-\frac{1}{2}}^{\mu,+}
	\right]\\
	&=\sum_{\mu=1}^q \omega_\mu^{G}\left[\omega_1^{+} \rho_{i+\frac{1}{2},j}^{-,\mu}+\omega_1^{-}
	\rho_{i-\frac{1}{2},j}^{+,\mu}+\omega_2^{+} \rho_{i,j+\frac{1}{2}}^{\mu,-}+\omega_2^{-}
	\rho_{i,j-\frac{1}{2}}^{\mu,+}\right]
	+\sum_{s=1}^S \omega_s \widetilde{\rho}_{\sigma}^{n}\left(x_{ij}^{(s)}, y_{ij}^{(s)}\right)\\
	&\quad+\frac{\Delta t_{n}}{\Delta x}\sum_{\mu=1}^{q}\omega_{\mu}^{G}\left[
	\left(\mathscr{V}_{1,i+\frac{1}{2}}^{\mu,-}-\hat{\alpha}_{1,i+\frac{1}{2}}^{\mu,-}\right)
	\rho_{i+\frac{1}{2},j}^{-,\mu}-\left(\mathscr{V}_{1,i-\frac{1}{2}}^{\mu,+}+
	\hat{\alpha}_{1,i-\frac{1}{2}}^{\mu,+}\right)
	\rho_{i-\frac{1}{2},j}^{+,\mu}
	\right]\\
	&\quad+\frac{\Delta t_{n}}{\Delta y}\sum_{\mu=1}^{q}\omega_{\mu}^{G}\left[
	\left(\mathscr{V}_{2,j+\frac{1}{2}}^{\mu,-}-\hat{\alpha}_{2,j+\frac{1}{2}}^{\mu,-}\right)
	\rho_{i,j+\frac{1}{2}}^{\mu,-}-\left(\mathscr{V}_{2,j-\frac{1}{2}}^{\mu,+}+
	\hat{\alpha}_{2,j-\frac{1}{2}}^{\mu,+}\right)
	\rho_{i,j-\frac{1}{2}}^{\mu,+}
	\right],\\
	&\geq\sum_{\mu=1}^{q}\omega_{\mu}^{G}\rho_{i+\frac{1}{2},j}^{-,\mu}\left(\omega_{1}^{+}+
	\frac{\Delta t_{n}}{\Delta x}\left(\mathscr{V}_{1,i+\frac{1}{2}}^{\mu,-}-\hat{\alpha}_{1,i+\frac{1}{2}}^{\mu,-}\right) \right)\\
	&\quad+\sum_{\mu=1}^{q}\omega_{\mu}^{G}\rho_{i-\frac{1}{2},j}^{+,\mu}\left(\omega_{1}^{-}-\frac{\Delta t_{n}}{\Delta x}\left(\mathscr{V}_{1,i-\frac{1}{2}}^{\mu,+}+
	\hat{\alpha}_{1,i-\frac{1}{2}}^{\mu,+}\right)\right)\\
	&\quad+\sum_{\mu=1}^{q}\omega_{\mu}^{G}\rho_{i,j+\frac{1}{2}}^{\mu,-}\left(\omega_{2}^{+}+\frac{\Delta t_{n}}{\Delta y}\left(\mathscr{V}_{2,j+\frac{1}{2}}^{\mu,-}-
	\hat{\alpha}_{2,j+\frac{1}{2}}^{\mu,-}\right)\right)\\
	&\quad+\sum_{\mu=1}^{q}\omega_{\mu}^{G}\rho_{i,j-\frac{1}{2}}^{\mu,+}\left(\omega_{2}^{-}-\frac{\Delta t_{n}}{\Delta y}\left(\mathscr{V}_{2,j-\frac{1}{2}}^{\mu,+}+
	\hat{\alpha}_{2,j-\frac{1}{2}}^{\mu,+}\right)\right)\\
	&\geq0,
\end{align*}
where we have used the convex decomposition (\ref{Equ:2DgeneralCD}) in the above equality, the condition (\ref{Equ:2DPPlimitercond}) in the second inequality, and the CFL condition (\ref{Equ:2DCFLcondi}) in the last inequality.

We then prove that $\bar{\boldsymbol{U}}_{ij}^{n+1}\cdot\boldsymbol{n}^{*}+\frac{{\left\|\boldsymbol{B}^{*} \right\|}^{2}}{2}>0$ for any $\boldsymbol{v}^{*}, \boldsymbol{B}^{*} \in \mathbb{R}^{3}$.
It follows from (\ref{Equ:rewrite2Davescheme}) that
\begin{equation}\label{Equ:vBequality}
	\bar{\boldsymbol{U}}_{ij}^{n+1} \cdot \boldsymbol{n}^*+\frac{\left\|\boldsymbol{B}^*\right\|^2}{2}=\Pi_0+\Pi_1+\Pi_2+\Pi_3+\Pi_4,
\end{equation}
where
$$
\Pi_0:=\bar{\boldsymbol{U}}_{ij}^{n} \cdot \boldsymbol{n}^*+\frac{\left\|\boldsymbol{B}^*\right\|^2}{2},
$$
$$
\begin{aligned}
	\Pi_1:=&\frac{\Delta t_{n}}{\Delta x}\sum_{\mu=1}^{q}\omega_{\mu}^{G}\left[
	\left(\mathscr{V}_{1,i+\frac{1}{2}}^{\mu,-}-\hat{\alpha}_{1,i+\frac{1}{2}}^{\mu,-}\right)
	\left(\boldsymbol{U}_{i+\frac{1}{2},j}^{-,\mu}\cdot \boldsymbol{n}^*+\frac{\left\|\boldsymbol{B}^*\right\|^2}{2}\right)\right.\\
	&\left.-\left(\mathscr{V}_{1,i-\frac{1}{2}}^{\mu,+}+
	\hat{\alpha}_{1,i-\frac{1}{2}}^{\mu,+}\right)
	\left(\boldsymbol{U}_{i-\frac{1}{2},j}^{+,\mu}\cdot \boldsymbol{n}^*+\frac{\left\|\boldsymbol{B}^*\right\|^2}{2}\right)
	\right]\\
	&+\frac{\Delta t_{n}}{\Delta y}\sum_{\mu=1}^{q}\omega_{\mu}^{G}\left[
	\left(\mathscr{V}_{2,j+\frac{1}{2}}^{\mu,-}-\hat{\alpha}_{2,j+\frac{1}{2}}^{\mu,-}\right)
	\left(\boldsymbol{U}_{i,j+\frac{1}{2}}^{\mu,-}\cdot \boldsymbol{n}^*+\frac{\left\|\boldsymbol{B}^*\right\|^2}{2}\right)\right.\\
	&\left.-\left(\mathscr{V}_{2,j-\frac{1}{2}}^{\mu,+}+
	\hat{\alpha}_{2,j-\frac{1}{2}}^{\mu,+}\right)
	\left(\boldsymbol{U}_{i,j-\frac{1}{2}}^{\mu,+}\cdot \boldsymbol{n}^*+\frac{\left\|\boldsymbol{B}^*\right\|^2}{2}\right)
	\right],\\
\end{aligned}
$$
$$
\begin{aligned}
	\Pi_2:=&\frac{\Delta t_{n}}{\Delta x }\sum_{\mu=1}^{q}\omega_{\mu}^{G}
	\left(\hat{\alpha}_{1,i+\frac{1}{2}}^{\mu,-}
	+\hat{\alpha}_{1,i-\frac{1}{2}}^{\mu,+}\right)
	\left(\bar{\boldsymbol{U}}_{1,ij}^{(\mu)}\cdot \boldsymbol{n}^*+\frac{\left\|\boldsymbol{B}^*\right\|^2}{2}\right)\\
	&+\frac{\Delta t_{n}}{\Delta y }\sum_{\mu=1}^{q}\omega_{\mu}^{G}\left(\hat{\alpha}_{2,j+\frac{1}{2}}^{\mu,-}+\hat{\alpha}_{2,j-\frac{1}{2}}^{\mu,+}
	\right)\left(\bar{\boldsymbol{U}}_{2,ij}^{(\mu)}\cdot \boldsymbol{n}^*+\frac{\left\|\boldsymbol{B}^*\right\|^2}{2}\right),
\end{aligned}
$$
$$
\begin{aligned}
	\Pi_3:=&-\frac{\Delta t_{n}}{\Delta x}\sum_{\mu=1}^{q}\omega_{\mu}^{G}\left[
	\mathscr{V}_{1,i+\frac{1}{2}}^{\mu,-}\left(\boldsymbol{H}_{1,i+\frac{1}{2}}^{\mu}\cdot \boldsymbol{n}^*+\frac{\left\|\boldsymbol{B}^*\right\|^2}{2}\right)
	-\mathscr{V}_{1,i-\frac{1}{2}}^{\mu,+}\left(\boldsymbol{H}_{1,i-\frac{1}{2}}^{\mu}\cdot \boldsymbol{n}^*+\frac{\left\|\boldsymbol{B}^*\right\|^2}{2}\right)\right]\\
	&-\frac{\Delta t_{n}}{\Delta y}\sum_{\mu=1}^{q}\omega_{\mu}^{G}\left[
	\mathscr{V}_{2,j+\frac{1}{2}}^{\mu,-}\left(\boldsymbol{H}_{2,j+\frac{1}{2}}^{\mu}\cdot \boldsymbol{n}^*+\frac{\left\|\boldsymbol{B}^*\right\|^2}{2}\right)
	-\mathscr{V}_{2,j-\frac{1}{2}}^{\mu,+}\left(\boldsymbol{H}_{2,j-\frac{1}{2}}^{\mu}\cdot \boldsymbol{n}^*+\frac{\left\|\boldsymbol{B}^*\right\|^2}{2}\right)\right],\\
\end{aligned}
$$
$$
\begin{aligned}
	\Pi_4:=&-\frac{\Delta t_{n}}{\Delta x}\sum_{\mu=1}^{q}\omega_{\mu}^{G}
	\left(\mathscr{B}_{1, i+\frac{1}{2}}^{-}
	\left(y_{j}^{(\mu)}\right)\boldsymbol{S}\left(\boldsymbol{U}_{i+\frac{1}{2},j}^{-,\mu}\right)\cdot \boldsymbol{n}^*+
	\mathscr{B}_{1, i-\frac{1}{2}}^{+}
	\left(y_{j}^{(\mu)}\right)\boldsymbol{S}\left(\boldsymbol{U}_{i-\frac{1}{2},j}^{+,\mu}\right)\cdot \boldsymbol{n}^*
	\right)\\
	&-\frac{\Delta t_{n}}{\Delta y}\sum_{\mu=1}^{q}\omega_{\mu}^{G}
	\left(\mathscr{B}_{2, j+\frac{1}{2}}^{-}
	\left(x_{i}^{(\mu)}\right)\boldsymbol{S}\left(\boldsymbol{U}_{i,j+\frac{1}{2}}^{\mu,-}\right)\cdot \boldsymbol{n}^*+
	\mathscr{B}_{2, j-\frac{1}{2}}^{+}
	\left(x_{i}^{(\mu)}\right)\boldsymbol{S}\left(\boldsymbol{U}_{i,j-\frac{1}{2}}^{\mu,+}\right)
	\cdot \boldsymbol{n}^*\right).
\end{aligned}
$$
We now estimate the lower bounds of $\Pi_0$, $\Pi_2$, $\Pi_3$, and $\Pi_4$, respectively.
Based on the the convex decomposition (\ref{Equ:2DgeneralCD}), we have
\begin{equation}\label{Equ:pi0bound}
	\begin{aligned}
		\Pi_0=&\sum_{\mu=1}^q \omega_\mu^{G}\left[\omega_1^{+} \left(\boldsymbol{U}_{i+\frac{1}{2},j}^{-,\mu}\cdot \boldsymbol{n}^*+\frac{\left\|\boldsymbol{B}^*\right\|^2}{2}\right)+\omega_1^{-}
		\left(\boldsymbol{U}_{i-\frac{1}{2},j}^{+,\mu}\cdot \boldsymbol{n}^*+\frac{\left\|\boldsymbol{B}^*\right\|^2}{2}\right)+\omega_2^{+} \left(\boldsymbol{U}_{i,j+\frac{1}{2}}^{\mu,-}\cdot \boldsymbol{n}^*+\frac{\left\|\boldsymbol{B}^*\right\|^2}{2}\right)\right.\\
		&\left.+\omega_2^{-}
		\left(\boldsymbol{U}_{i,j-\frac{1}{2}}^{\mu,+}\cdot \boldsymbol{n}^*+\frac{\left\|\boldsymbol{B}^*\right\|^2}{2}\right)\right]+\sum_{s=1}^S \omega_s \left(\widetilde{\boldsymbol{U}}_{\sigma}^{n}\left(x_{ij}^{(s)}, y_{ij}^{(s)}\right)\cdot \boldsymbol{n}^*+\frac{\left\|\boldsymbol{B}^*\right\|^2}{2}\right),\\
		&\geq\sum_{\mu=1}^q \omega_\mu^{G}\left[\omega_1^{+} \left(\boldsymbol{U}_{i+\frac{1}{2},j}^{-,\mu}\cdot \boldsymbol{n}^*+\frac{\left\|\boldsymbol{B}^*\right\|^2}{2}\right)+\omega_1^{-}
		\left(\boldsymbol{U}_{i-\frac{1}{2},j}^{+,\mu}\cdot \boldsymbol{n}^*+\frac{\left\|\boldsymbol{B}^*\right\|^2}{2}\right)\right.\\
		&\quad\left.+\omega_2^{+} \left(\boldsymbol{U}_{i,j+\frac{1}{2}}^{\mu,-}\cdot \boldsymbol{n}^*+\frac{\left\|\boldsymbol{B}^*\right\|^2}{2}\right)
		+\omega_2^{-}
		\left(\boldsymbol{U}_{i,j-\frac{1}{2}}^{\mu,+}\cdot \boldsymbol{n}^*+\frac{\left\|\boldsymbol{B}^*\right\|^2}{2}\right)\right],
	\end{aligned}
\end{equation}
where the inequality follows from the condition (\ref{Equ:2DPPlimitercond}) and  Lemma \ref{Lem:admistatequiv}.
According to Lemma \ref{Thm:2DPPstatetheorem}, for any  $\boldsymbol{v}^{*},\boldsymbol{B}^{*}\in \mathbb{R}^{3}$, we have 
$$
\bar{\boldsymbol{U}}_{1,ij}^{(\mu)}\cdot\boldsymbol{n}^{*}+
\frac{\left|\boldsymbol{B}^{*}\right|^{2}}{2}\geq
\frac{-\boldsymbol{v}^{*}\cdot\boldsymbol{B}^{*}}
{\hat{\alpha}_{1,i+\frac{1}{2}}^{\mu,-}
	+\hat{\alpha}_{1,i-\frac{1}{2}}^{\mu,+}}
\left(\left(B_{1}\right)_{i+\frac{1}{2},j}^{-,\mu}- \left(B_{1}\right)_{i-\frac{1}{2},j}^{+,\mu}\right),
$$
$$
\bar{\boldsymbol{U}}_{2,ij}^{(\mu)}\cdot\boldsymbol{n}^{*}+
\frac{\left|\boldsymbol{B}^{*}\right|^{2}}{2}\geq
\frac{-\boldsymbol{v}^{*}\cdot\boldsymbol{B}^{*}} {\hat{\alpha}_{2,j+\frac{1}{2}}^{\mu,-}+\hat{\alpha}_{2,j-\frac{1}{2}}^{\mu,+}
}\left(\left(B_{2}\right)_{i,j+\frac{1}{2}}^{\mu,-}- \left(B_{2}\right)_{i,j-\frac{1}{2}}^{\mu,+}\right).
$$
Therefore, the lower bound of $\Pi_2$ can be given by
\begin{equation}\label{Equ:pi2bound}
	\begin{aligned}
		\Pi_2& \geq -\frac{\Delta t_{n}}{\Delta x }\left(\boldsymbol{v}^{*}\cdot\boldsymbol{B}^{*} \right)\sum_{\mu=1}^{q}\omega_{\mu}^{G}
		\left( \left(B_{1}\right)_{i+\frac{1}{2},j}^{-,\mu}- \left(B_{1}\right)_{i-\frac{1}{2},j}^{+,\mu}\right)\\
		&\quad -\frac{\Delta t_{n}}{\Delta y }\left(\boldsymbol{v}^{*}\cdot\boldsymbol{B}^{*} \right)\sum_{\mu=1}^{q}\omega_{\mu}^{G}
		\left(\left(B_{2}\right)_{i,j+\frac{1}{2}}^{\mu,-}- \left(B_{2}\right)_{i,j-\frac{1}{2}}^{\mu,+}\right)\\
		&=-\Delta t_{n}\left(\boldsymbol{v}^{*}\cdot\boldsymbol{B}^{*} \right)\text{div}_{\text{local}}\widetilde{\boldsymbol{B}}_{\sigma}^{n}
		=0,
	\end{aligned}
\end{equation}
where the last equality uses the LDF property of the magnetic field within $I_{ij}$: 
\begin{equation*} 
	\begin{aligned}
		\text{div}_{\text{local}}\widetilde{\boldsymbol{B}}_{\sigma}^{n}& =\frac{1}{\Delta x}\sum_{\mu=1}^{q}\omega_{\mu}^{G}
		\left[{\left(B_{1}\right)}_{\sigma}^{n}\left(x_{i+\frac{1}{2}}^{-},
		y_{j}^{(\mu)}\right)-
		{\left(B_{1}\right)}_{\sigma}^{n}\left(x_{i-\frac{1}{2}}^{+}, y_{j}^{(\mu)}\right)\right]\\
		&\quad +\frac{1}{\Delta y}\sum_{\mu=1}^{q}\omega_{\mu}^{G}
		\left[{\left(B_{2}\right)}_{\sigma}^{n}\left(x_{i}^{(\mu)}, y_{j+\frac{1}{2}}^{-}\right)-
		{\left(B_{2}\right)}_{\sigma}^{n}\left(x_{i}^{(\mu)},y_{j-\frac{1}{2}}^{+}
		\right)\right]
		\\
		& =  \frac{1}{\Delta x \Delta y}\oint_{\partial I_{ij}} \widetilde{\boldsymbol{B}}_{ij}^{n}\left( \boldsymbol{x}\right)  \cdot
		\boldsymbol{n}\mathrm{~d}s = \frac{1}{\Delta x \Delta y}\int_{I_{ij}}
		\nabla \cdot \widetilde{\boldsymbol{B}}_{ij}^{n}\left( \boldsymbol{x}\right) \mathrm{~d} \boldsymbol{x} = 0. 
	\end{aligned}
\end{equation*}
By using the inequality (\ref{Equ:2DPPHinequal}) in Lemma \ref{Thm:2DPPfluxtheorem}, we give a lower bound of $\Pi_3$ as
\begin{equation}\label{Equ:pi3bound}
	\begin{aligned}
		\Pi_3\geq & -\frac{\Delta t_{n}}{\Delta x}\sum_{\mu=1}^{q}\omega_{\mu}^{G}\left[
		\mathscr{V}_{1,i+\frac{1}{2}}^{\mu,-}\left(-\frac{\boldsymbol{v}^*\cdot\boldsymbol{B}^*}
		{\mathscr{V}_{1,i+\frac{1}{2}}^{\mu,+}-\mathscr{V}_{1,i+\frac{1}{2}}^{\mu,-}}
		\left(\left(\boldsymbol{B}_{1}\right)_{i+\frac{1}{2},j}^{+,\mu}-
		\left(\boldsymbol{B}_{1}\right)_{i+\frac{1}{2},j}^{-,\mu}\right)\right)\right.\\
		&\left.-\mathscr{V}_{1,i-\frac{1}{2}}^{\mu,+}\left(-\frac{\boldsymbol{v}^*\cdot\boldsymbol{B}^*}
		{\mathscr{V}_{1,i-\frac{1}{2}}^{\mu,+}-\mathscr{V}_{1,i-\frac{1}{2}}^{\mu,-}}
		\left(\left(\boldsymbol{B}_{1}\right)_{i-\frac{1}{2},j}^{+,\mu}-
		\left(\boldsymbol{B}_{1}\right)_{i-\frac{1}{2},j}^{-,\mu}\right)\right)\right]\\
		&-\frac{\Delta t_{n}}{\Delta y}\sum_{\mu=1}^{q}\omega_{\mu}^{G}\left[
		\mathscr{V}_{2,j+\frac{1}{2}}^{\mu,-}\left(-\frac{\boldsymbol{v}^*\cdot\boldsymbol{B}^*}
		{\mathscr{V}_{2,j+\frac{1}{2}}^{\mu,+}-\mathscr{V}_{2,j+\frac{1}{2}}^{\mu,-}}
		\left(\left(\boldsymbol{B}_{2}\right)_{i,j+\frac{1}{2}}^{\mu,+}-
		\left(\boldsymbol{B}_{2}\right)_{i,j+\frac{1}{2}}^{\mu,-}\right)\right)\right.\\
		&-\left.\mathscr{V}_{2,j-\frac{1}{2}}^{\mu,+}\left(-\frac{\boldsymbol{v}^*\cdot\boldsymbol{B}^*}
		{\mathscr{V}_{2,j-\frac{1}{2}}^{\mu,+}-\mathscr{V}_{2,j-\frac{1}{2}}^{\mu,-}}
		\left(\left(\boldsymbol{B}_{2}\right)_{i,j-\frac{1}{2}}^{\mu,+}-
		\left(\boldsymbol{B}_{2}\right)_{i,j-\frac{1}{2}}^{\mu,-}\right)\right)\right].
	\end{aligned}
\end{equation}
Following the inequality (\ref{Equ:sourceinequiv1}) in Lemma \ref{Lem:sourcetermequiv}, we estimate the lower bound of $\Pi_4$ as
\begin{align}
		\Pi_4&\geq\frac{\Delta t_{n}}{\Delta x}\sum_{\mu=1}^{q}\omega_{\mu}^{G}\left(
		\mathscr{B}_{1, i+\frac{1}{2}}^{-}
		\left(y_{j}^{(\mu)}\right)\left(\boldsymbol{v}^*\cdot\boldsymbol{B}^*\right)-
		\frac{\left|\mathscr{B}_{1, i+\frac{1}{2}}^{-}
			\left(y_{j}^{(\mu)}\right) \right|}{\sqrt{\rho_{i+\frac{1}{2},j}^{-,\mu}}}\left( \boldsymbol{U}_{i+\frac{1}{2},j}^{-,\mu}\cdot\boldsymbol{n}^*+
		\frac{\left\|\boldsymbol{B}^*\right\|^2}{2}\right)\right)  \nonumber   \\
		&+\frac{\Delta t_{n}}{\Delta x}\sum_{\mu=1}^{q}\omega_{\mu}^{G}\left(
		\mathscr{B}_{1, i-\frac{1}{2}}^{+}
		\left(y_{j}^{(\mu)}\right)\left(\boldsymbol{v}^*\cdot\boldsymbol{B}^*\right)-
		\frac{\left|\mathscr{B}_{1, i-\frac{1}{2}}^{+}
			\left(y_{j}^{(\mu)}\right) \right|}{\sqrt{\rho_{i-\frac{1}{2},j}^{+,\mu}}}\left( \boldsymbol{U}_{i-\frac{1}{2},j}^{+,\mu}\cdot\boldsymbol{n}^*+
		\frac{\left\|\boldsymbol{B}^*\right\|^2}{2}\right)\right)  \nonumber  \\
		&+\frac{\Delta t_{n}}{\Delta y}\sum_{\mu=1}^{q}\omega_{\mu}^{G}\left(
		\mathscr{B}_{2, j+\frac{1}{2}}^{-}
		\left(x_{i}^{(\mu)}\right)\left(\boldsymbol{v}^*\cdot\boldsymbol{B}^*\right)-
		\frac{\left|\mathscr{B}_{2, j+\frac{1}{2}}^{-}
			\left(x_{i}^{(\mu)}\right) \right|}{\sqrt{\rho_{i,j+\frac{1}{2}}^{\mu,-}}}\left( \boldsymbol{U}_{i,j+\frac{1}{2}}^{\mu,-}\cdot\boldsymbol{n}^*+
		\frac{\left\|\boldsymbol{B}^*\right\|^2}{2}\right)\right) \nonumber  \\
		&+\frac{\Delta t_{n}}{\Delta y}\sum_{\mu=1}^{q}\omega_{\mu}^{G}\left(
		\mathscr{B}_{2, j-\frac{1}{2}}^{+}
		\left(x_{i}^{(\mu)}\right)\left(\boldsymbol{v}^*\cdot\boldsymbol{B}^*\right)-
		\frac{\left|\mathscr{B}_{2, j-\frac{1}{2}}^{+}
			\left(x_{i}^{(\mu)}\right) \right|}{\sqrt{\rho_{i,j-\frac{1}{2}}^{\mu,+}}}\left( \boldsymbol{U}_{i,j-\frac{1}{2}}^{\mu,+}\cdot\boldsymbol{n}^*+
		\frac{\left\|\boldsymbol{B}^*\right\|^2}{2}\right)\right). 
	\label{Equ:pi4bound}
\end{align}
Combining the lower bounds in (\ref{Equ:pi0bound}), (\ref{Equ:pi2bound}), (\ref{Equ:pi3bound}), (\ref{Equ:pi4bound}) with (\ref{Equ:vBequality}), we obtain
\begin{align*}
	&\bar{\boldsymbol{U}}_{ij}^{n+1}\cdot\boldsymbol{n}^*+
	\frac{\left\|\boldsymbol{B}^*\right\|^2}{2}\\
	&\geq\sum_{\mu=1}^{q}\omega_{\mu}^{G}\left(
	\omega_{1}^{+}+\frac{\Delta t_{n}}{\Delta x}\left(
	\mathscr{V}_{1,i+\frac{1}{2}}^{\mu,-}-\hat{\alpha}_{1,i+\frac{1}{2}}^{\mu,-}-
	\frac{\left|\mathscr{B}_{1, i+\frac{1}{2}}^{-}
		\left(y_{j}^{(\mu)}\right) \right|}{\sqrt{\rho_{i+\frac{1}{2},j}^{-,\mu}}}
	\right)\right)\left( \boldsymbol{U}_{i+\frac{1}{2},j}^{-,\mu}\cdot\boldsymbol{n}^*+
	\frac{\left\|\boldsymbol{B}^*\right\|^2}{2}\right)\\
	&+\sum_{\mu=1}^{q}\omega_{\mu}^{G}\left(
	\omega_{1}^{-}-\frac{\Delta t_{n}}{\Delta x}\left(
	\mathscr{V}_{1,i-\frac{1}{2}}^{\mu,+}+\hat{\alpha}_{1,i-\frac{1}{2}}^{\mu,+}+
	\frac{\left|\mathscr{B}_{1, i-\frac{1}{2}}^{+}
		\left(y_{j}^{(\mu)}\right) \right|}{\sqrt{\rho_{i-\frac{1}{2},j}^{+,\mu}}}
	\right)\right)\left( \boldsymbol{U}_{i-\frac{1}{2},j}^{+,\mu}\cdot\boldsymbol{n}^*+
	\frac{\left\|\boldsymbol{B}^*\right\|^2}{2}\right)\\
	&+\sum_{\mu=1}^{q}\omega_{\mu}^{G}\left(
	\omega_{2}^{+}+\frac{\Delta t_{n}}{\Delta y}\left(
	\mathscr{V}_{2,j+\frac{1}{2}}^{\mu,-}-\hat{\alpha}_{2,j+\frac{1}{2}}^{\mu,-}-
	\frac{\left|\mathscr{B}_{2, j+\frac{1}{2}}^{-}
		\left(x_{i}^{(\mu)}\right) \right|}{\sqrt{\rho_{i,j+\frac{1}{2}}^{\mu,-}}}
	\right)\right)\left( \boldsymbol{U}_{i,j+\frac{1}{2}}^{\mu,-}\cdot\boldsymbol{n}^*+
	\frac{\left\|\boldsymbol{B}^*\right\|^2}{2}\right)\\
	&+\sum_{\mu=1}^{q}\omega_{\mu}^{G}\left(
	\omega_{2}^{-}-\frac{\Delta t_{n}}{\Delta y}\left(
	\mathscr{V}_{2,j-\frac{1}{2}}^{\mu,+}+\hat{\alpha}_{2,j-\frac{1}{2}}^{\mu,+}+
	\frac{\left|\mathscr{B}_{2, j-\frac{1}{2}}^{+}
		\left(x_{i}^{(\mu)}\right) \right|}{\sqrt{\rho_{i,j-\frac{1}{2}}^{\mu,+}}}
	\right)\right)\left( \boldsymbol{U}_{i,j-\frac{1}{2}}^{\mu,+}\cdot\boldsymbol{n}^*+
	\frac{\left\|\boldsymbol{B}^*\right\|^2}{2}\right)\\
	&>0
\end{align*}
where the last inequality follows from the condition (\ref{Equ:2DPPlimitercond}) and CFL condition (\ref{Equ:2DCFLcondi}).
In summary, we have
$$
\bar{\rho}_{ij}^{n+1}>0,\quad
\bar{\boldsymbol{U}}_{ij}^{n+1}\cdot\boldsymbol{n}^*+
\frac{\left\|\boldsymbol{B}^*\right\|^2}{2}>0 \quad\forall \boldsymbol{v}^*, \boldsymbol{B}^* \in
\mathbb{R}^{3},
$$
which implies $\bar{\boldsymbol{U}}_{ij}^{n+1}\in\mathcal{G}$ by Lemma \ref{Lem:admistatequiv}.
The proof is completed. \qed

As direct consequences of Theorem \ref{Thm:2DPPpropertytheorem}, we have the following two corollaries.
\begin{corollary}[PP via Zhang--Shu convex decomposition]\label{Cor:2DzhangshuCFL}
	Let the OEDG solution polynomial vectors $\left\{ \widetilde{\boldsymbol{U}}_{ij}^{n}\left( \boldsymbol{x}\right)\right\}$ be LDF and satisfy the following condition after the PP limiting procedure: 
	$$
	\widetilde{\boldsymbol{U}}_{ij}^{n}\left(\boldsymbol{x}\right)\in\mathcal{G}
	\quad \forall \boldsymbol{x}\in \mathbb{S}_{ij}^{\text{classic}},
	$$
	where $\mathbb{S}_{ij}^{\text{classic}}$ denotes the set of all nodes of the classic Zhang--Shu convex decomposition (\ref{Equ:2DZhangShuCD}). 
	If the wave speeds in the HLL flux satisfy the condition (\ref{Equ:2DPPwavespeedcon}), then  the scheme (\ref{Equ:updatecellave}) preserves $\bar{\boldsymbol{U}}_{ij}^{n+1}\in\mathcal{G}$ under the CFL-type condition
	\begin{equation*}\label{Equ:2DCFLzhangshucondi}
		\begin{aligned}
			&\Delta t_{n} \leq 
			\frac{\Delta x\Delta y \hat{\omega}_{1}^{GL}}{a_{2}\Delta x+a_{1}\Delta y},
		\end{aligned}
	\end{equation*}
	where $a_1 = \max \{  \alpha_{1,i+\frac{1}{2}}^{\mu}\}$,
	$a_2 = \max \{  \alpha_{2,j+\frac{1}{2}}^{\mu} \}$, and 
	$\hat{\omega}_{1}^{GL}=\frac{1}{L\left(L-1\right)}$ with $L=\lceil\frac{k+3}{2}\rceil$.
\end{corollary}

\noindent \textbf{\emph{Proof }}
	This follows from Theorem \ref{Thm:2DPPpropertytheorem} by taking the convex decomposition as the Zhang--Shu convex decomposition (\ref{Equ:2DZhangShuCD}) with $a_1 = \max \{  \alpha_{1,i+\frac{1}{2}}^{\mu}\}$ and 
	$a_2 = \max \{  \alpha_{2,j+\frac{1}{2}}^{\mu} \}$. 
\qed

\begin{corollary}[PP via optimal convex decomposition for $\mathbb P^2$ and $\mathbb P^3$]\label{Cor:2DoptimalCFL}
	Consider the $\mathbb P^2$- and $\mathbb P^3$-based LDF OEDG method. 
		Let the OEDG solution polynomial vectors $\left\{ \widetilde{\boldsymbol{U}}_{ij}^{n}\left( \boldsymbol{x}\right)\right\}$ be LDF and satisfy the following condition after the PP limiting procedure: 
			$$
		\widetilde{\boldsymbol{U}}_{ij}^{n}\left(\boldsymbol{x}\right)\in\mathcal{G}
		\quad \forall \boldsymbol{x}\in \mathbb{S}_{ij}^{\text{optimal}},
		$$
		where $\mathbb{S}_{ij}^{\text{optimal}}$ denotes the set of all nodes of optimal convex decomposition (\ref{Equ:2DoptimalCD}). 
	If the wave speeds in the HLL flux satisfy the condition (\ref{Equ:2DPPwavespeedcon}),
	then the scheme (\ref{Equ:updatecellave}) preserves $\bar{\boldsymbol{U}}_{ij}^{n+1}\in\mathcal{G}$ under the CFL-type condition
	\begin{equation*}\label{Equ:2DCFLoptimalcondi}
		\begin{aligned}
			&\Delta t_{n} \leq 
			\frac{\Delta x \Delta y}{\max\left\{ 6a_{1}\Delta y+2a_{2}\Delta x,  6a_{2}\Delta x+2a_{1}\Delta y\right\}},
		\end{aligned}
	\end{equation*}
where $a_1 = \max \{  \alpha_{1,i+\frac{1}{2}}^{\mu}\}$ and 
$a_2 = \max \{  \alpha_{2,j+\frac{1}{2}}^{\mu} \}$. 
\end{corollary}

\noindent \textbf{\emph{Proof }}
	This follows from Theorem \ref{Thm:2DPPpropertytheorem} by adopting the optimal convex decomposition (\ref{Equ:2DoptimalCD}) with $a_1 = \max \{  \alpha_{1,i+\frac{1}{2}}^{\mu}\}$ and 
	$a_2 = \max \{  \alpha_{2,j+\frac{1}{2}}^{\mu} \}$. 
\qed 

\begin{remark}
	For the $\mathbb P^k$-based LDF OEDG method with higher $k\ge 4$, the PP CFL condition derived from the optimal convex  decomposition can be expressed as 
	 $$
	 \Delta t_n \left(\frac{a_1}{\Delta x} + \frac{a_2}{\Delta y}\right) \le  \bar \omega_\star( \delta, \mathbb P^k ),
	 $$
	 where weight $\bar \omega_\star=\bar \omega_\star( \delta, \mathbb P^k )$ depends on $\delta :=  \frac{a_1 \Delta y - a_2 \Delta x}{ a_1 \Delta y + a_2 \Delta x } $ and the space $\mathbb P^k$. 
	 For example, $\bar \omega_\star( \delta, \mathbb P^1 ) = \frac12$, 
	 $\bar \omega_\star( \delta, \mathbb P^2 )= \bar \omega_\star( \delta, \mathbb P^3 )=\frac{1}{4+2|\delta|}$, and 
	 	for $k\in \{4,5\}$, 
	 	$$\bar \omega_\star( \delta, \mathbb P^k )= \left[\frac{14}{3}+\frac{2}{3}\sqrt{78\,\delta^2+46} \cos \left( \frac{1}{3} \arccos\frac{1476\,\delta^2-244}{(78\,\delta^2+46)^{\frac{3}{2}}} \right) \right]^{-1}.$$ 
	 For higher degree $k$; see \cite{Cui2024On} for more details. 
\end{remark}

Thanks to the ``two-state'' inequalities established in Lemma \ref{Thm:2DPPstatetheorem} and the optimal convex decomposition, our above PP CFL condition \eqref{Equ:2DCFLoptimalcondi} is notably improved, compared to the estimates in \cite{Wu2019Provably} for general meshes when reduced to Cartesian meshes.

In our above analysis, our focus has been on employing the forward Euler time discretization.  Given that a high-order SSP method is essentially a convex combination of the forward Euler approach, our PP analysis of the proposed schemes persists as applicable for high-order SSP time discretization owing to the convexity $\mathcal{G}$.

\section{Numerical tests}\label{Sec:numericaltest}

This section gives several 1D and 2D MHD numerical examples to validate the accuracy, the essentially non-oscillatory shock-capturing capability, and the robustness of the proposed PP LDF OEDG method on uniform 1D meshes and 2D rectangular meshes. Without loss of generality, we focus on the third-order ($k=2$) OEDG method, coupled with the third-order SSP Runge--Kutta time discretization.
We use the HLL flux and set the CFL number to $0.12$.

\subsection{1D MHD tests}

This subsection assesses the OEDG scheme through various 1D MHD examples, including a smooth sine wave problem, three shock tube problems, and the MHD version of the Leblanc problem.

\subsubsection{1D smooth sine wave problem}

We begin by simulating a smooth MHD sine wave problem to evaluate the accuracy of the PP OEDG scheme in the 1D setting. This problem involves a sine wave propagating with low density, and its exact solution is described by:
$$
\left(\rho, \boldsymbol{u}, p, \boldsymbol{B}\right)\left(x, t\right)=
\left(1+0.99\sin\left(x-t\right), 1, 0, 0, 1, 0.1, 0, 0\right), \quad
x\in\left[0, 2\pi\right], t\geq 0.
$$
Periodic boundary conditions are applied, and $\gamma=1.4$.
Table \ref{Tab:1DsmoothsineMHD} lists the numerical errors and the corresponding convergence rates at $t=0.1$ for the density. The results confirm that the 1D $\mathbb P^2$-based PP OEDG scheme achieves the expected third-order convergence rate. This confirms that the OE procedure and PP limiter do not compromise the optimal convergence order of the DG solutions.

\begin{table}[!htb]  
  \centering
  \setlength{\abovecaptionskip}{0pt}%
  \setlength{\belowcaptionskip}{10pt}%
  \caption{Numerical errors and the corresponding convergence rates at $t=0.1$ for the density in the 1D smooth sine wave problem.}
  \label{Tab:1DsmoothsineMHD}
  \begin{tabular}{lllllll}
    \hline
     \text{Mesh} & $l^{1}$-error  & Order  & $l^{2}$-error  & Order & $l^{\infty}$-error  & Order\\
     \hline
     100 $\times$ 100 & 7.3204E-06 & - & 3.3012E-06 & - &
     2.2240E-06 & -  \\
     200 $\times$ 200 & 8.0642E-07 & 3.1823 & 3.5962E-07 & 3.1984 &
     2.3122E-07 & 3.2658  \\
     400 $\times$ 400 & 9.5388E-08 & 3.0796 & 4.2327E-08 & 3.0868 &
     2.6784E-08 & 3.1098  \\
     800 $\times$ 800 & 1.1621E-08 & 3.0371 & 5.1440E-09 & 3.0406 &
     3.2338E-09 & 3.0500  \\
     1600 $\times$ 1600 &  1.4340E-09 & 3.0186 & 6.3399E-10 & 3.0204 &
     3.9760E-10 & 3.0239  \\
     3200 $\times$ 3200 &  1.7806E-10 &  3.0096 & 7.8681E-11 & 3.0104 &
     4.9302E-11 & 3.0116  \\
     \hline
  \end{tabular}
\end{table}

\subsubsection{Shock tube problems}

\begin{figure}[!htb]
  \centering
  \includegraphics[scale=0.8]{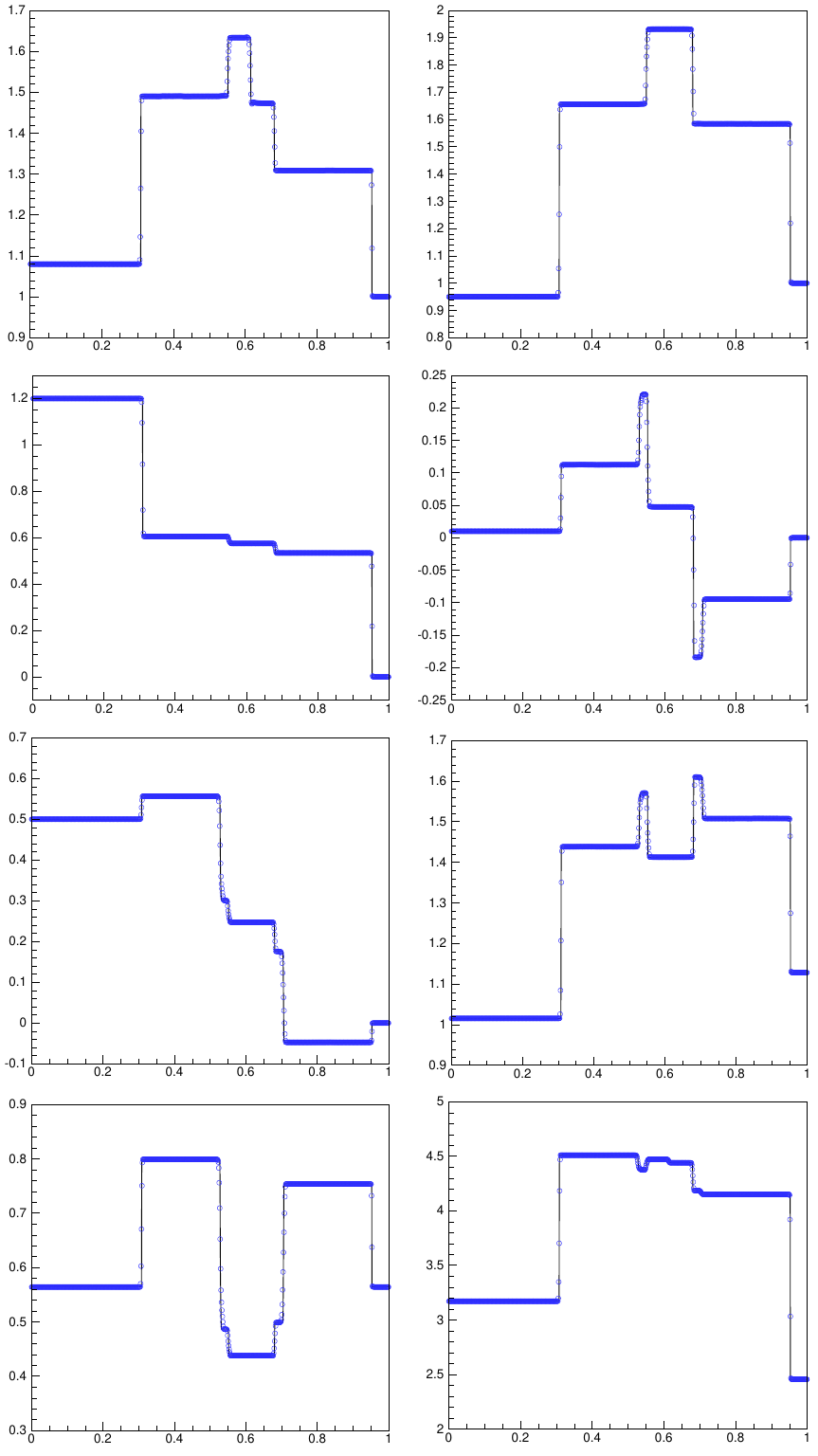}
  \caption{Numerical solutions of the first MHD shock tube problem at time $t = 0.2$ with $800$ cells (symbols ``$\circ$") and $4000$ cells (solid lines).
  Left column from top to bottom: $\rho$, $u_{1}$, $u_{3}$, $B_{3}$.
  Right column from top to bottom: $p$, $u_{2}$, $B_{2}$, $E$.
  }
  \label{Fig:1Dshocktube}
\end{figure}

\begin{figure}[!htb]
  \centering
  \includegraphics[scale=0.8]{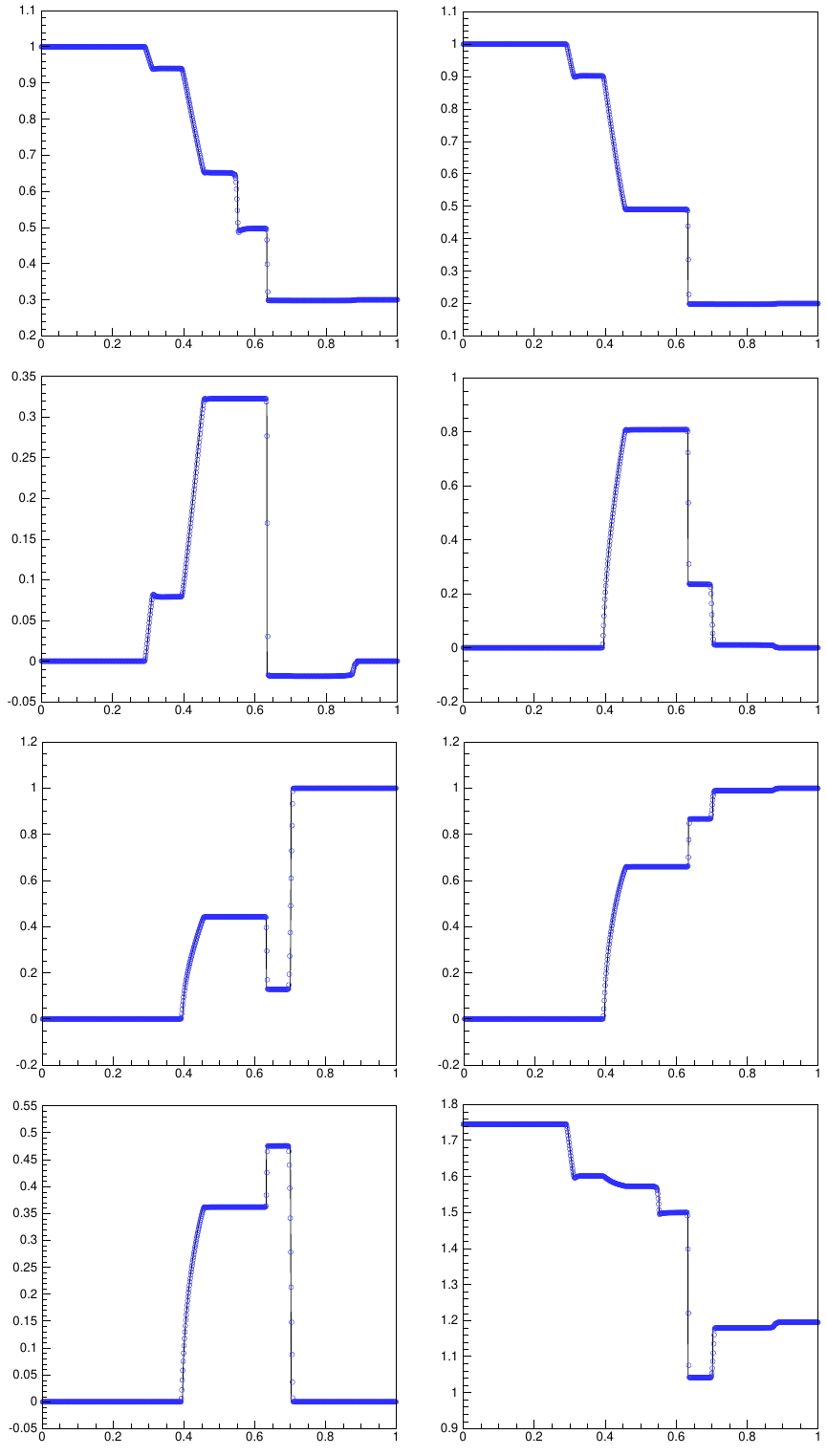}
  \caption{Numerical solutions of the second MHD shock tube problem at time $t = 0.16$ with $800$ cells (symbols ``$\circ$") and $4000$ cells (solid lines).
  Left column from top to bottom: $\rho$, $u_{1}$, $u_{3}$, $B_{3}$.
  Right column from top to bottom: $p$, $u_{2}$, $B_{2}$, $E$.
  }
  \label{Fig:shocktube1D2}
\end{figure}

\begin{figure}[!htb]
  \centering
  \includegraphics[scale=0.8]{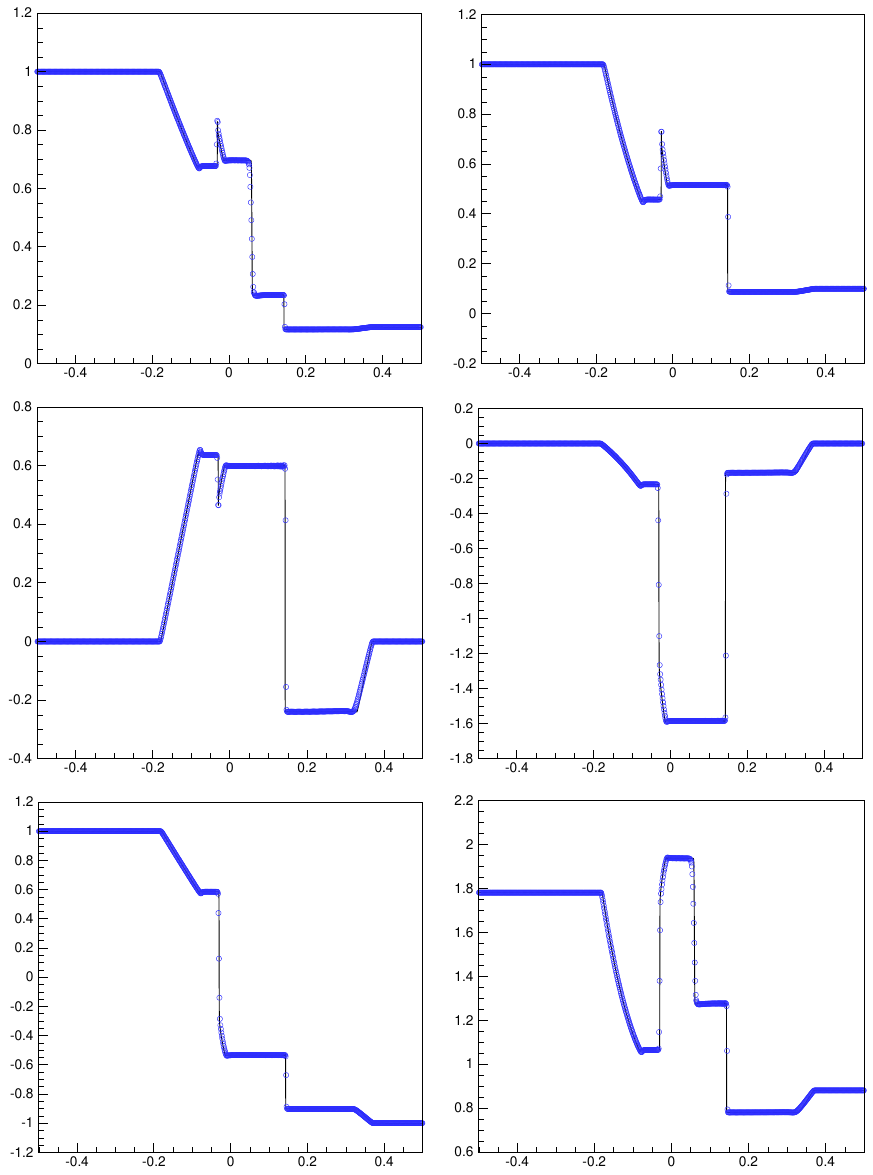}
  \caption{Numerical solutions of the third MHD shock tube problem at time $t = 0.1$ computed by the third-order OEDG scheme with $800$ cells (symbols ``$\circ$") and $4000$ cells (solid lines).
  Left column from top to bottom: $\rho$, $u_{1}$, $B_{2}$.
  Right column from top to bottom: $p$, $u_{2}$, $E$.
  }
  \label{Fig:1Dshocktube4}
\end{figure}

To evaluate the shock-capturing and non-oscillatory properties of the PP OEDG method, we simulate three 1D MHD shock tube problems. For the first two problems, the adiabatic index $\gamma$ is set to $5/3$, and for the third problem, it is set to $2$. All simulations employ the third-order PP OEDG scheme with $800$ uniform cells.
The initial conditions of the first shock tube problem \cite{Ryu1995Numerical}
 are  given by
$$
\left(\rho, \boldsymbol{u}, \boldsymbol{B}, p\right)\left(x, 0\right)=
\left\{
    \begin{aligned}
    &\left(1.08, 1.2, 0.01, 0.5, \frac{2}{\sqrt{4\pi}}, \frac{3.6}{\sqrt{4\pi}}, \frac{2}{\sqrt{4\pi}}, 0.95\right), & x<0.5,\\
    &\left(1, 0, 0, 0, \frac{2}{\sqrt{4\pi}}, \frac{4}{\sqrt{4\pi}}, \frac{2}{\sqrt{4\pi}}, 1\right), & x>0.5,\\
    \end{aligned}
\right.
$$
and the numerical results at time $t=0.2$ are presented in Figure \ref{Fig:1Dshocktube}.
The second shock tube problem encompasses a hydrodynamic rarefaction, a switch-on slow shock, a contact discontinuity, a slow shock, a rotational discontinuity, and a fast rarefaction, following the setup in \cite{Ryu1995Numerical}:
$$
\left(\rho, \boldsymbol{u}, \boldsymbol{B}, p\right)\left(x, 0\right)=
\left\{
    \begin{aligned}
    &\left(1, 0, 0, 0, 0.7, 0, 0, 1\right), & x<0.5,\\
    &\left(0.3, 0, 0, 1, 0.7, 1, 0, 0.2\right), & x>0.5,\\
    \end{aligned}
\right.
$$
and
Figure \ref{Fig:shocktube1D2} shows the numerical results at time $t=0.16$.
The third shock tube problem, initially proposed by Brio and Wu \cite{BRIO1988400}, features the following initial conditions:
$$
\left(\rho, \boldsymbol{u}, \boldsymbol{B}, p\right)\left(x, 0\right)=
\left\{
    \begin{aligned}
    &\left(1, 0, 0, 0, 0.75, 1, 0, 1\right), & x<0,\\
    &\left(0.125, 0, 0, 0, 0.75, -1, 0, 0.1\right), & x>0.\\
    \end{aligned}
\right.
$$
The numerical results at time $t=0.1$ are shown in Figure \ref{Fig:1Dshocktube4}.
From these results, one can see that the PP OEDG method accurately captures
all the discontinuities such as shocks with high resolution and  without any noticeable spurious oscillations.

\subsubsection{Leblanc problem}

\begin{figure}[!htb]
  \centering
  \includegraphics[scale=1.0]{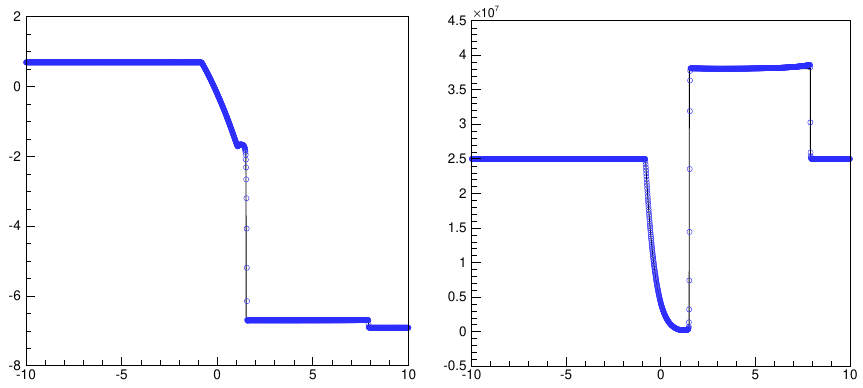}
  \caption{Density logarithm (left) and magnetic pressure (right) of the 1D Leblanc problem at time $t=0.00003$ with $2000$ cells (symbols ``$\circ$") and $10000$ cells (solid lines).
  }
  \label{Fig:1DLeblanc}
\end{figure}

To assess the robustness of the 1D PP OEDG scheme, we test a challenging MHD Riemann problem with a strong magnetic field, a large pressure jump, and an extremely low plasma-beta ($\beta=4\times10^{-8}$).
This problem
is a variant of the Leblanc problem \cite{ZHANG20108918} of the gas
dynamics and was proposed in \cite{Wu2019Provably}.
The initial conditions are defined as follows:
$$
\left(\rho, \boldsymbol{u}, p, \boldsymbol{B}\right)\left(x, 0\right)=
\left\{
    \begin{aligned}
    &\left(2, 0, 0, 0, 10^9, 0, 5000, 5000\right), & x<0,\\
    &\left(0.001, 0, 0, 0, 1, 0, 5000, 5000\right), & x>0.\\
    \end{aligned}
\right.
$$
The computational domain spans $\left[-10, 10\right]$, and the adiabatic index $\gamma$ is set to $1.4$. Figure \ref{Fig:1DLeblanc} presents the numerical results for density and magnetic pressure at $t=0.00003$. As anticipated, no negative density or pressure values are observed in the OEDG solution, demonstrating the robustness and stability of the 1D PP OEDG scheme.

\subsection{2D MHD tests}

This subsection conducts six 2D MHD tests on the third-order PP LDF OEDG method, which include simulations of a smooth sine wave, the Orszag--Tang vortex, a rotor problem, shock-cloud interaction, an MHD blast wave, and MHD jet flows.

\subsubsection{2D smooth sine wave problem}

To evaluate the convergence order of the PP LDF OEDG method, we simulate a 2D smooth sine wave problem \cite{Wu2018A}, which involves an MHD sine wave propagating through a medium with low density. The exact solution of this problem is expressed as follows:
$$
\left(\rho, \boldsymbol{u}, p, \boldsymbol{B}\right)\left(x, y, t\right)=
\left(1+0.99\sin\left(x+y-2t\right), 1, 1, 0, 1, 0.1, 0.1, 0\right), \quad
x\in\left[0, 2\pi\right], y\in\left[0, 2\pi\right], t\geq 0,
$$
and the adiabatic index $\gamma=1.4$.
Table \ref{Tab:2DsineMHD} presents the numerical errors and corresponding convergence rates for density at time $t=0.1$. The results demonstrate that the $\mathbb{P}^{2}$-based PP LDF OEDG scheme consistently achieves third-order accuracy. That is, the optimal convergence rate is not destroyed by the LDF OE procedure and the PP limiter.

\begin{table}[!htb]  
  \centering
  \caption{Numerical errors and corresponding convergence orders at $t=0.1$ for density in the 2D smooth sine wave problem.}
  \label{Tab:2DsineMHD}
  \begin{tabular}{lllllll}
    \hline
     \text{Mesh} & $l^{1}$-error  & Order  & $l^{2}$-error  & Order & $l^{\infty}$-error  & Order\\
     \hline
     15 $\times$ 15 & 3.7023E-01 & - & 8.3500E-02 & - &
     7.2613E-02 & -  \\
     30 $\times$ 30 & 4.8828E-02 & 2.9227 & 1.1535E-02 & 2.8557 &
     1.0510E-02 & 2.7885  \\
     60 $\times$ 60 & 2.5913E-03 & 4.2360 & 6.1398E-04 & 4.2317 &
     4.8343E-04 & 4.4423  \\
     120 $\times$ 120 & 1.0602E-04 & 4.6113 & 2.4135E-05 & 4.6690 &
     2.3836E-05 & 4.3421  \\
     240 $\times$ 240 &  9.0494E-06 & 3.5503 & 2.1438E-06 & 3.4929 &
     2.3849E-06 & 3.3212  \\
     480 $\times$ 480 &  9.8980E-07 &  3.1926 & 2.3377E-07 & 3.1970 &
     2.6805E-07 & 3.1533  \\
     \hline
  \end{tabular}
\end{table}

\subsubsection{Orszag--Tang problem}

We next examine the benchmark Orszag--Tang vortex problem \cite{Orszag1979Small}. This problem begins with smooth initial data, rapidly forming shocks, and eventually transitioning into 2D MHD turbulence. The initial conditions are set as
$$
\begin{array}{l}
  \rho\left(x, y\right)=\gamma^{2}, \quad u_{1}\left(x, y\right)=-\sin y, \quad
  u_{2}\left(x, y\right)= \sin x,  \quad u_{3}\left(x, y\right)=0,\\
  p\left(x, y\right)=\gamma, \quad B_{1}\left(x, y\right)=-\sin y, \quad
  B_{2}\left(x, y\right)=\sin \left(2x\right), \quad B_{3}\left(x, y\right)=0,
\end{array}
$$
where the adiabatic index $\gamma$ is $5/3$. The computational domain spans $\left[0, 2\pi\right]\times\left[0, 2\pi\right]$, with periodic boundary conditions applied on all edges.

\begin{figure}[!htb]
  \centering
  \includegraphics[scale=0.9]{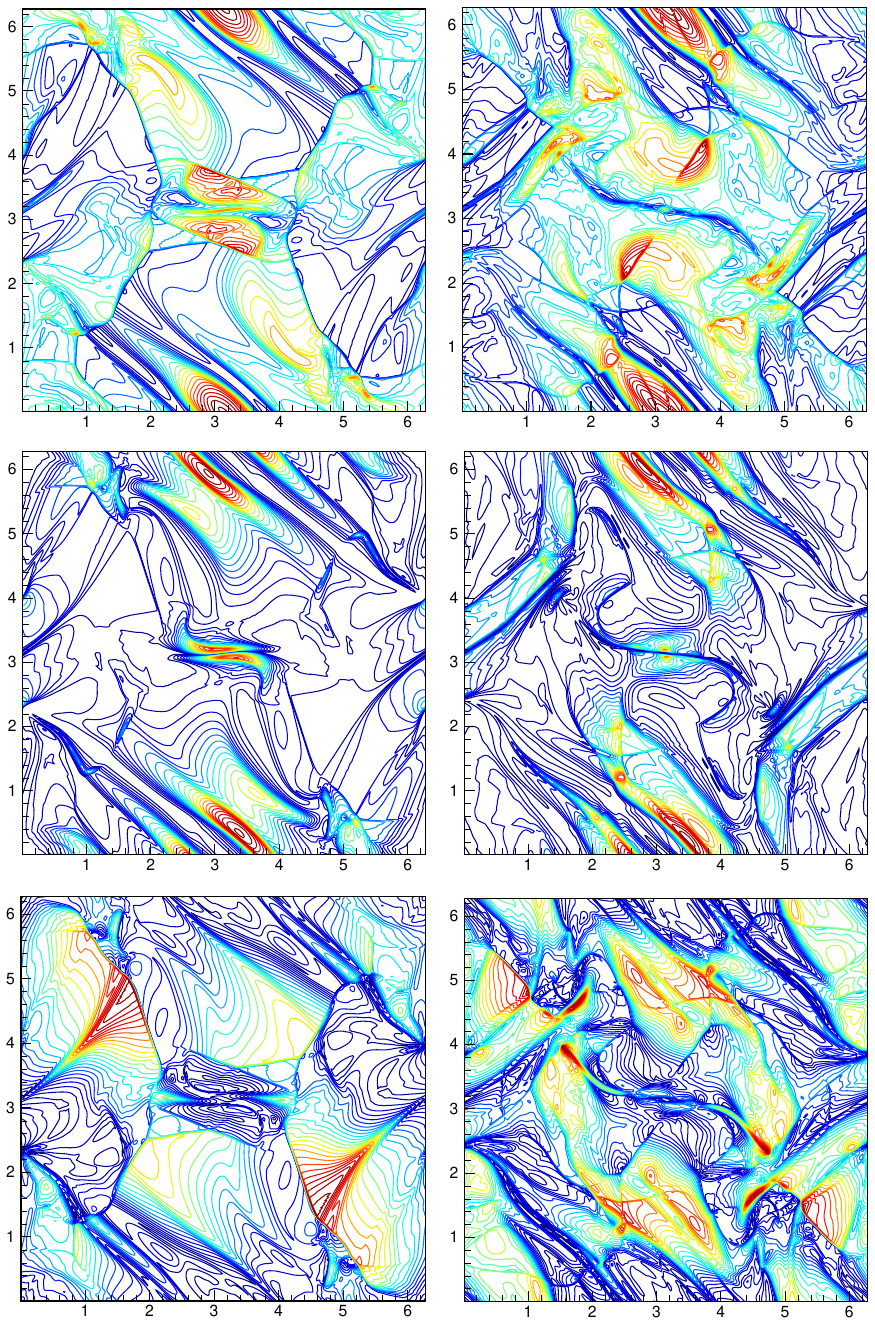}
  \caption{Contour plots for the Orszag--Tang problem with $400\times400$ cells.
  From top to bottom: density $\rho$, magnetic pressure
  $\|\boldsymbol{B}\|^{2}/2$, and Mach number $\|\boldsymbol{u}\|/c$.
  Left: $t=3$; right: $t=4$.}
  \label{Fig:OTcontour}
\end{figure}

\begin{figure}[!htb]
  \centering
  \includegraphics[scale=1.0]{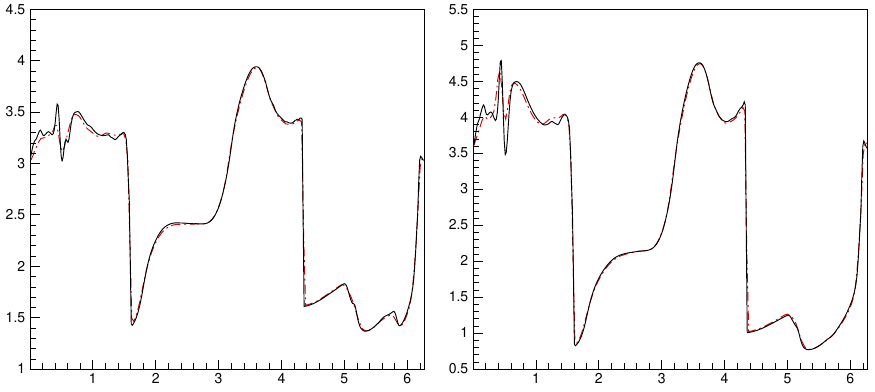}
  \caption{Profiles of density (left) and thermal pressure (right) along line $y=0.625\pi$ at time $t=3$ for the Orszag--Tang problem  with $200\times200$ cells (dash dot lines) and $400\times400$ cells (solid lines).}
  \label{Fig:OTprofile}
\end{figure}

Figure \ref{Fig:OTcontour} displays contour plots of density, magnetic pressure, and Mach number at times $t=3$ and $t=4$, obtained by the third-order PP LDF OEDG scheme with $400\times400$ cells.
Our method accurately resolves the shocks and smooth flows, which are
in good agreement with those in \cite{LI20114828,LI20122655,Guillet2019}.
For this problem, Jiang and Wu \cite{JIANG1999561} previously reported that  negative pressure was occurred at $t\approx3.9$ in their computation.
Thanks to the PP property of our scheme, negative pressure was not observed throughout our simulation.
To further verify the convergence of our scheme, Figure \ref{Fig:OTprofile} plots the
profiles of density and thermal pressure along the line of $y=0.625\pi$ at time $t=3$ under $200\times200$ and $400\times400$ resolutions.
As shown in Figure \ref{Fig:OTprofile}, the
shock discontinuities for density and thermal pressure are formed
near $x=0.4, 0.5, 1.6$ and $4.4$, which are consistent with the results in \cite{YANG2017561,LIU2021110694}.
As discussed by Yang et al.~in \cite{YANG2017561},
high-resolution simulations typically exhibit sharper shock profiles, particularly noticeable at $x=0.4$ and $x=0.5$.
In Figure \ref{Fig:OTprofile}, it is observed that the shock discontinuities around $x=0.4$ and $x=0.5$ with $400\times400$ resolution are much sharper than those with $200\times200$ resolution, indicating  the convergence and effectiveness of our scheme.

\subsubsection{Rotor problem}

We further evaluate the PP LDF OEDG method by simulating the rotor problem  \cite{Balsara1999A}. This problem involves a dense disk of fluid embedded in a static fluid background, with a velocity tapering layer between the dense disk's edge and the static ambient fluid. As time progresses, the disk rotates within the ambient fluid. The initial configurations are set as follows:
$$
\left(u_{3}, B_{1}, B_{2}, B_{3}, p\right)=(0,2.5 / \sqrt{4 \pi}, 0,0,0.5),
$$
and
$$
\left(\rho, u_{1}, u_{2}\right)=
\left\{
\begin{aligned}
&\left(10, -(y-0.5) / r_{0}, (x-0.5) / r_{0}\right) & \text { if } r<r_{0}, \\
&\left(1+9f, -f(y-0.5) / r, f(x-0.5) / r\right) & \text { if } r_{0}<r<r_{1}, \\
&(1, 0, 0) & \text { if } r>r_{1},
\end{aligned}
\right.
$$
where $r_{0}=0.1$, $r_{1}=0.115$, $r = \sqrt{(x-0.5)^2+(y-0.5)^2}$, and $f=(r_{1}-r)/(r_{1}-r_{0})$.
The adiabatic index $\gamma$ is $5/3$, and the computational domain is $\left[0, 1\right]\times\left[0, 1\right]$ with periodic boundary conditions.

\begin{figure}[!htb]
  \centering
  \includegraphics[scale=0.9]{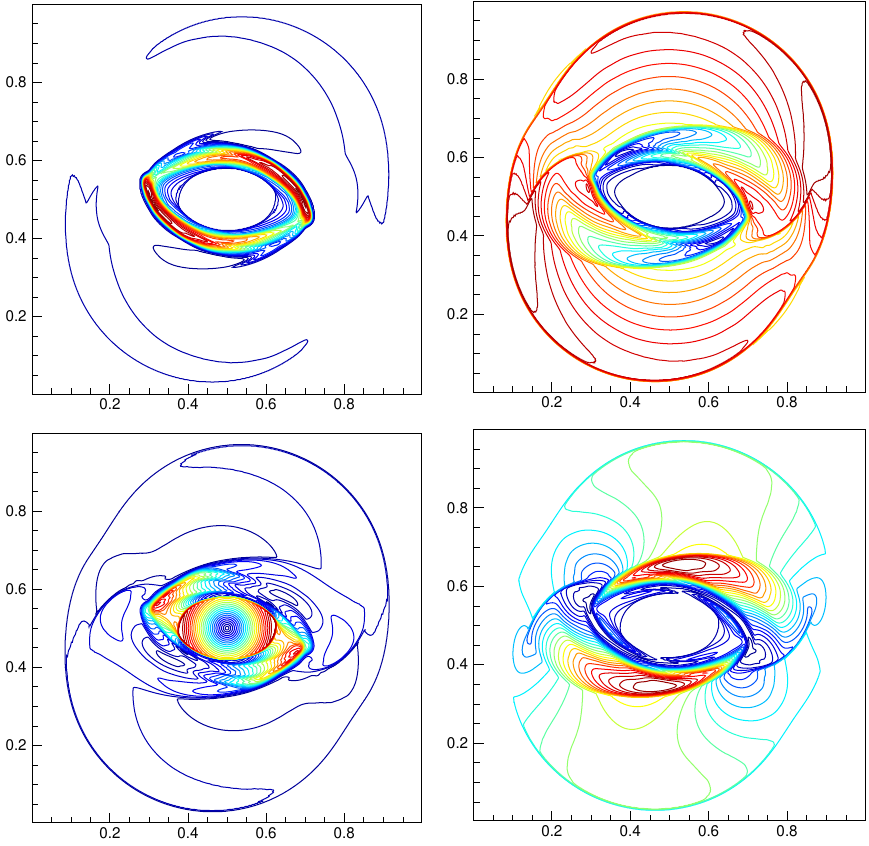}
  \caption{Contour plots of density $\rho$ (top left), thermal pressure $p$ (top right), Mach number $\|\boldsymbol{u}\|/c$ (bottom left) and magnetic pressure $\|\boldsymbol{B}\|^{2}/2$ (bottom right) at time $t=0.295$ for the rotor problem with $400\times400$ cells.
  Twenty-five equally spaced contours are used for density with range $[0.8, 10.2]$, thermal pressure with range $[0.01, 1.1]$, Mach number with range $[0, 2.2]$, and magnetic pressure with range $[0.02, 0.68]$.}
  \label{Fig:RTcontour}
\end{figure}

\begin{figure}[!htb]
  \centering
  \includegraphics[scale=1.0]{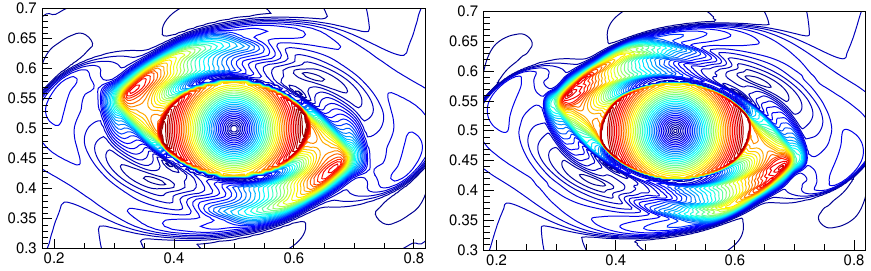}
  \caption{Zoom-in central parts of Mach number at time $t=0.295$ for the rotor problem with $200\times200$ cells (left) and
   $400\times400$ cells (right).
  Forty equally spaced contours are utilized within the range $[0, 2.2]$.
  }
  \label{Fig:RTmach}
\end{figure}

\begin{figure}[!htb]
  \centering
  \includegraphics[scale=1.0]{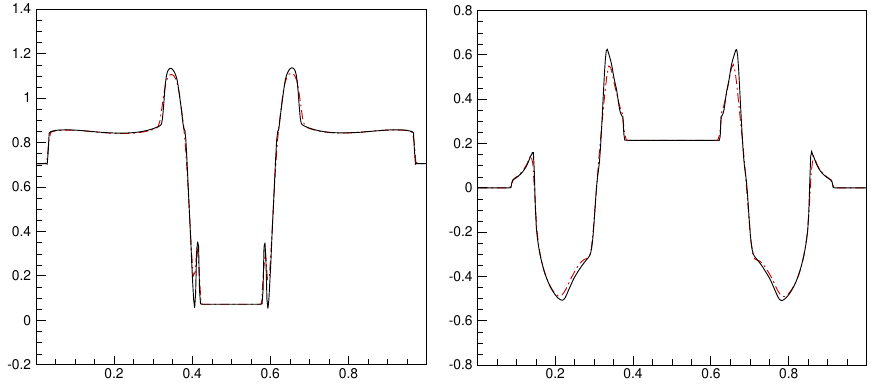}
  \caption{Slices of $B_{1}$ (left) and $B_{2}$ (right) along the lines $x=0.5$ and $y=0.5$ at time $t=0.295$ for the rotor problem with $200\times200$ cells (dash dot lines) and $400\times400$ cells (solid lines).}
  \label{Fig:RTprofile}
\end{figure}

Figure \ref{Fig:RTcontour} presents the contours of density, thermal pressure, Mach number, and magnetic pressure at $t=0.295$, computed using the third-order PP LDF OEDG method with $400\times400$ cells. The contours show good conservation of the circular rotation pattern, which T\'{o}th \cite{TOTH2000605}  found challenging for some other MHD numerical schemes. Distortions may appear in the numerical solution around the central part of the Mach number if there are large divergence errors in the magnetic field, as observed in \cite{Li2005}. 
However, our simulation results show no distortions around the central almond-shaped disk region, suggesting effective control of divergence errors by our LDF OEDG scheme. 
To further assess the performance of our scheme, Figure \ref{Fig:RTmach} zooms into the central part of our computed Mach number at different resolutions, showing no observed distortions. Figure \ref{Fig:RTprofile} presents two slices of the magnetic field at time $t=0.295$ across horizontal and vertical cuts through the center of the domain at different mesh resolutions, with results that align well with those reported by  \cite{LIU2021110694,LIU2023111687}.
It was found that in \cite{LIU2021110694} higher-order scheme exhibits higher peaks near $x=0.35, 0.65$ of $B_{1}$ and $x=0.325, 0.675$ of $B_{2}$, as well as sharper shock discontinuities around $x=0.4, 0.6$ of $B_{1}$.
In Figure \ref{Fig:RTprofile}, the same trend is observed in our high-resolution results compared to the low-resolution ones, reflecting the convergence behavior of our scheme.
According to Liu et al.~\cite{LIU2021110694},
the overshoots at $y=0.05, 0.95$ for $B_{1}$ and at $x=0.38, 0.62$ for $B_{2}$
were generated by their third-order finite volume scheme.
Our results in Figure \ref{Fig:RTprofile} show no overshoots, further verifying that the LDF OE procedure effectively suppresses nonphysical oscillations.

\subsubsection{Shock-cloud interaction}

To further validate the shock-capturing and PP capabilities
of our structure-preserving OEDG scheme, we simulate the shock-cloud interaction problem \cite{DAI1998331}.
This test involves strong MHD shocks interacting with a dense cloud, which is disrupted by the interaction. The initial conditions feature left and right states separated by a discontinuity along the line $x=0.6$:
$$
\left(\rho, \mathbf{u}, p, \boldsymbol{B}\right)=
\left\{
\begin{aligned}
&\left(3.86859,0,0,0,167.345,0,2.1826182,-2.1826182\right), & x<0.6, \\ &\left(1,-11.2536,0,0,1,0,0.56418958,0.56418958\right), & x>0.6,\\
\end{aligned}
\right.
$$
and includes a circular cloud with a radius of $0.15$ centered at $(0.8, 0.5)$  on the right side of the discontinuity. The cloud has the same state
as the surrounding medium but with a higher density of $10$. The computational domain is $[0, 1]\times[0,1]$, and adiabatic index $\gamma$ is $5/3$.
For this problem, supersonic inflow condition is employed on the right boundary, while outflow conditions are utilized on the others.

\begin{figure}[!htb]
  \centering
  \includegraphics[scale=0.9]{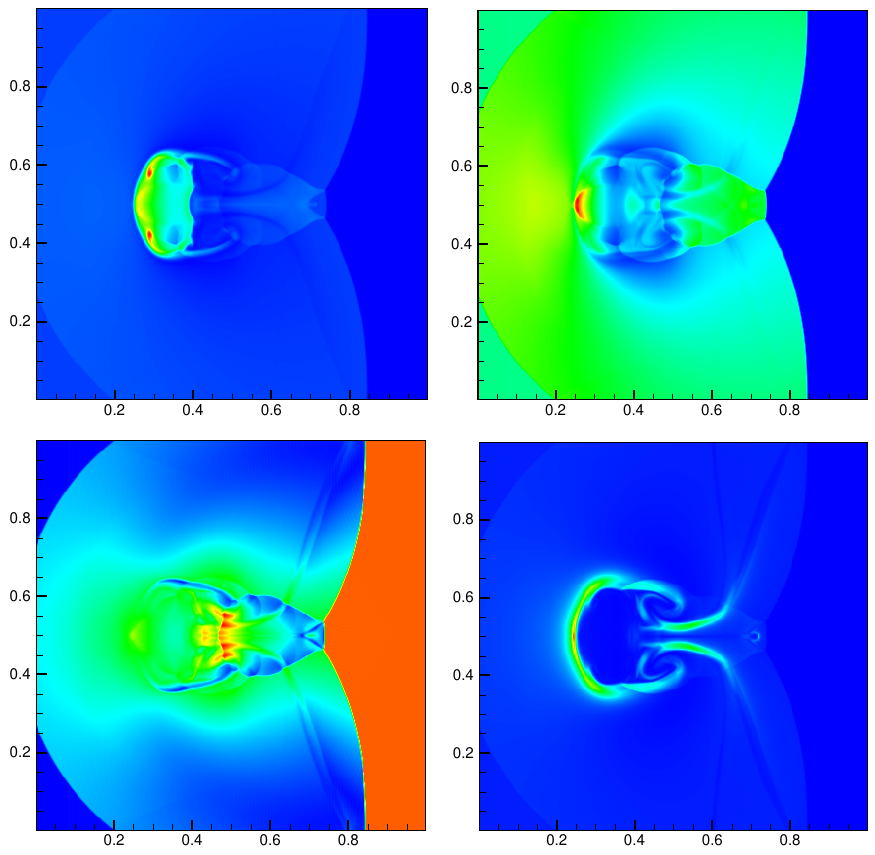}
  \caption{The density $\rho$ (top left), thermal pressure $p$ (top right), velocity $\|\boldsymbol{u}\|$ (bottom left) and magnetic pressure $\|\boldsymbol{B}\|^{2}/2$ (bottom right) at time $t=0.06$ for the shock cloud interaction with $400\times400$ cells.}
  \label{Fig:SCcontour}
\end{figure}

Figure \ref{Fig:SCcontour} displays the numerical results at time $t=0.06$ for density, thermal pressure, velocity magnitude, and magnetic pressure, computed using our third-order PP LDF OEDG scheme on a $400\times400$ mesh. The results successfully capture complex flow structures and discontinuities at high resolution, demonstrating good agreement with those in the literature, e.g., \cite{TOTH2000605,Balbas2006Non,Wu2018A,Wu2019Provably}.
We note that disabling the LDF OE procedure in this simulation leads to failure, underscoring the critical role of our OE technique. Furthermore, removal of either the Godunov--Powell source terms or the PP limiter resulted in negative pressure values in the cell averages of numerical solutions, reaffirming the necessity of these elements for ensuring positivity in the simulation.

\subsubsection{MHD blast wave}

We now turn to a 2D blast wave test to further verify the capability of our structure-preserving OEDG scheme in handling strong shocks and rarefaction waves. This test involves an over-pressured region situated at the center of a strongly magnetized medium with low plasma $\beta$. The scenario leads to the formation of an MHD blast wave that drives fast outward magnetic shocks, compressing both the plasma and magnetic fields ahead, while internally generating a rarefaction wave.
We utilize the setup in \cite{Stone2008A}, with the initial conditions specified as
$$
 \rho=1,~  u_1=u_2=u_3=0, ~ B_1=B_2=\frac1 { \sqrt{2}},~  B_3=0,~
 p=\left\{
\begin{aligned}
&10, & \text { if } \sqrt{(x-1.5)^2+y^2}<0.1, \\
&0.1, & \text { otherwise}.\\
\end{aligned}
\right.
$$
The adiabatic index $\gamma$ is set to $5/3$, and Dirichlet boundary conditions matching the initial values are applied. The simulation is conducted over the domain $\left[1, 2\right]\times\left[-0.5, 0.5\right]$ and runs until a final time of $t=0.2$.

\begin{figure}[!htb]
  \centering
  \includegraphics[scale=0.9]{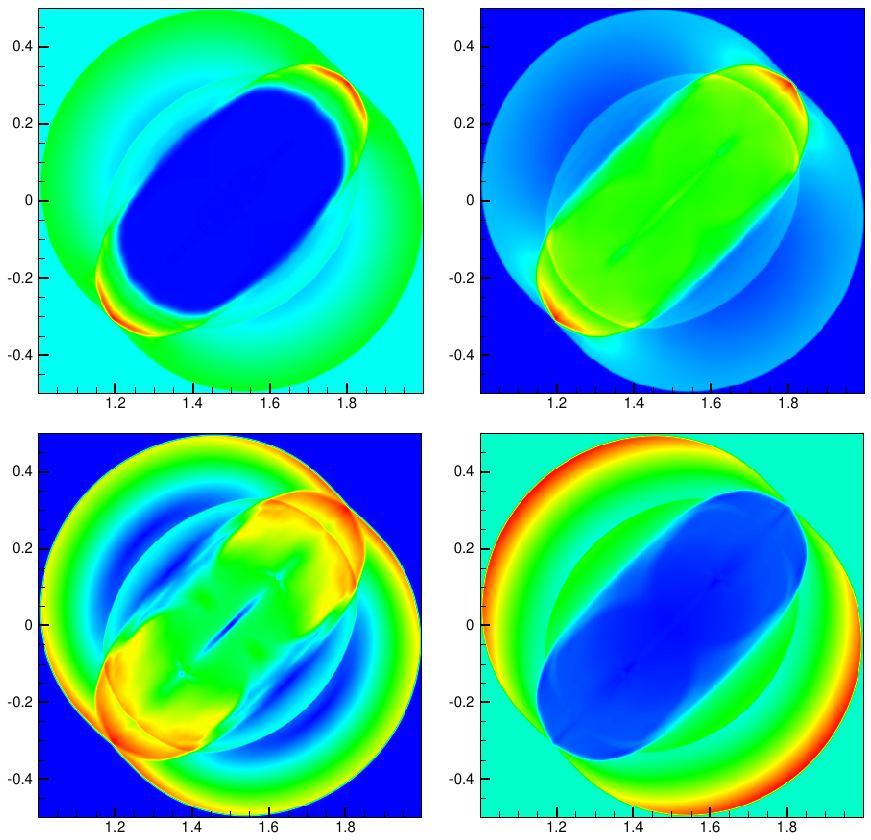}
  \caption{The density $\rho$ (top left), thermal pressure $p$ (top right), velocity $\|\boldsymbol{u}\|$ (bottom left) and magnetic pressure $\|\boldsymbol{B}\|^{2}/2$ (bottom right) at time $t=0.2$ for the MHD blast wave problem with $400\times400$ cells.}
  \label{Fig:BLcontour}
\end{figure}

Figure \ref{Fig:BLcontour} presents the results at $t=0.2$ computed using our third-order PP LDF OEDG scheme on a $400\times400$ mesh, showing the density, thermal pressure, velocity, and magnetic pressure. The results align well with those found in  \cite{Jiang2010AMR,YANG2017561}.
The visualizations clearly illustrate that a circular blast wave propagates outward, while a rarefaction wave moves inward. The numerical solution maintains symmetry excellently.
Our method successfully captures the shocks with high resolution and remains free from nonphysical oscillations, demonstrating its robustness and reliability in handling highly magnetized shock configurations. This performance underscores the capability of our structure-preserving OEDG scheme to effectively manage complex dynamic interactions in MHD environments.

\subsubsection{MHD jet flows}

In our final example, we explore two challenging MHD jet problems to assess our scheme's capability in capturing shocks and preserving positivity under demanding conditions.
These MHD jet flows were proposed in \cite{Wu2018A} by adding a magnetic field in the gas dynamical jet of \cite{BALSARA2012Self}.
Since there are strong shock waves, shear flows and interface instabilities in high Mach number jet with strong magnetic field and
huge kinetic energy, negative pressure is very likely to be produced in the numerical simulations.

Our computational domain is $[-0.5, 0.5]\times[0, 1.5]$, and the adiabatic index is $\gamma=1.4$.
Initially, we set $(\rho, p)=(0.1\gamma, 1)$ in the static ambient medium with a magnetic field $(0, B_{2}, 0)$.
As observed in \cite{Wu2018A}, the larger $B_{2}$ is set, the more challenging the test becomes.
On the bottom boundary, a dense jet described by $(\rho, p, u_{1}, u_{2}, u_{3})=(\gamma, 1, 0, u_{2}, 0)$ is injected through the nozzle $\{y=0, |x|<0.05\}$, with the Mach number determined by $u_2$. We employ fixed inflow conditions at the nozzle and outflow conditions on other boundaries. Reflecting boundary conditions are used at $x=0$, with the domain limited to $[0, 0.5] \times [0, 1.5]$ on the right half, divided into $200 \times 600$ cells.

\begin{figure}[!htb]
  \centering
  \includegraphics[scale=0.85]{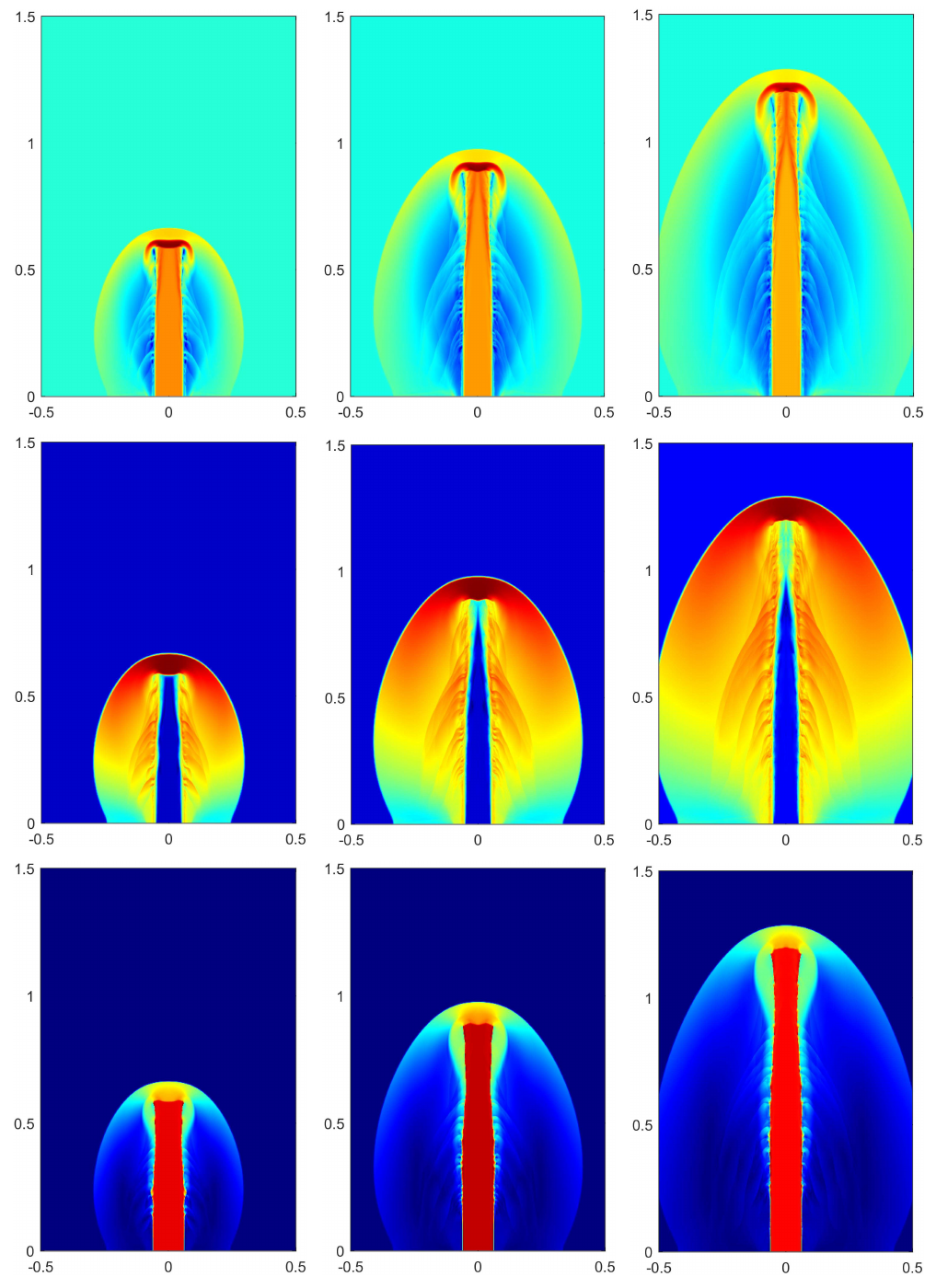}
  \caption{The density logarithm (top), thermal pressure logarithm (middle) and velocity (bottom) for the Mach 800 jet problem with $B_{2}=\sqrt{2000}$. From left to right: $t=0.001, 0.0015, 0.002$.}
  \label{Fig:jet2contour}
\end{figure}

First, we simulate a Mach 800 dense jet in a magnetized environment with $B_2 = \sqrt{2000}$. The numerical results, showcasing density logarithm $\log_{10}\rho$, thermal pressure logarithm $\log_{10}p$, and velocity magnitude $\|\boldsymbol{u}\|$, are displayed in Figure \ref{Fig:jet2contour}. These results vividly depict the evolution of jet flow structures, consistent with those reported in \cite{Wu2018A,Wu2019Provably}.
Our scheme accurately captures the Mach shock wave at the jet head and the beam/cocoon interface with high resolution.
Specifically, the shear flows near the beam/cocoon interfaces of jets are clearly presented.

\begin{figure}[!htb]
  \centering
  \includegraphics[scale=0.85]{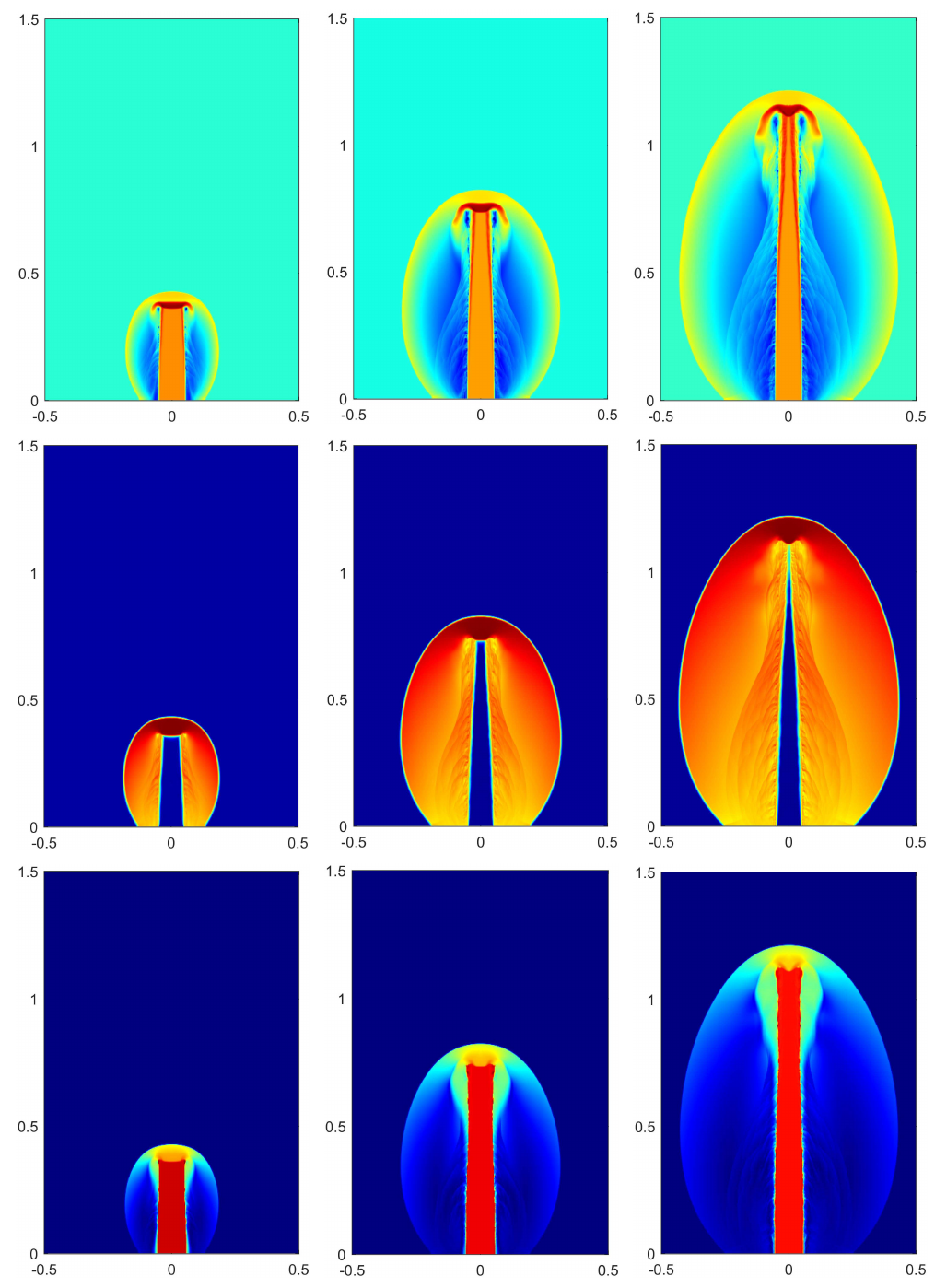}
  \caption{The density logarithm (top), thermal pressure logarithm (middle), and velocity (bottom) for the Mach 10000 jet problem with $B_{2}=\sqrt{20000}$. From left to right: $t=0.00005, 0.0001, 0.00015$.}
  \label{Fig:jet3contour}
\end{figure}

We also conduct a more demanding test with a Mach 10000 dense jet and a stronger magnetic field characterized by $B_2 = \sqrt{20000}$. Figure \ref{Fig:jet3contour} illustrates the time resolution of flow structures, including $\log_{10}\rho$, $\log_{10}p$, and $\|\boldsymbol{u}\|$. Throughout all simulations, no instances of negative pressure or density were observed, highlighting our scheme's robustness even under extreme conditions. However, if we disable the LDF OE procedure, the PP limiting procedure, or the upwind discrete Godunov--Powell source term, the code immediately fails in these challenging tests due to negative pressure in the numerical solution. This underscores the essential roles of the LDF OE procedure, PP limiting procedure, and upwind discrete Godunov--Powell source terms in ensuring the PP property and overall robustness of our scheme.

\section{Concluding remarks}\label{Sec:conclude}

In this paper, we have proposed a structure-preserving, oscillation-eliminating discontinuous Galerkin (OEDG) method for ideal magnetohydrodynamics (MHD). 
This method stands out for its capability to effectively suppresses potential spurious oscillations in the DG solution while preserving the key physical structures, including the divergence-free constraint of  magnetic field and the positivity of density and pressure. 
Based on a novel damping equation, the locally divergence-free (LDF)  
 OE procedure is designed to eliminate spurious oscillations 
while automatically maintaining an LDF magnetic field. This procedure is implemented after each Runge--Kutta stage in a non-intrusive way, and thus can be easily integrated into existing DG codes as an independent module. Its integration into the LDF DG framework maintains beneficial properties such as conservation, local compactness, and optimal convergence rates. 
A rigorous positivity-preserving (PP) analysis has been conducted for the LDF OEDG method with the HLL flux. This analysis utilizes the geometric quasi-linearization (GQL) approach to transform the nonlinear constraint of pressure positivity into a set of equivalent linear constraints. It has been proven that the LDF OEDG method is PP if an ``upwind" discrete Godunov--Powell source term is added solely to the evolution equations of the cell averages, under a condition achievable through a simple PP limiter. 
We have derived the PP CFL condition via a general convex decomposition of the cell averages and have also obtained the sharp PP CFL constraint through the optimal convex decomposition. 
The effectiveness of our method is validated through extensive one- and two-dimensional MHD tests. The numerical results demonstrate that the proposed method not only achieves the expected optimal convergence order but also exhibits robust performance in challenging MHD scenarios involving low density, low pressure, strong discontinuities, or low plasma-beta.
Our future work includes extending the PP LDF OEDG schemes to numerical simulations of MHD on unstructured meshes.

\section*{Acknowledgments}

This work was partially supported by National Natural Science Foundation of China (Grant No.~12171227) and Shenzhen Science and Technology Program (Grant No.~RCJC20221008092757098).

\bibliography{mybibfile}

\end{document}